\documentclass[12pt,usenames,dvipsnames]{article}
\usepackage{amssymb,amsmath,amsthm,tikz,multirow,xcolor,pdflscape}
\usepackage[normalem]{ulem}
\usetikzlibrary{arrows}

\title{Tilings of the Sphere by Congruent Pentagons III: Edge Combination $a^5$}

\author{
Yohji Akama, Tohoku University \\
Erxiao Wang\thanks{Corresponding author. Research was supported by ZJNU Shuang-Long Distinguished Professorship Fund No. YS304319159.}, 
Zhejiang Normal University \\
Min Yan\thanks{Research was supported by Hong Kong RGC General Research Fund 16303515 and 16305920.}, 
Hong Kong University of Science and Technology}

\usepackage[hidelinks, hyperindex]{hyperref}

\hyperref[sec:function]{}

\newcommand{\thin}{\hspace{0.1em}\rule{0.7pt}{0.8em}\hspace{0.1em}}

\DeclareMathOperator{\opt}{\; |\;}
\newcommand{\pentagon}{\tikz \foreach \a in {0,...,4} \draw[rotate=72*\a] (18:0.17) -- (90:0.17);}

\newtheorem{theorem}{Theorem}
\newtheorem{lemma}[theorem]{Lemma}

\newtheorem*{theorem*}{Theorem}

\theoremstyle{definition}
\newtheorem*{definition*}{Definition}
\newtheorem*{case*}{Case}
\newtheorem*{subcase*}{Subcase}

\theoremstyle{remark}

\numberwithin{equation}{section}

\begin{document}

\maketitle

\begin{abstract}
There are exactly eight edge-to-edge tilings of the sphere by congruent equilateral pentagons: three pentagonal subdivision tilings with 12, 24, 60 tiles; four earth map tilings with 16, 20, 24, 24 tiles; and one flip modification of the earth map tiling with 20 tiles. 

\end{abstract}

\section{Introduction}

In an edge-to-edge tiling of the sphere by congruent pentagons, the pentagon can have five possible edge combinations \cite[Lemma 9]{wy1}: $a^2b^2c,a^3bc,a^3b^2,a^4b,a^5$. We classified tilings for $a^2b^2c,a^3bc,a^3b^2$ in \cite{wy1,wy2}. In this paper, we classify tilings for $a^5$, i.e., tilings of the sphere by congruent equilateral pentagons.

\begin{theorem*}
There are eight edge-to-edge tilings of the sphere by congruent equilateral pentagons: 
\begin{enumerate}
\item Three pentagonal subdivision tilings with $12,24,60$ tiles.
\item Four earth map tilings with $16,20,24,24$ tiles.
\item One flip modification of the earth map tiling with $20$ tiles.
\end{enumerate}
\end{theorem*}

The {\em pentagonal subdivision tiling} is introduced in \cite[Section 3.1]{wy1}. We apply the operation to Platonic solids to get  three two-parameter families of tilings of the sphere by congruent pentagons with the edge combination $a^2b^2c$. The tilings in the first part of the main theorem are the equilateral reduction $a=b=c$ of the pentagonal subdivision tiling. By \cite{ay1,gsy} or direct argument, the equilateral pentagonal subdivision tiling with $f=12$ tiles is exactly the regular dodecahedron. The equilateral pentagonal subdivision tilings with $f=24,60$ tiles are derived in Section \ref{case42c}. 

\begin{figure}[htp]
\centering
\begin{tikzpicture}[>=latex,scale=1]


\foreach \a in {0,...,4}
\draw[rotate=-72*\a]
	(-54:0.33) -- (18:0.33) -- (18:0.65) -- (-18:0.8) -- (-54:0.65)
	(54:0.8) -- (54:1.2) -- (-18:1.2);


\foreach \a in {0,1,2,3}
{
\begin{scope}[xshift=3.5 cm, rotate=90*\a]

\draw
	(0,0) -- (0.4,0) -- (0.6,0.25) -- (0.25,0.6) -- (0,0.4)
	(-15:1) -- (15:1) -- (45:1) -- (75:1)
	(0.6,0.25) -- (15:1)
	(0.6,-0.25) -- (-15:1)
	(0:1.3) -- (30:1.3) -- (60:1.3) -- (90:1.3)
	(-15:1) -- (0:1.3)
	(45:1) -- (30:1.3)
	(60:1.3) -- (60:1.6);

\fill
	(0,0) circle (0.1)
	(75:1) circle (0.1)
	(60:1.6) -- ++(-75:0.1) arc (-75:-165:0.1);
	
\end{scope}
}


\foreach \a in {0,...,4}
{
\begin{scope}[xshift=8.2 cm, rotate=72*\a]

\draw
	(0,0) -- (18:0.45) -- (36:0.7) -- (72:0.7) -- (90:0.45)
	(0:0.7) -- (6:1.1)
	(30:1.1) -- (36:0.7)
	(6:1.1) -- (30:1.1) -- (54:1.1) -- (78:1.1)
	(-9:1.4) -- (6:1.1) -- (21:1.4) -- (39:1.4) -- (54:1.1)
	(-9:1.4) -- (0:1.6) -- (12:1.6) -- (21:1.4)
	(39:1.4) -- (48:1.6) -- (63:1.4)
	(12:1.6) -- (30:1.9) -- (48:1.6)
	(54:1.9) -- (30:1.9) -- (6:1.9) -- (-18:1.9) -- (0:1.6)
	(6:1.9) -- (-6:2.2)
	(30:1.9) -- (42:2.2)
	(-54:2.2) -- (-30:2.2) -- (-6:2.2) -- (18:2.2) -- (18:2.5);

\fill
	(0,0) circle (0.1)
	(6:1.1) circle (0.1)
	(30:1.9) circle (0.1)
	(18:2.5) -- ++(234:0.1) arc (234:162:0.1);

\end{scope}
}

\end{tikzpicture}
\caption{Pentagonal subdivision tilings, $\bullet=\epsilon^4$ or $\epsilon^5$.}
\label{subdivision_tiling}
\end{figure}

Figure \ref{subdivision_tiling} gives the three equilateral pentagonal subdivision tilings. The first is the regular dodecahedron, with all angles equal to $\frac{2}{3}\pi$. The second and third are the pentagonal subdivisions of the octahedron and icosahedron. The tiles are the first of Figure \ref{pentagon}. All $\epsilon$ angles appear at the vertex $\bullet$ (including the ``fractional'' $\bullet$ at the infinity in the second and third pictures). Then the $\epsilon$ angles determine the locations of the other angles according to the first of Figure \ref{pentagon}. We have $\delta=\frac{2}{3}\pi$ and $\epsilon=\frac{1}{2}\pi$ for the pentagonal subdivision of the octahedron, and the edge length $a$ is given by the biggest root of the quartic equation \eqref{edge24}. We have $\delta=\frac{2}{3}\pi$ and $\epsilon=\frac{2}{5}\pi$ for the pentagonal subdivision of the icosahedron, and the edge length $a$ is given by a quadratic factor of \eqref{edge60}.

\begin{figure}[htp]
\centering
\begin{tikzpicture}[>=latex,scale=1]


\foreach \a in {0,...,4}
\fill[gray!50, xshift=3cm, rotate=72*\a]
	(0,0) -- (18:1) -- (90:1);

\foreach \a in {0,...,4}
\foreach \b in {0,1}
{
\begin{scope}[xshift=3*\b cm, rotate=-72*\a]

\draw
	(18:1) -- (-54:1);

\node at (-18:1) {\small $a$};

\end{scope}
}

\node at (90:0.75) {$\gamma$};
\node at (162:0.75) {$\epsilon$};
\node at (15:0.75) {$\delta$};
\node at (234:0.75) {$\alpha$};
\node at (-54:0.75) {$\beta$};

\node at (0,0) {$+$};

\begin{scope}[xshift=3 cm]

\node at (90:0.75) {$\gamma$};
\node at (162:0.75) {$\delta$};
\node at (15:0.75) {$\epsilon$};
\node at (234:0.75) {$\beta$};
\node at (-54:0.75) {$\alpha$};

\node at (0,0) {$-$};

\end{scope}

\end{tikzpicture}
\caption{Equilateral pentagon.}
\label{pentagon}
\end{figure}

The angles in Figure \ref{pentagon} are arranged in a strange order (switching $\delta,\gamma$ recovers the normal order). This is the result of the sequence of argument adopted in this paper. We use this order in the introduction for the reader to easily match the values of angles and edges in the body of the paper. 

The second of Figure \ref{pentagon} is the flip of the first. We regard the arrangement of angles to be positively oriented in the first, and negatively oriented in the second. Instead of showing all the angles, it is much easier to present tilings by the locations of $\epsilon$ and the orientations of tiles. All tiles in Figure \ref{subdivision_tiling} are positively oriented. 

The {\em earth map tiling} is introduced in \cite{yan} as the edge-to-edge tilings of the sphere by pentagons, such that there are exactly two vertices of degree $>3$. This is a combinatorial concept. There are five families of such tilings distinguished by the distance between the two vertices of degree $>3$, and the four tilings in the theorem are the earth map tilings of distance $5$. Figure \ref{earth} gives the earth map tiling with $20$ tiles. All the edges pointing north converge at one vertex (north {\em pole}), and all the edges pointing south converge at another vertex (south pole). Then the left and right boundaries are glued together. The tiling consists of five countries (tiles) in the north, five countries in the south, and ten countries along the equator. 

The earth map tiling of 20 countries clearly has five {\em timezones}, between five {\em meridians} connecting the two poles. Each timezone consists of four countries (one northern, one southern, and two equatorial). If we add one more timezone, then we get the earth map tiling with $24$ tiles. If we delete one timezone, then we get the earth map tiling with $16$ tiles.

\begin{figure}[htp]
\centering
\begin{tikzpicture}[>=latex,scale=1]

\foreach \a in {0,...,4}
\fill[gray!50, xshift=1.6*\a cm]
	(-0.4,-1) -- (-0.4,-0.6) -- (0,-0.3) -- (0.8,-0.3) -- (1.2,-0.6) -- (1.2,-1)
	(0.8,-0.3) -- (0.8,0.3) -- (1.6,0.3) -- (1.6,-0.3) -- (1.2,-0.6);	

\foreach \a in {0,...,5}
\draw[xshift=1.6*\a cm]
	(0.4,1) -- (0.4,0.6) -- (0,0.3) -- (0,-0.3) -- (-0.4,-0.6) -- (-0.4,-1);

\foreach \a in {0,...,4}
{
\begin{scope}[xshift=1.6*\a cm]

\draw
	(0,-0.3) -- (0.8,-0.3) -- (1.2,-0.6)
	(0.4,0.6) -- (0.8,0.3) -- (1.6,0.3)
	(0.8,0.3) -- (0.8,-0.3);
	
\node at (0.85,0.45) {\small $\epsilon$};	
\node at (0.65,0.25) {\small $\epsilon$};
\node at (1.45,-0.25) {\small $\epsilon$};
\node at (0,-0.45) {\small $\epsilon$};	

\end{scope}
}

\foreach \a in {0,1}
\draw[very thick, xshift=8*\a cm]
	(0.4,1) -- (0.4,0.6) -- (0,0.3) -- (0,-0.3) -- (-0.4,-0.6) -- (-0.4,-1);

\draw[very thick]
	(4.4,-1) -- (4.4,-0.6) -- (4,-0.3) -- (4,0.3) -- (3.6,0.6) -- (3.6,1);
	
\node at (4,1.3) {\small north};
\node at (4,-1.3) {\small south};

\draw[->, very thick]
	(-0.7,0) -- ++(0.5,0);

\node at (-1.4,0) {\small equator};

\end{tikzpicture}
\caption{Earth map tiling.}
\label{earth}
\end{figure}

For the earth map tilings in the second part of the main theorem, we need all the edges to be great arcs, and we need to provide all the angles in Figure \ref{earth}. The angle arrangements are given by $\epsilon$, and the orientations of tiles according to Figure \ref{pentagon}. The pentagons are unique, and the angles have specific values.

The two earth map tilings with $24$ tiles are derived in Sections \ref{case55}, \ref{case12e}, \ref{casespecial}. The first has $\beta=\frac{4}{3}\pi,\delta=\frac{1}{3}\pi,\epsilon=\frac{5}{6}\pi$, and the edge length is given by \eqref{earth24A}. The second has $\gamma=\frac{1}{2}\pi,\delta=\frac{1}{3}\pi,\epsilon=\frac{5}{6}\pi$, and the edge length is given by \eqref{earth24B}. 

The earth map tiling with $16$ tiles is derived in Section \ref{case15b}. The pentagon has $\delta=\frac{1}{2}\pi,\epsilon=\frac{3}{4}\pi$, and the edge length is given by the biggest root of \eqref{length16}.

The earth map tiling with $20$ tiles is derived in Section \ref{case14e}. The pentagon has $\delta=\frac{2}{5}\pi,\epsilon=\frac{4}{5}\pi$, and the edge length is given by the biggest root of \eqref{length20}.

Next, we explain the {\em flip modification} in the third part of the main theorem. In Figure \ref{earth}, we may use the thick lines to divide the earth map tiling with $20$ tiles into two halves (eastern and western hemispheres). The left of Figure \ref{flipmod} shows the boundary between these two halves. Due to $2\delta=\epsilon$, we may flip the interior half with respect to the line $L$, and still keep the angle sums of all the vertices along the boundary to be $2\pi$. Therefore the flip modification is still a tiling.

\begin{figure}[htp]
\centering
\begin{tikzpicture}[>=latex]

\draw[gray]
	(126:2.2) -- (-54:2.2);

\node at (1.5,-1.7) {\small $L$};

\foreach \b in {0,1}
{
\begin{scope}[xshift=5*\b cm]

\foreach \a in {0,...,9}
\draw[rotate=36*\a]
	(90:1.5) -- (54:1.5);

\draw
	(54:1.5) -- (54:1.9)
	(-18:1.5) -- (-18:1.9)
	(126:1.5) -- (126:1.9)
	(198:1.5) -- (198:1.9);

\node at (90:1.7) {\small $\delta^2$};
\node at (60:1.65) {\small $\gamma$};
\node at (48:1.65) {\small $\alpha$};
\node at (18:1.7) {\small $\beta$};
\node at (-11:1.65) {\small $\delta$};
\node at (-24:1.65) {\small $\epsilon$};
\node at (-54:1.7) {\small $\gamma$};

\node at (118:1.65) {\small $\beta$};
\node at (132:1.65) {\small $\alpha$};
\node at (162:1.65) {\small $\epsilon$};
\node at (192:1.65) {\small $\gamma$};
\node at (204:1.65) {\small $\alpha$};
\node at (234:1.7) {\small $\beta$};
\node at (-90:1.7) {\small $\delta^3$};

\end{scope}
}

\draw[->, very thick]
	(2,0) -- ++(1,0);
	
\node at (2.5,0.3) {flip};

	
\draw
	(18:1.5) -- (18:1)
	(-54:1.5) -- (-54:1)
	(162:1.5) -- (162:1)
	(234:1.5) -- (234:1);
	
\node at (90:1.25) {\small $\delta^3$};
\node at (54:1.25) {\small $\beta$};
\node at (26:1.25) {\small $\alpha$};
\node at (10:1.25) {\small $\gamma$};
\node at (-18:1.3) {\small $\epsilon$};
\node at (-46:1.25) {\small $\alpha$};
\node at (-63:1.25) {\small $\beta$};

\node at (126:1.25) {\small $\gamma$};
\node at (154:1.3) {\small $\epsilon$};
\node at (170:1.25) {\small $\delta$};
\node at (198:1.25) {\small $\beta$};
\node at (226:1.25) {\small $\alpha$};
\node at (242:1.25) {\small $\gamma$};
\node at (-90:1.25) {\small $\delta^2$};


\begin{scope}[xshift=5 cm]

\draw
	(18:1.5) -- (18:1)
	(-54:1.5) -- (-54:1)
	(90:1.5) -- (90:1)
	(234:1.5) -- (234:1);
	
\node at (126:1.25) {\small $\gamma$};
\node at (98:1.3) {\small $\epsilon$};
\node at (162:1.25) {\small $\delta^3$};
\node at (198:1.25) {\small $\beta$};
\node at (226:1.25) {\small $\alpha$};
\node at (242:1.25) {\small $\gamma$};
\node at (-90:1.3) {\small $\epsilon$};

\node at (82:1.25) {\small $\delta$};
\node at (54:1.25) {\small $\beta$};
\node at (26:1.25) {\small $\alpha$};
\node at (10:1.25) {\small $\gamma$};
\node at (-18:1.25) {\small $\delta^2$};
\node at (-62:1.25) {\small $\alpha$};
\node at (-44:1.25) {\small $\beta$};

\end{scope}

\end{tikzpicture}
\caption{Flip modification of earth map tiling.}
\label{flipmod}
\end{figure}

The earth map tilings in the main theorem actually fit into a family of tilings of the sphere by congruent almost equilateral (i.e., edge combination $a^4b$) pentagons \cite{ly1}. The number $f$ of tiles can be any multiple of $4$, and for each such $f$, the tiling allows one free parameter. Moreover, if $\frac{f}{4}$ is odd, then we may apply the flip modification to get another tiling. 

For exactly four cases, the almost equilateral earth map tiling becomes equilateral. Moreover, for the single case of $f=20$, we may apply the flip modification. 

Next, we outline the technique for proving the main theorem. Due to lack of edge length variation, we cannot use much of the technique in the other two papers in the series \cite{wy1,wy2}. A completely new idea is needed: A general pentagon is determined by the free choices of $4$ edge lengths and $3$ angles, yielding $7$ degrees of freedom. The requirement that all $5$ edges are equal imposes $4$ equations, leaving $7-4=3$ degrees of freedom for equilateral pentagons. Therefore $3$ more independent equations are enough to completely determine equilateral pentagons. Once we know the pentagons, it is then easy to derive tilings.

In Section \ref{deg3}, for general tilings of any surface with up to five distinct angles, we derive all the possible angle combinations at degree $3$ vertices. The result is summarised in Lemma \ref{deg3v} and Table \ref{deg3AVC}. In Section \ref{anglecombo}, we specialise Lemma \ref{deg3v} to tilings of the sphere by congruent equilateral pentagons and find that all except one case has three degree $3$ vertices. The angle sums of the three vertices provide three independent equations that we can use to determine the equilateral pentagon. The remaining case has two degree $3$ vertices, and we consider the possibility of adding another vertex of degree $4$ or $5$. With the additional vertex, we again get three equations that we can use to determine the equilateral pentagon. After all these cases, we have one remaining exceptional case of two specific degree $3$ vertices, and no vertices of degree $4$ or $5$ (i.e., any vertex of degree $>3$ has degree $\ge 6$). 

In Section \ref{calculation}, we calculate equilateral pentagons for all the cases in Section \ref{anglecombo}. If we take into account of various angle arrangements in the pentagon, altogether we have more than 400 cases to calculate. The calculation is based on symbolic manipulations of polynomials, and can be numerically arbitrarily accurate. Therefore all our calculations for equilateral pentagons suitable for tilings are absolutely rigorous. At the end, we get $24$ cases (some sharing the same pentagon) yielding equilateral pentagons that might be suitable for tiling the sphere. In Section \ref{tiling}, we construct tilings for these cases, plus the exceptional case.  

Throughout the paper, we will assume $f$ is an even number $\ge 16$. As explained in \cite[Section 2.1]{wy1}, the assumption follows from \cite{ay1,gsy,yan}. 

The values of angles and edge lengths are expressed as multiples of $\pi$. When the multiple is a rational fraction, such as $\alpha=\frac{2}{3}\pi$, we mean the precise value. When the multiple is a decimal expression, such as $\alpha=0.29539\pi$, we mean an approximate value. In other words, $0.29539\pi\le\alpha < 0.29540\pi$.

Finally, we would like to thank Hoi Ping Luk for the preliminary work \cite{luk} that leads to Lemma 1.

\section{Angle Combination at Degree $3$ Vertex}
\label{deg3}

We study the collection $\text{AVC}_3$ (for {\em anglewise vertex combination}) of angle combinations at degree $3$ vertices in an edge-to-edge tiling of any surface, under the assumption of up to five distinct angles appearing at degree $3$ vertices. We distinguish angles by their values. 

Each angle combination at a vertex gives an angle sum equality. If there are too many angle combinations, then the corresponding equalities will force some distinct angles to have the same value. Our goal is to list all the possible collections. Since only the angle sums of vertices are used, our list actually applies to any edge-to-edge tiling of any surface. 

We say an angle combination at a vertex is of {\em $111$-type} if it is $\lambda\mu\nu$ for three distinct angles $\lambda,\mu,\nu$. Moreover, an {\em $12$-type} vertex is $\lambda\mu^2$, and a {\em $3$-type} vertex is $\lambda^3$. 

\subsubsection*{Case. One angle $\alpha$ at degree $3$ vertices}

The only possible degree $3$ vertex is $\alpha^3$. We denote this by 
\[
\text{AVC}_3=\{\alpha^3\}.
\]

\subsubsection*{Case. Two distinct angles $\alpha,\beta$ at degree $3$ vertices}

There are four possible degree $3$ vertices $\alpha^3,\alpha^2\beta,\alpha\beta^2,\beta^3$. If any two appear in a tiling, then their angle sums imply $\alpha=\beta$. For example, if both $\alpha^2\beta$ and $\alpha\beta^2$ are vertices, then we have $2\alpha+\beta=2\pi=\alpha+2\beta$, which implies $\alpha=\beta$. Therefore only one of the four can appear. In order for both $\alpha$ and $\beta$ to appear, the vertex is either $\alpha^2\beta$ or $\alpha\beta^2$. Up to the symmetry $\alpha\leftrightarrow \beta$, we get
\[
\text{AVC}_3=\{\alpha\beta^2\}.
\]

\subsubsection*{Case. Three distinct angles $\alpha,\beta,\gamma$ at degree $3$ vertices}

If $\alpha\beta\gamma$ is a vertex, then all three angles already appear, and we get $\{\alpha\beta\gamma\}\subset \text{AVC}_3$. Then we find whether the $\text{AVC}_3$ can admit additional degree $3$ vertices. An additional vertex of $12$-type always forces some angles to become equal. For example, the angle sums of $\alpha\beta\gamma$ and $\alpha\gamma^2$ imply $\beta=\gamma$. On the other hand, it is possible for one of $\alpha^3,\beta^3,\gamma^3$ to be a vertex while still keeping $\alpha,\beta,\gamma$ distinct. Moreover, if two of $\alpha^3,\beta^3,\gamma^3$ are vertices, then the corresponding angles become equal. Up to the symmetry of permuting $\alpha,\beta,\gamma$, therefore, we get 
\[
\{\alpha\beta\gamma\}
\subset \text{AVC}_3 \subset
\{\alpha\beta\gamma,\alpha^3\}.
\]
We denote this by writing 
\[
\text{AVC}_3=\{\alpha\beta\gamma \opt \alpha^3\}.
\]
Here $\alpha\beta\gamma$ is the {\em necessary} part to make sure all three angles appear at degree $3$ vertices, and $\alpha^3$ is the {\em optional} part that can be added without forcing distinct angles to become equal. The necessary part is the lower bound for the $\text{AVC}_3$ and the necessary plus optional part is the upper bound for the $\text{AVC}_3$.

Next we assume that there are no $111$-type vertices, and there are $12$-type vertices. Up to symmetry, we may assume  $\alpha\beta^2$ is a vertex. Then $\gamma$ may appear as $\alpha^2\gamma,\beta\gamma^2,\gamma^3$, without forcing some angles to become equal. Since $\{\alpha\beta^2,\alpha^2\gamma\}$ can be transformed to $\{\alpha\beta^2,\beta\gamma^2\}$ via $\alpha\to\beta\to\gamma\to\alpha$, up to symmetry, we get two possible necessary parts 
\[
\{\alpha\beta^2,\alpha^2\gamma\},\quad
\{\alpha\beta^2,\gamma^3\}.
\]
It can be easily verified that neither allow optional vertices.

Finally we assume that there are no $111$-type and $12$-type vertices. This means that only $\alpha^3,\beta^3,\gamma^3$ can appear. Since we cannot have all three angles to appear in this way without forcing them to become equal, there is no more $\text{AVC}_3$.

\subsubsection*{Case. Four distinct angles $\alpha,\beta,\gamma,\delta$ at degree $3$ vertices}

If $\alpha\beta\gamma$ is a vertex, then up to symmetry, the angle $\delta$ may appear as $\alpha\delta^2,\alpha^2\delta,\delta^3$. This gives three possible necessary parts 
\[
\{\alpha\beta\gamma,\alpha\delta^2\},\quad
\{\alpha\beta\gamma,\alpha^2\delta\},\quad
\{\alpha\beta\gamma,\delta^3\}.
\]
By $\alpha\beta\gamma$, the necessary parts do not allow  optional $111$-type vertices. The necessary part $\{\alpha\beta\gamma,\alpha\delta^2\}$ only allows any one of $\beta^2\delta,\gamma^2\delta,\beta^3,\gamma^3$ to be optional, and does not allow any two (i.e., the appearance of two from the four will imply distinct angles becoming equal). The necessary part $\{\alpha\beta\gamma,\alpha^2\delta\}$ only allows any one of $\beta\delta^2,\gamma\delta^2,\beta^3,\gamma^3$ to be optional, and does not allow any two. The necessary part $\{\alpha\beta\gamma,\delta^3\}$ does not allow optional vertices. Therefore up to symmetry, we get five possible $\text{AVC}_3$s
\[
\{\alpha\beta\gamma,\alpha\delta^2\opt \beta^2\delta\},\quad
\{\alpha\beta\gamma,\alpha\delta^2\opt \beta^3\},
\]
\[
\{\alpha\beta\gamma,\alpha^2\delta\opt \beta\delta^2\},\quad
\{\alpha\beta\gamma,\alpha^2\delta\opt \beta^3\},\quad
\{\alpha\beta\gamma,\delta^3\}.
\]

Next we assume that there are no $111$-type vertices, and $\alpha\beta^2$ is a vertex. Then $\gamma$ may appear as $\alpha^2\gamma,\beta\gamma^2,\gamma\delta^2,\gamma^2\delta,\gamma^3$. Up to symmetry, we may drop $\beta\gamma^2$ and $\gamma^2\delta$. Moreover, for the combinations $\{\alpha\beta^2,\alpha^2\gamma\}$ and $\{\alpha\beta^2,\gamma^3\}$, we need to further consider the way $\delta$ appears. Up to symmetry, this leads to six possible necessary parts
\[
\{\alpha\beta^2,\alpha^2\gamma,\beta\delta^2\},\quad
\{\alpha\beta^2,\alpha^2\gamma,\gamma^2\delta\},\quad
\{\alpha\beta^2,\alpha^2\gamma,\delta^3\},
\]
\[
\{\alpha\beta^2,\gamma\delta^2\},\quad
\{\alpha\beta^2,\gamma^3,\alpha^2\delta\},\quad
\{\alpha\beta^2,\gamma^3,\beta\delta^2\}.
\] 
The first is a special case of the fourth via $\alpha\to\beta\to\gamma\to\alpha$. The second is a special case of the fourth via $\gamma\leftrightarrow \delta$. The fifth becomes the third via $\gamma\leftrightarrow \delta$. The sixth becomes the third via $\alpha\to\gamma\to\delta\to\beta\to\alpha$ Therefore we only need to consider the third and fourth necessary parts. 

The third does not allow optional vertices. Under the assumption of no $111$-type vertices, the fourth only allows $\alpha^2\delta,\beta\gamma^2$ to be optional, and does not allow two. Up to symmetry, we get two possible $\text{AVC}_3$s
\[
\{\alpha\beta^2,\alpha^2\gamma,\delta^3\},\quad
\{\alpha\beta^2,\gamma\delta^2 \opt \alpha^2\delta\}.
\]

Finally, it is easy to see that we cannot have all vertices to be of $3$-type.

\subsubsection*{Case. Five distinct angles $\alpha,\beta,\gamma,\delta,\epsilon$ at degree $3$ vertices}

If all vertices are of $3$-type, then all angles are equal. Therefore either there are $111$-type vertices, or there are $12$-type vertices. Moreover, given five distinct angles, there can be at most two $111$-type vertices. This leads to three subcases.

\subsubsection*{Subcase. There are two $111$-type vertices}

Up to symmetry, we may assume that $\alpha\beta\gamma$ and $\alpha\delta\epsilon$ are all the $111$-type vertices. We get one necessary part $\{\alpha\beta\gamma,\alpha\delta\epsilon\}$, and the optional vertices are of either $12$-type or $3$-type. It is easy to see that the only possible optional vertex involving $\alpha$ is $\alpha^3$. 

We look for optional $12$-type vertices. Since such a vertex cannot involve $\alpha$, up to symmetry that preserves the collection $\{\alpha\beta\gamma,\alpha\delta\epsilon\}$, such an optional vertex is $\beta\delta^2$. Then $\{\alpha\beta\gamma,\alpha\delta\epsilon\opt\beta\delta^2\}$  further allows any one of the 12-type vertices $\beta^2\epsilon,\gamma\epsilon^2,\gamma^2\delta,\gamma^2\epsilon$ to be optional, but not two. Therefore we get $\{\alpha\beta\gamma,\alpha\delta\epsilon \opt \beta\delta^2,\beta^2\epsilon\}$, $\{\alpha\beta\gamma,\alpha\delta\epsilon \opt \beta\delta^2,\gamma\epsilon^2\}$, $\{\alpha\beta\gamma,\alpha\delta\epsilon \opt \beta\delta^2,\gamma^2\delta\}$, $\{\alpha\beta\gamma,\alpha\delta\epsilon \opt \beta\delta^2,\gamma^2\epsilon\}$.  Then we ask whether any of the four $\text{AVC}_3$s allows 3-type optional vertices. The answer is that only the second one allows $\alpha^3$, and $\{\alpha\beta\gamma,\alpha\delta\epsilon \opt \beta\delta^2,\gamma\epsilon^2,\alpha^3\}$ does not allow any more optional vertices. Since $\{\alpha\beta\gamma,\alpha\delta\epsilon \opt \beta\delta^2,\gamma^2\delta\}$ is changed to $\{\alpha\beta\gamma,\alpha\delta\epsilon \opt \beta\delta^2,\beta^2\epsilon\}$ via $\delta\to \beta\to \epsilon\to\gamma\to\delta$, we get three possible $\text{AVC}_3$s 
\[
\{\alpha\beta\gamma,\alpha\delta\epsilon \opt \beta\delta^2,\beta^2\epsilon\},\quad
\{\alpha\beta\gamma,\alpha\delta\epsilon \opt \beta\delta^2,\gamma\epsilon^2,\alpha^3\},\quad
\{\alpha\beta\gamma,\alpha\delta\epsilon \opt \beta\delta^2,\gamma^2\epsilon\}.
\]

After exhausting further optional $12$-type vertices for $\{\alpha\beta\gamma,\alpha\delta\epsilon \opt \beta\delta^2\}$, we still need to look for further optional $3$-type vertices, under the assumption of no further optional $12$-type vertices. The answer is $\alpha^3,\gamma^3,\epsilon^3$. Since $\alpha^3$ is already included in the second $\text{AVC}_3$ above, we get two more possible $\text{AVC}_3$s 
\[
\{\alpha\beta\gamma,\alpha\delta\epsilon \opt \beta\delta^2,\gamma^3\},\quad
\{\alpha\beta\gamma,\alpha\delta\epsilon \opt \beta\delta^2,\epsilon^3\}.
\]

It remains to consider the case there are no $12$-type optional vertices. In other words, all optional vertices are of $3$-type. Up to symmetry, we get $\{\alpha\beta\gamma,\alpha\delta\epsilon \opt\alpha^3\}$, $\{\alpha\beta\gamma,\alpha\delta\epsilon \opt\beta^3\}$. The first is included in one of the five $\text{AVC}_3$s above, and the second is also included via $\beta\leftrightarrow \gamma$.

\subsubsection*{Subcase. There is only one $111$-type vertex}

Up to symmetry, we may assume that $\alpha\beta\gamma$ is the only $111$-type vertex.

Since $\delta^3$ and $\epsilon^3$ cannot both be vertices, one of $\delta,\epsilon$ must appear in an $12$-type vertex. Up to symmetry, we may assume that one of $\alpha\delta^2,\alpha^2\delta,\delta\epsilon^2$ is a vertex.

If $\alpha\beta\gamma$ and $\alpha\delta^2$ are vertices, then under the assumption of no more $111$-type vertices, and up to the symmetry $\beta\leftrightarrow\gamma$, the angle $\epsilon$ may appear as $\alpha^2\epsilon,\beta\epsilon^2,\beta^2\epsilon,\delta\epsilon^2,\epsilon^3$. Similarly, if $\alpha\beta\gamma$ and $\alpha^2\delta$ are vertices, then $\epsilon$ may appear as $\alpha\epsilon^2,\beta\epsilon^2,\beta^2\epsilon,\delta^2\epsilon,\epsilon^3$. We get total of ten combinations, such that all five angles appear. Up to symmetry, the ten combinations are reduced to eight possible necessary parts. We also list all the allowed optional vertices under the assumption of no more $111$-type vertices:
\begin{align*}
\{\alpha\beta\gamma,\alpha\delta^2,\alpha^2\epsilon\}
&\colon
\beta\epsilon^2,\beta^2\delta,\text{\sout{\ensuremath{\gamma\epsilon^2}}},\text{\sout{\ensuremath{\gamma^2\delta}}},\beta^3,\text{\sout{\ensuremath{\gamma^3}}} 
\\
\{\alpha\beta\gamma,\alpha\delta^2,\delta\epsilon^2\}
&\colon
\beta^2\epsilon,\text{\sout{\ensuremath{\gamma^2\epsilon}}},\beta^3,\text{\sout{\ensuremath{\gamma^3}}} 
\\
\{\alpha\beta\gamma,\alpha^2\delta,\delta^2\epsilon\}
&\colon
\beta\epsilon^2,\text{\sout{\ensuremath{\gamma\epsilon^2}}},\beta^3,\text{\sout{\ensuremath{\gamma^3}}} 
\\
\{\alpha\beta\gamma,\alpha\delta^2,\beta\epsilon^2\}
&\colon 
\alpha^2\epsilon,\text{\sout{\ensuremath{\beta^2\delta}}},\gamma^2\delta,\text{\sout{\ensuremath{\gamma^2\epsilon}}},\gamma^3 
\\
\{\alpha\beta\gamma,\alpha\delta^2,\beta^2\epsilon\}
&\colon
\gamma\epsilon^2,\gamma^2\delta,\gamma^3,\delta\epsilon^2 
\\
\{\alpha\beta\gamma,\alpha^2\delta,\beta^2\epsilon\}
&\colon
\alpha\epsilon^2,\text{\sout{\ensuremath{\beta\delta^2}}},\gamma\delta^2,\text{\sout{\ensuremath{\gamma\epsilon^2}}},\gamma^3 
\\
\{\alpha\beta\gamma,\alpha\delta^2,\epsilon^3\}
&\colon
\beta^2\delta,\text{\sout{\ensuremath{\gamma^2\delta}}} 
\\
\{\alpha\beta\gamma,\alpha^2\delta,\epsilon^3\}
&\colon 
\beta\delta^2,\text{\sout{\ensuremath{\gamma\delta^2}}} 
\end{align*}
The necessary part $\{\alpha\beta\gamma,\alpha\delta^2,\alpha^2\epsilon\}$ allows any one of the six optional vertices, but not two. Up to the symmetry $\beta\leftrightarrow \gamma$ of the necessary part, the list may be reduced to $\beta\epsilon^2,\beta^2\delta,\beta^3$ (indicated by $\text{\sout{\ensuremath{\gamma\epsilon^2}}},\text{\sout{\ensuremath{\gamma^2\delta}}},\text{\sout{\ensuremath{\gamma^3}}}$). Therefore we get three possible $\text{AVC}_3$s 
\[
\{\alpha\beta\gamma,\alpha\delta^2,\alpha^2\epsilon \opt \beta\epsilon^2\},\quad
\{\alpha\beta\gamma,\alpha\delta^2,\alpha^2\epsilon \opt \beta^2\delta\},\quad
\{\alpha\beta\gamma,\alpha\delta^2,\alpha^2\epsilon \opt \beta^3\}.
\]
Similar discussion applies to all the other necessary parts. In particular, we cannot have two optional vertices. 
Up to symmetry, we get total of nineteen $\text{AVC}_3$s, each with one vertex that is not stricken out as the optional vertex.

It remains to consider the case $\delta\epsilon^2$ is a vertex. After considering the vertices $\alpha\delta^2$ or $\alpha^2\delta$, we may additionally assume that $\alpha,\beta,\gamma$ is not combined with $\delta,\epsilon$ to form a degree $3$ vertex. Therefore $\delta\epsilon^2$ is the only degree $3$ vertex involving $\delta,\epsilon$, and the optional vertices allowed by the necessary part $\{\alpha\beta\gamma,\delta\epsilon^2\}$ can only involve $\alpha,\beta,\gamma$. By the discussion about the $\text{AVC}_3$ of three distinct angles with the necessary part $\{\alpha\beta\gamma\}$, we get one possible $\text{AVC}_3$
\[
\{\alpha\beta\gamma,\delta\epsilon^2 \opt \alpha^3\}.
\]

\subsubsection*{Subcase. There are no $111$-type vertices}

Since there can be at most one $3$-type vertex, we may assume $\alpha^3,\beta^3,\gamma^3,\delta^3$ are not vertices. This means $\alpha,\beta,\gamma,\delta$ all appear at $12$-type vertices. This implies that, under the assumption of no $111$-type vertices, and up to the symmetry of exchanging $\alpha,\beta,\gamma,\delta$, we may assume $\alpha\beta^2$ and $\gamma\delta^2$ are vertices. Moreover, up to the symmetry $\alpha\leftrightarrow\gamma$ and $\beta\leftrightarrow\delta$, we may further assume $\epsilon$ appears as $\alpha^2\epsilon,\beta\epsilon^2,\epsilon^3$. This gives three possible necessary parts
\[
\{\alpha\beta^2,\gamma\delta^2,\alpha^2\epsilon\},\quad
\{\alpha\beta^2,\gamma\delta^2,\beta\epsilon^2\},\quad
\{\alpha\beta^2,\gamma\delta^2,\epsilon^3\}.
\]
Since the second becomes the first via $\alpha\to\epsilon\to\beta\to\alpha$, we only need to consider the first and the third.

By no $111$-type vertices, the necessary part $\{\alpha\beta^2,\gamma\delta^2,\alpha^2\epsilon\}$ only allows $\beta\gamma^2,\delta\epsilon^2$ to be optional, and does not allow two. Then we get two possible $\text{AVC}_3$s
\[
\{\alpha\beta^2,\gamma\delta^2,\alpha^2\epsilon \opt \beta\gamma^2\},\quad
\{\alpha\beta^2,\gamma\delta^2,\alpha^2\epsilon \opt \delta\epsilon^2\}.
\]
The necessary part $\{\alpha\beta^2,\gamma\delta^2,\epsilon^3\}$ only allows $\alpha^2\delta$, $\beta\gamma^2$ to be optional, and does not allow two. Up to symmetry, we get one possible $\text{AVC}_3$
\[
\{\alpha\beta^2,\gamma\delta^2,\epsilon^3 \opt \alpha^2\delta\}.
\]

The following is the summary of our discussion.

\begin{lemma}\label{deg3v}
If an edge-to-edge tiling of a surface has at most five distinct angles at degree $3$ vertices, then after suitable relabeling of the distinct angles, the anglewise vertex combinations at degree $3$ vertices is in Table \ref{deg3AVC}.
\end{lemma}

\begin{table}[h]
\centering
\begin{tabular}{|c|c|c|c|}
\hline 
\multicolumn{3}{|c|}{Necessary} & Optional  \\
\hline 
\hline  
$\alpha^3$ & \multicolumn{1}{r}{} & & \\
\hline 
\hline  
$\alpha\beta^2$ & \multicolumn{1}{r}{} & & \\
\hline 
\hline  
\multicolumn{2}{|c|}{$\alpha\beta\gamma$} & & $\alpha^3$ \\
\cline{1-2} \cline{4-4} 
\multirow{2}{*}{$\alpha\beta^2$} & $\alpha^2\gamma$ & & \\
\cline{2-2}  
 & $\gamma^3$ & & \\
\hline 
\hline  
\multicolumn{2}{|c|}{\multirow{5}{*}{$\alpha\beta\gamma$}}& \multirow{2}{*}{$\alpha\delta^2$}  &  $\beta^2\delta$ \\
\cline{4-4}
\multicolumn{2}{|c|}{} &&  $\beta^3$  \\
\cline{3-4} 
\multicolumn{2}{|c|}{} & \multirow{2}{*}{$\alpha^2\delta$}  & $\beta\delta^2$ \\
\cline{4-4}
\multicolumn{2}{|c|}{} &&  $\beta^3$ \\
\cline{3-4} 
\multicolumn{2}{|c|}{} & $\delta^3$ &  \\
\hline 
\multirow{2}{*}{$\alpha\beta^2$} & \multicolumn{2}{|c|}{$\gamma\delta^2$}  & $\alpha^2\delta$ \\
\cline{2-4} 
 & \multicolumn{2}{|c|}{$\alpha^2\gamma,\delta^3$} &  \\
\hline
\hline 
\multirow{5}{*}{$\alpha\beta\gamma$} & 
\multicolumn{2}{|c|}{\multirow{5}{*}{$\alpha\delta\epsilon$}}
& $\beta\delta^2,\beta^2\epsilon$ \\
\cline{4-4}
& \multicolumn{2}{|c|}{} &  $\beta\delta^2,\gamma\epsilon^2,\alpha^3$ \\
\cline{4-4}
& \multicolumn{2}{|c|}{} &  $\beta\delta^2,\gamma^2\epsilon$ \\
\cline{4-4}
& \multicolumn{2}{|c|}{} &  $\beta\delta^2,\gamma^3$ \\
\cline{4-4}
& \multicolumn{2}{|c|}{} &  $\beta\delta^2,\epsilon^3$ \\
\hline
\end{tabular}
\quad
\begin{tabular}{|c|c|c|c|}
\hline 
\multicolumn{3}{|c|}{Necessary} & Optional  \\
\hline 
\hline 
\multirow{19}{*}{$\alpha\beta\gamma$} 
& \multirow{12}{*}{$\alpha\delta^2$} 
& \multirow{3}{*}{$\alpha^2\epsilon$} 
& $\beta\epsilon^2$ \\
\cline{4-4}
&&& $\beta^2\delta$ \\
\cline{4-4}
&&& $\beta^3$ \\ 
\cline{3-4} 
&& \multirow{3}{*}{$\beta\epsilon^2$} 
& $\alpha^2\epsilon$ \\
\cline{4-4}
&&& $\gamma^2\delta$ \\
\cline{4-4}
&&& $\gamma^3$ \\
\cline{3-4} 
&& \multirow{4}{*}{$\beta^2\epsilon$}  
& $\gamma\epsilon^2$ \\
\cline{4-4}
&&& $\gamma^2\delta$ \\
\cline{4-4}
&&& $\gamma^3$ \\
\cline{4-4}
&&& $\delta\epsilon^2$ \\
\cline{3-4} 
&& \multirow{2}{*}{$\delta\epsilon^2$} 
& $\beta^2\epsilon$ \\
\cline{4-4}
&&& $\beta^3$ \\
\cline{3-4} 
&& $\epsilon^3$ & $\beta^2\delta$ \\
\cline{2-4} 
& \multirow{6}{*}{$\alpha^2\delta$}
& \multirow{3}{*}{$\beta^2\epsilon$} 
& $\alpha\epsilon^2$ \\
\cline{4-4}
&&& $\gamma\delta^2$ \\
\cline{4-4}
&&& $\gamma^3$ \\
\cline{3-4} 
&& \multirow{2}{*}{$\delta^2\epsilon$} 
& $\beta\epsilon^2$ \\
\cline{4-4}
&&& $\beta^3$ \\
\cline{3-4} 
&& $\epsilon^3$ & $\beta\delta^2$ \\
\cline{2-4} 
& \multicolumn{2}{|c|}{$\delta\epsilon^2$} &   $\alpha^3$ \\
\hline 
\multicolumn{2}{|c|}{\multirow{3}{*}{$\alpha\beta^2,\gamma\delta^2$}}
& \multirow{2}{*}{$\alpha^2\epsilon$} &  $\beta\gamma^2$ \\
\cline{4-4} 
\multicolumn{2}{|c|}{} &&  $\delta\epsilon^2$ \\
\cline{3-4} 
\multicolumn{2}{|c|}{} & $\epsilon^3$ &  $\alpha^2\delta$ \\
\hline 
\end{tabular}
\caption{Anglewise vertex combinations at degree $3$ vertices.}
\label{deg3AVC}
\end{table}

We note that the argument leading to the table follows a specific sequence of cases, and the discussion of a case often assumes the exclusion of the earlier cases. Following a different sequence of cases may lead to a different table, and not excluding earlier cases may introduce many overlappings between various cases.

\section{Angle in Pentagonal Tiling of the Sphere}
\label{anglecombo}

Lemma \ref{deg3v} gives all the collections of angle combinations at degree $3$ vertices in a general tiling. By taking the sphere and the pentagon into consideration, we get more information. We use the following basic facts about pentagonal tilings of the sphere, from Lemma 4 through Lemma 8 of \cite{wy1}.

\begin{lemma}\label{anglesum}
If all tiles in a tiling of the sphere by $f$ pentagons have the same five angles $\alpha,\beta,\gamma,\delta,\epsilon$, then 
\[
\alpha+\beta+\gamma+\delta+\epsilon
=(3 + \tfrac{4}{f})\pi.
\]
\end{lemma}

\begin{lemma}\label{deg3a}
If an angle appears at every degree $3$ vertex in a tiling of the sphere by pentagons with the same angle combination, then the angle appears at least twice in the pentagon.
\end{lemma}

\begin{lemma}\label{deg3b}
If an angle appears at least twice at every degree $3$ vertex in a tiling of the sphere by pentagons with the same angle combination, then the angle appears at least three times in the pentagon. 
\end{lemma}

\begin{lemma}\label{deg3c}
If two angles together appear at least twice at every degree $3$ vertex in a tiling of the sphere by pentagons with the same angle combination, then the two angles together appear at least three times in the pentagon. 
\end{lemma}

\begin{lemma}\label{hdeg}
Suppose an angle $\theta$ does not appear at degree $3$ vertices in a tiling of the sphere by pentagons with the same angle combination.
\begin{enumerate}
\item There can be at most one such angle $\theta$.
\item The angle $\theta$ appears only once in the pentagon.
\item $2v_4+v_5\ge 12$.
\item One of $\theta^3\rho,\theta^4,\theta^5$ is a vertex, where $\rho\ne\theta$.
\end{enumerate}
\end{lemma}

The lemmas so far do not use any edge information. In fact, they do not even require edges to be great arcs. If we require all tiles to be congruent, then by \cite[Lemma 1]{gsy}, the boundary of any tile is a simple closed curve. For spherical pentagon (edges are great arcs) with simple boundary, we have the following constraint from \cite[Lemma 1]{wy2}.

\begin{lemma}\label{geometry1}
If the spherical pentagon in Figure \ref{geom1} is simple and has two pairs of edges of equal lengths $a$ and $b$, then $\beta>\gamma$ is equivalent to $\delta<\epsilon$.
\end{lemma}

\begin{figure}[htp]
\centering
\begin{tikzpicture}

\draw
	(0,0.7) -- node[fill=white,inner sep=0.5] {\small $b$}
	(-1,0) -- node[fill=white,inner sep=0.5] {\small $a$}
	(-0.7,-1) -- 
	(0.7,-1) -- node[fill=white,inner sep=0.5] {\small $a$}
	(1,0) -- node[fill=white,inner sep=0.5] {\small $b$}
	cycle;

\node at (0,0.5) {\small $\alpha$};	
\node at (-0.75,-0.05) {\small $\beta$};
\node at (0.8,-0.1) {\small $\gamma$};	
\node at (-0.6,-0.8) {\small $\delta$};	
\node at (0.6,-0.8) {\small $\epsilon$};	

\end{tikzpicture}
\caption{Pentagon with two pairs of equal edges.}
\label{geom1}
\end{figure}

An immediate consequence of Lemma \ref{geometry1} is that the number of distinct angles in an equilateral pentagon is $1,3$ or $5$. If the number is $1$, then all the angles are equal, and we get the regular dodecahedron tiling.

\subsubsection*{Case. Three distinct angles}

Suppose the number of distinct angles in the equilateral pentagon is $3$. Denote the three distinct angles by $\alpha,\beta,\gamma$. By Lemma \ref{geometry1}, the pentagon has the angle combination $\alpha^2\beta^2\gamma,\alpha^2\beta\gamma^2,\alpha\beta^2\gamma^2$. By Lemma \ref{hdeg}, at least two of $\alpha,\beta,\gamma$ should appear at degree $3$ vertices. If only two appears, then up to the symmetry of exchanging the angles, by the two angle part of Table \ref{deg3AVC}, we may assume that $\alpha\beta^2$ is the only degree $3$ vertex. However, by Lemmas \ref{deg3a} and \ref{deg3b}, this implies that $\alpha$ appears at least twice in the pentagon, and $\beta$ appears at least three times. This implies $\gamma$ cannot appear in the pentagon. The contradiction  implies that all three angles must appear at degree $3$ vertices. This leads to the three angle part of Table \ref{deg3AVC}. 

For $\text{AVC}_3=\{\alpha\beta\gamma\opt\alpha^3\}$, if $\alpha^3$ is not a vertex, then each of $\alpha,\beta,\gamma$ appears at every degree $3$ vertex. By Lemma \ref{deg3a}, this implies each appears at least twice in the pentagon. Since this contradicts the pentagon, we conclude $\alpha^3$ must appear. Therefore $\text{AVC}_3=\{\alpha\beta\gamma,\alpha^3\}$, where $\alpha^3$ becomes necessary. Now $\alpha$ appears at every degree $3$ vertex. By Lemma \ref{deg3a}, this implies $\alpha$ appears at least twice in the pentagon. Then up to the exchange $\beta\leftrightarrow \gamma$, we may denote the case by $\{\alpha^2\beta^2\gamma\colon \alpha\beta\gamma,\alpha^3\}$. 

For $\text{AVC}_3=\{\alpha\beta^2,\alpha^2\gamma\}$, by Lemma \ref{deg3a}, $\alpha$ again appears at least twice in the pentagon. Therefore the pentagon has the angle combination $\alpha^2\beta^2\gamma$ or $\alpha^2\beta\gamma^2$. If we exchange $\beta$ and $\gamma$ in the second case $\{\alpha^2\beta\gamma^2\colon \alpha\beta^2,\alpha^2\gamma\}$, then we get two possible cases $\{\alpha^2\beta^2\gamma\colon \alpha\beta^2,\alpha^2\gamma\}$ and $\{\alpha^2\beta^2\gamma\colon \alpha^2\beta,\alpha\gamma^2\}$. 

For $\text{AVC}_3=\{\alpha\beta^2,\gamma^3\}$, if the pentagon has the angle combination $\alpha\beta^2\gamma^2$, then by the angle sums of $\alpha\beta^2,\gamma^3$ and the angle sum for pentagon (Lemma \ref{anglesum}), we have
\[
3\pi+\tfrac{4}{f}\pi
=\alpha+2\beta+2\gamma
=(\alpha+2\beta)+\tfrac{2}{3}\cdot 3\gamma
=2\pi+\tfrac{2}{3}\cdot 2\pi.
\]
Since this implies $f=12$, we may dismiss the case. Therefore the pentagon has the angle combination $\alpha^2\beta^2\gamma$ or $\alpha^2\beta\gamma^2$. If we exchange $\beta$ and $\gamma$ in the second case, then we get two possible cases  $\{\alpha^2\beta^2\gamma\colon \alpha\beta^2,\gamma^3\}$ and $\{\alpha^2\beta^2\gamma\colon \alpha\gamma^2,\beta^3\}$.

In summary, we get the following complete list of three distinct angles in equilateral pentagon.
\begin{description}
\item[3.1] $\{\alpha^2\beta^2\gamma\colon \alpha\beta\gamma,\alpha^3\}$. 
\item[3.2a] $\{\alpha^2\beta^2\gamma\colon \alpha\beta^2,\alpha^2\gamma\}$.
\item[3.2b] $\{\alpha^2\beta^2\gamma\colon \alpha^2\beta,\alpha\gamma^2\}$.
\item[3.3a] $\{\alpha^2\beta^2\gamma\colon \alpha\beta^2,\gamma^3\}$.
\item[3.3b] $\{\alpha^2\beta^2\gamma\colon \alpha\gamma^2,\beta^3\}$.
\end{description}

In the labels, the first digit $3$ means only three distinct angles appear at degree $3$ vertices. The second digit refers to the angle combinations at degree $3$ vertices according to Table \ref{deg3AVC}. The last alphabet refers to the possible variations by the consideration of angle combinations in the pentagon.

By Lemma \ref{geometry1}, a pentagon with the angle combination $\alpha^2\beta^2\gamma$ has two possible angle arrangements, given in Figure \ref{pent3}. Combined with five combinations above, we get total of $10$ cases.

\begin{figure}[htp]
\centering
\begin{tikzpicture}


\foreach \a in {0,...,4}
\foreach \x in {0,1}
\draw[xshift=3.5*\x cm, rotate=72*\a]
	(18:1) -- (90:1);
	
	
\node at (90:0.7) {\small $\gamma$};
\node at (18:0.7) {\small $\beta$};
\node at (162:0.7) {\small $\beta$};
\node at (234:0.7) {\small $\alpha$};
\node at (-54:0.7) {\small $\alpha$};

\node at (0,-1.2) {$[\alpha,\alpha,\beta,\gamma,\beta]$};


\begin{scope}[xshift=3.5cm]

\node at (90:0.7) {\small $\gamma$};
\node at (18:0.7) {\small $\alpha$};
\node at (162:0.7) {\small $\alpha$};
\node at (234:0.7) {\small $\beta$};
\node at (-54:0.7) {\small $\beta$};

\node at (0,-1.2) {$[\alpha,\gamma,\alpha,\beta,\beta]$};

\end{scope}

\end{tikzpicture}
\caption{Two angle arrangements for the pentagon $\alpha^2\beta^2\gamma$.}
\label{pent3}
\end{figure}

\subsubsection*{Case. Five distinct angles}

Suppose the number of distinct angles in the equilateral pentagon is $5$. We denote the distinct angles by $\alpha,\beta,\gamma,\delta,\epsilon$, and denote the pentagon by $\alpha\beta\gamma\delta\epsilon$. By Lemma \ref{deg3a}, there is no angle appearing at all degree $3$ vertices. By Lemma \ref{hdeg}, we may assume that $\alpha,\beta,\gamma,\delta$ appear at degree $3$ vertices. The discussion may be further divided by whether $\epsilon$ appears at degree $3$ vertices.

\subsubsection*{Subcase. Four angles at degree $3$}

Suppose $\epsilon$ does not appear at degree $3$ vertices. Then we may use the four angle part of Table \ref{deg3AVC} for what may appear at degree $3$ vertices.

For $\text{AVC}_3=\{\alpha\beta\gamma,\alpha\delta^2\opt \beta^2\delta\text{ or }\beta^3\}$, by $\alpha$ not appearing at all degree $3$ vertices, we know one of the two optional vertices necessarily appear, and we get two $\text{AVC}_3$s
\[
\{\alpha\beta\gamma,\alpha\delta^2,\beta^2\delta\},\quad
\{\alpha\beta\gamma,\alpha\delta^2,\beta^3\}.
\]
For $\text{AVC}_3=\{\alpha\beta\gamma,\alpha^2\delta\opt \beta^2\delta\text{ or }\beta^3\}$, the same reason gives two $\text{AVC}_3$s
\[
\{\alpha\beta\gamma,\alpha^2\delta,\beta\delta^2\},\quad
\{\alpha\beta\gamma,\alpha^2\delta,\beta^3\}.
\]
Since $\alpha\leftrightarrow\beta$ exchanges $\{\alpha\beta\gamma,\alpha\delta^2,\beta^2\delta\}$ and $\{\alpha\beta\gamma,\alpha^2\delta,\beta\delta^2\}$, the four $\text{AVC}_3$s may be reduced to three. 

For $\text{AVC}_3=\{\alpha\beta\gamma,\delta^3\}$, we use Lemma \ref{hdeg} to conclude that one of $\alpha\epsilon^3$, $\beta\epsilon^3$, $\gamma\epsilon^3$, $\delta\epsilon^3$, $\epsilon^4$, $\epsilon^5$ must be a vertex. Up to the symmetry of exchanging $\alpha,\beta,\gamma$, we may assume that one of $\alpha\epsilon^3$, $\delta\epsilon^3$, $\epsilon^4$, $\epsilon^5$ is a vertex.

For $\text{AVC}_3=\{\alpha\beta^2,\gamma\delta^2\opt \alpha^2\delta\}$, we may apply Lemma \ref{deg3c} to conclude that, if $\alpha^2\delta$ does not appear, then $\beta$ and $\delta$ together appear at least three times in the pentagon, a contradiction to five distinct angles. Therefore $\alpha^2\delta$ necessarily appears, and we get $\text{AVC}_3=\{\alpha\beta^2,\gamma\delta^2,\alpha^2\delta\}$. 

Finally, $\text{AVC}_3=\{\alpha\beta^2,\alpha^2\gamma,\delta^3\}$ already has three vertices.

In summary, we get the following complete list of five distinct angles in the pentagon, and four distinct angles at degree $3$ vertices. 
\begin{description}
\item[4.1a] 
$\{\alpha\beta\gamma\delta\epsilon\colon 
\alpha\beta\gamma,\alpha\delta^2,\beta^2\delta\}$,
$12$ arrangements.
\item[4.1b] 
$\{\alpha\beta\gamma\delta\epsilon\colon 
\alpha\beta\gamma,\alpha\delta^2,\beta^3\}$,
$12$ arrangements.
\item[4.1c] 
$\{\alpha\beta\gamma\delta\epsilon\colon 
\alpha\beta\gamma,\alpha^2\delta,\beta^3\}$,
$12$ arrangements.
\item[4.2a] 
$\{\alpha\beta\gamma\delta\epsilon\colon 
\alpha\beta\gamma,\delta^3,\alpha\epsilon^3\}$,
$6$ arrangements by $\beta\leftrightarrow\gamma$.
\item[4.2b] 
$\{\alpha\beta\gamma\delta\epsilon\colon 
\alpha\beta\gamma,\delta^3,\delta\epsilon^3\}$,
$2$ arrangements by $\alpha,\beta,\gamma$ exchange.
\item[4.2c] 
$\{\alpha\beta\gamma\delta\epsilon\colon 
\alpha\beta\gamma,\delta^3,\epsilon^4\}$,
$2$ arrangements by $\alpha,\beta,\gamma$ exchange.
\item[4.2d] 
$\{\alpha\beta\gamma\delta\epsilon\colon 
\alpha\beta\gamma,\delta^3,\epsilon^5\}$,
$2$ arrangements by $\alpha,\beta,\gamma$ exchange.
\item[4.3] 
$\{\alpha\beta\gamma\delta\epsilon\colon 
\alpha\beta^2,\gamma\delta^2,\alpha^2\delta\}$,
$12$ arrangements.
\item[4.4] 
$\{\alpha\beta\gamma\delta\epsilon\colon 
\alpha\beta^2,\alpha^2\gamma,\delta^3\}$,
$12$ arrangements.
\end{description}

Similar to the three distinct angle case, the first digit $4$ means only four angles appear at degree $3$ vertices. The second digit refers to the angle combinations at degree $3$ vertices according to Table \ref{deg3AVC}. 

We also remark that, up to the symmetry of flip and rotation, the pentagon $\alpha\beta\gamma\delta\epsilon$ has twelve possible angle arrangements in general
\begin{align*}
&
[\alpha, \beta, \gamma, \delta, \epsilon], 
[\alpha, \beta, \gamma, \epsilon, \delta], 
[\alpha, \beta, \delta, \gamma, \epsilon], 
[\alpha, \beta, \delta, \epsilon, \gamma], \\
&
[\alpha, \beta, \epsilon, \gamma, \delta], 
[\alpha, \beta, \epsilon, \delta, \gamma], 
[\alpha, \gamma, \beta, \delta, \epsilon], 
[\alpha, \gamma, \beta, \epsilon, \delta], \\
&
[\alpha, \gamma, \delta, \beta, \epsilon], 
[\alpha, \gamma, \epsilon, \beta, \delta], 
[\alpha, \delta, \beta, \gamma, \epsilon], 
[\alpha, \delta, \gamma, \beta, \epsilon].
\end{align*}
We denote the arrangements by
\[
A1=[\alpha, \beta, \gamma, \delta, \epsilon],\;\dotsc,\;
A12=[\alpha, \delta, \gamma, \beta, \epsilon].
\]
Then $A1$ is the arrangement in Figure \ref{pent}, and $A3$ is the arrangement in Figure \ref{pentagon}.

Further symmetries in some $\text{AVC}_3$s may reduce the number of arrangements we need to consider. For example, Case 4.2a is symmetric with respect to the exchange $\beta\leftrightarrow\gamma$. Since $\beta\leftrightarrow\gamma$ takes $A1=[\alpha, \beta, \gamma, \delta, \epsilon]$ to $A7=[\alpha, \gamma, \beta, \delta, \epsilon]$, we only need to consider $A1$ and do not need to consider $A7$. We also note that $\beta\leftrightarrow\gamma$ takes $A4=[\alpha, \beta, \delta, \epsilon, \gamma]$ to $[\alpha, \gamma, \delta, \epsilon, \beta]=[\alpha, \beta, \epsilon, \delta, \gamma]=A6$, where we get the first equality by flipping the pentagon. Taking into account of the various arrangements, we get total of $72$ cases.

\subsubsection*{Subcase. Five angles at degree $3$}

Suppose all five distinct angles appear at degree $3$ vertices. Then we may use the five angle part of Table \ref{deg3AVC}.

We know no angle can appear at all the degree $3$ vertices. For the first $\text{AVC}_3=\{\alpha\beta\gamma,\alpha\delta\epsilon\opt \cdots\}$ in the five angle part of Table \ref{deg3AVC}, therefore, some optional vertex not involving $\alpha$ must appear. Up to the symmetry of symbols, we may assume that either $\gamma\epsilon^2$ or $\gamma^3$ appears. This gives $\text{AVC}_3=\{\alpha\beta\gamma,\alpha\delta\epsilon, \gamma\epsilon^2 \opt \cdots\}$ or $\text{AVC}_3=\{\alpha\beta\gamma,\alpha\delta\epsilon, \gamma^3 \opt \cdots\}$.

For the other $\text{AVC}_3$s in the five angle part of Table \ref{deg3AVC}, we always get three necessary vertices with the only exception of $\text{AVC}_3=\{\alpha\beta\gamma,\delta\epsilon^2\opt \alpha^3\}$. If we include the case that $\alpha^3$ also appears for this particular $\text{AVC}_3$, then we get the following complete list of five distinct angles at degree $3$ vertices. 
\begin{description}
\item[5.1a] 
$\{\alpha\beta\gamma\delta\epsilon\colon 
\alpha\beta\gamma,\alpha\delta\epsilon,\gamma\epsilon^2\}$,
$12$ arrangements.
\item[5.1b] 
$\{\alpha\beta\gamma\delta\epsilon\colon 
\alpha\beta\gamma,\alpha\delta\epsilon,\gamma^3\}$,
$6$ arrangements by $\delta\leftrightarrow\epsilon$.
\item[5.2] 
$\{\alpha\beta\gamma\delta\epsilon\colon 
\alpha\beta\gamma,\alpha\delta^2,\alpha^2\epsilon\}$,
$6$ arrangements by $\beta\leftrightarrow\gamma$.
\item[5.3] 
$\{\alpha\beta\gamma\delta\epsilon\colon 
\alpha\beta\gamma,\alpha\delta^2,\beta\epsilon^2\}$,
$8$ arrangements by $(\alpha,\delta)\leftrightarrow(\beta,\epsilon)$.
\item[5.4] 
$\{\alpha\beta\gamma\delta\epsilon\colon 
\alpha\beta\gamma,\alpha\delta^2,\beta^2\epsilon\}$,
$12$ arrangements.
\item[5.5] 
$\{\alpha\beta\gamma\delta\epsilon\colon 
\alpha\beta\gamma,\alpha\delta^2,\delta\epsilon^2\}$,
$6$ arrangements by $\beta,\gamma$ exchange.
\item[5.6] 
$\{\alpha\beta\gamma\delta\epsilon\colon 
\alpha\beta\gamma,\alpha\delta^2,\epsilon^3\}$,
$6$ arrangements by $\beta,\gamma$ exchange.
\item[5.7] 
$\{\alpha\beta\gamma\delta\epsilon\colon 
\alpha\beta\gamma,\alpha^2\delta,\beta^2\epsilon\}$,
$8$ arrangements by $(\alpha,\delta),(\beta,\epsilon)$ exchange.
\item[5.8] 
$\{\alpha\beta\gamma\delta\epsilon\colon 
\alpha\beta\gamma,\alpha^2\delta,\delta^2\epsilon\}$,
$6$ arrangements by $\beta\leftrightarrow\gamma$.
\item[5.9] 
$\{\alpha\beta\gamma\delta\epsilon\colon 
\alpha\beta\gamma,\alpha^2\delta,\epsilon^3\}$,
$6$ arrangements by $\beta\leftrightarrow\gamma$.
\item[5.10] 
$\{\alpha\beta\gamma\delta\epsilon\colon 
\alpha\beta\gamma,\delta\epsilon^2,\alpha^3\}$,
$6$ arrangements by $\beta\leftrightarrow\gamma$.
\item[5.11] 
$\{\alpha\beta\gamma\delta\epsilon\colon 
\alpha\beta^2,\gamma\delta^2,\alpha^2\epsilon\}$,
$12$ arrangements.
\item[5.12] 
$\{\alpha\beta\gamma\delta\epsilon\colon 
\alpha\beta^2,\gamma\delta^2,\epsilon^3\}$,
$8$ arrangements by $(\alpha,\beta)\leftrightarrow(\gamma,\delta)$.
\end{description}
Taking into account of various arrangements, we get total of 102 cases.

\subsubsection*{Subcase. Extend $\text{AVC}_3=\{\alpha\beta\gamma,\delta\epsilon^2\}$ by degree $4$}

The cases listed so far have enough number of angle equations for completely determining the equilateral pentagon. The only remaining case we need to consider is $\text{AVC}_3=\{\alpha\beta\gamma,\delta\epsilon^2\}$, which means that $\alpha\beta\gamma$ and $\delta\epsilon^2$ are the only degree $3$ vertices. To get the third angle sum equation (necessary for determining the equilateral pentagon), we consider all the possible vertices of degree $4$ or $5$. After exhausting all these cases, we are left with the exceptional case that $\text{AVC}_3=\{\alpha\beta\gamma,\delta\epsilon^2\}$ and $v_4=v_5=0$. 

If $\text{AVC}_3=\{\alpha\beta\gamma,\delta\epsilon^2\}$, then by the angle sum, $\alpha\beta\gamma\cdots$ and $\delta\epsilon^2\cdots$ cannot be vertices of degree $>3$. Moreover, the angle sums of $\alpha\beta\gamma,\delta\epsilon^2$ and the angle sum for pentagon (Lemma \ref{anglesum}) imply
\[
\alpha+\beta+\gamma=2\pi,\quad
\delta=\tfrac{8}{f}\pi,\quad
\epsilon=(1-\tfrac{4}{f})\pi.
\]
By $f>12$, we get $\epsilon>\frac{2}{3}\pi$. Therefore $\epsilon^3\cdots$ is not a vertex. 

Suppose $\text{AVC}_3=\{\alpha\beta\gamma,\delta\epsilon^2\}$, and there is a degree $4$ vertex. Then we add degree $4$ angle combinations that are not $\alpha\beta\gamma\cdots,\delta\epsilon^2\cdots,\epsilon^3\cdots$ to get the following complete list, up to the symmetry of exchanging $\alpha,\beta,\gamma$.
\begin{description}
\item[1.1] 
$\{\alpha\beta\gamma\delta\epsilon\colon 
\alpha\beta\gamma,\delta\epsilon^2,\alpha\beta\delta\epsilon\}$,
$6$ arrangements by $\alpha\leftrightarrow\beta$.
\item[1.2a] 
$\{\alpha\beta\gamma\delta\epsilon\colon 
\alpha\beta\gamma,\delta\epsilon^2,\alpha\beta^2\delta\}$,
$12$ arrangements. 
\item[1.2b] 
$\{\alpha\beta\gamma\delta\epsilon\colon 
\alpha\beta\gamma,\delta\epsilon^2,\alpha\beta^2\epsilon\}$,
$12$ arrangements.
\item[1.2c] 
$\{\alpha\beta\gamma\delta\epsilon\colon 
\alpha\beta\gamma,\delta\epsilon^2,\alpha\beta\delta^2\}$,
$6$ arrangements by $\alpha\leftrightarrow\beta$.
\item[1.2d] 
$\{\alpha\beta\gamma\delta\epsilon\colon 
\alpha\beta\gamma,\delta\epsilon^2,\alpha\beta\epsilon^2\}$,
$6$ arrangements by $\alpha\leftrightarrow\beta$.
\item[1.2e] 
$\{\alpha\beta\gamma\delta\epsilon\colon 
\alpha\beta\gamma,\delta\epsilon^2,\alpha\delta^2\epsilon\}$,
$6$ arrangements by $\beta\leftrightarrow\gamma$.
\item[1.2f] 
$\{\alpha\beta\gamma\delta\epsilon\colon 
\alpha\beta\gamma,\delta\epsilon^2,\alpha^2\delta\epsilon\}$,
$6$ arrangements by $\beta\leftrightarrow\gamma$.
\item[1.3a] 
$\{\alpha\beta\gamma\delta\epsilon\colon 
\alpha\beta\gamma,\delta\epsilon^2,\alpha^2\beta^2\}$,
$6$ arrangements by $\alpha\leftrightarrow\beta$.
\item[1.3b] 
$\{\alpha\beta\gamma\delta\epsilon\colon 
\alpha\beta\gamma,\delta\epsilon^2,\alpha^2\delta^2\}$,
$6$ arrangements by $\beta\leftrightarrow\gamma$.
\item[1.3c] 
$\{\alpha\beta\gamma\delta\epsilon\colon 
\alpha\beta\gamma,\delta\epsilon^2,\alpha^2\epsilon^2\}$,
$6$ arrangements by $\beta\leftrightarrow\gamma$.
\item[1.4a] 
$\{\alpha\beta\gamma\delta\epsilon\colon 
\alpha\beta\gamma,\delta\epsilon^2,\alpha\beta^3\}$,
$12$ arrangements.
\item[1.4b] 
$\{\alpha\beta\gamma\delta\epsilon\colon 
\alpha\beta\gamma,\delta\epsilon^2,\alpha\delta^3\}$,
$6$ arrangements by $\beta\leftrightarrow\gamma$.
\item[1.4c] 
$\{\alpha\beta\gamma\delta\epsilon\colon 
\alpha\beta\gamma,\delta\epsilon^2,\alpha^3\delta\}$,
$6$ arrangements by $\beta\leftrightarrow\gamma$.
\item[1.4d] 
$\{\alpha\beta\gamma\delta\epsilon\colon 
\alpha\beta\gamma,\delta\epsilon^2,\alpha^3\epsilon\}$,
$6$ arrangements by $\beta\leftrightarrow\gamma$.
\item[1.4e] 
$\{\alpha\beta\gamma\delta\epsilon\colon 
\alpha\beta\gamma,\delta\epsilon^2,\delta^3\epsilon\}$,
$2$ arrangements by $\alpha,\beta,\gamma$ exchange.
\item[1.5a] 
$\{\alpha\beta\gamma\delta\epsilon\colon 
\alpha\beta\gamma,\delta\epsilon^2,\alpha^4\}$,
$6$ arrangements by $\beta\leftrightarrow\gamma$.
\item[1.5b] 
$\{\alpha\beta\gamma\delta\epsilon\colon 
\alpha\beta\gamma,\delta\epsilon^2,\delta^4\}$,
$2$ arrangements by $\alpha,\beta,\gamma$ exchange.
\end{description}

In the labels, the first digit $1$ (because only $3,4,5$ have been used so far) means the inclusion of a degree $4$ vertex. The second digit refers to the types $1111,112,22,13,4$ of degree $4$ vertices. The last alphabet refers to the possible variations. We combine the cases that are the same after exchanging the angles ($\{\alpha\beta\gamma\delta\epsilon\colon 
\alpha\beta\gamma,\delta\epsilon^2,\alpha\gamma\delta\epsilon\}$ is the same as $\{\alpha\beta\gamma\delta\epsilon\colon 
\alpha\beta\gamma,\delta\epsilon^2,\alpha\beta\delta\epsilon\}$ after $\beta\leftrightarrow\gamma$). Taking into account of various arrangements, we get total of 112 cases.

\subsubsection*{Subcase. Extend $\text{AVC}_3=\{\alpha\beta\gamma,\delta\epsilon^2\}$ by degree $5$}

Suppose $\text{AVC}_3=\{\alpha\beta\gamma,\delta\epsilon^2\}$, and there is a degree $5$ vertex. Then we add degree $5$ angle combinations that are not $\alpha\beta\gamma\cdots,\delta\epsilon^2\cdots,\epsilon^3\cdots$ to get the following complete list.
\begin{description}
\item[2.1a] 
$\{\alpha\beta\gamma\delta\epsilon\colon 
\alpha\beta\gamma,\delta\epsilon^2,\alpha\beta^2\delta\epsilon\}$,
$12$ arrangements.
\item[2.1b] 
$\{\alpha\beta\gamma\delta\epsilon\colon 
\alpha\beta\gamma,\delta\epsilon^2,\alpha\beta\delta^2\epsilon\}$,
$6$ arrangements by $\alpha\leftrightarrow\beta$.
\item[2.2a] 
$\{\alpha\beta\gamma\delta\epsilon\colon 
\alpha\beta\gamma,\delta\epsilon^2,\alpha^2\beta^2\delta\}$,
$6$ arrangements by $\alpha\leftrightarrow\beta$.
\item[2.2b] 
$\{\alpha\beta\gamma\delta\epsilon\colon 
\alpha\beta\gamma,\delta\epsilon^2,\alpha^2\beta^2\epsilon\}$,
$6$ arrangements by $\alpha\leftrightarrow\beta$.
\item[2.2c] 
$\{\alpha\beta\gamma\delta\epsilon\colon 
\alpha\beta\gamma,\delta\epsilon^2,\alpha^2\delta^2\epsilon\}$,
$6$ arrangements by $\beta\leftrightarrow\gamma$.
\item[2.2d] 
$\{\alpha\beta\gamma\delta\epsilon\colon 
\alpha\beta\gamma,\delta\epsilon^2,\alpha\beta^2\delta^2\}$,
$12$ arrangements.
\item[2.2e] 
$\{\alpha\beta\gamma\delta\epsilon\colon 
\alpha\beta\gamma,\delta\epsilon^2,\alpha\beta^2\epsilon^2\}$,
$12$ arrangements.
\item[2.3a] 
$\{\alpha\beta\gamma\delta\epsilon\colon 
\alpha\beta\gamma,\delta\epsilon^2,\alpha^3\delta\epsilon\}$,
$6$ arrangements by $\beta\leftrightarrow\gamma$.
\item[2.3b] 
$\{\alpha\beta\gamma\delta\epsilon\colon 
\alpha\beta\gamma,\delta\epsilon^2,\alpha\beta^3\delta\}$,
$12$ arrangements.
\item[2.3c] 
$\{\alpha\beta\gamma\delta\epsilon\colon 
\alpha\beta\gamma,\delta\epsilon^2,\alpha\beta^3\epsilon\}$,
$12$ arrangements.
\item[2.3d] 
$\{\alpha\beta\gamma\delta\epsilon\colon 
\alpha\beta\gamma,\delta\epsilon^2,\alpha\delta^3\epsilon\}$,
$6$ arrangements by $\beta\leftrightarrow\gamma$.
\item[2.3e] 
$\{\alpha\beta\gamma\delta\epsilon\colon 
\alpha\beta\gamma,\delta\epsilon^2,\alpha\beta\delta^3\}$,
$6$ arrangements by $\alpha\leftrightarrow\beta$.
\item[2.4a] 
$\{\alpha\beta\gamma\delta\epsilon\colon 
\alpha\beta\gamma,\delta\epsilon^2,\alpha^2\beta^3\}$,
$12$ arrangements.
\item[2.4b] 
$\{\alpha\beta\gamma\delta\epsilon\colon 
\alpha\beta\gamma,\delta\epsilon^2,\alpha^2\delta^3\}$,
$6$ arrangements by $\beta\leftrightarrow\gamma$.
\item[2.4c] 
$\{\alpha\beta\gamma\delta\epsilon\colon 
\alpha\beta\gamma,\delta\epsilon^2,\alpha^3\delta^2\}$,
$6$ arrangements by $\beta\leftrightarrow\gamma$.
\item[2.4d] 
$\{\alpha\beta\gamma\delta\epsilon\colon 
\alpha\beta\gamma,\delta\epsilon^2,\alpha^3\epsilon^2\}$,
$6$ arrangements by $\beta\leftrightarrow\gamma$.
\item[2.5a] 
$\{\alpha\beta\gamma\delta\epsilon\colon 
\alpha\beta\gamma,\delta\epsilon^2,\alpha\beta^4\}$,
$12$ arrangements.
\item[2.5b] 
$\{\alpha\beta\gamma\delta\epsilon\colon 
\alpha\beta\gamma,\delta\epsilon^2,\alpha\delta^4\}$,
$6$ arrangements by $\beta\leftrightarrow\gamma$.
\item[2.5c] 
$\{\alpha\beta\gamma\delta\epsilon\colon 
\alpha\beta\gamma,\delta\epsilon^2,\alpha^4\delta\}$,
$6$ arrangements by $\beta\leftrightarrow\gamma$.
\item[2.5d] 
$\{\alpha\beta\gamma\delta\epsilon\colon 
\alpha\beta\gamma,\delta\epsilon^2,\alpha^4\epsilon\}$,
$6$ arrangements by $\beta\leftrightarrow\gamma$.
\item[2.5e] 
$\{\alpha\beta\gamma\delta\epsilon\colon \alpha\beta\gamma,\delta\epsilon^2,\delta^4\epsilon\}$,
$2$ arrangements by $\alpha,\beta,\gamma$ exchange.
\item[2.6a] 
$\{\alpha\beta\gamma\delta\epsilon\colon \alpha\beta\gamma,\delta\epsilon^2,\alpha^5\}$,
$6$ arrangements by $\beta\leftrightarrow\gamma$.
\item[2.6b] 
$\{\alpha\beta\gamma\delta\epsilon\colon \alpha\beta\gamma,\delta\epsilon^2,\delta^5\}$,
$2$ arrangements by $\alpha,\beta,\gamma$ exchange.
\end{description}

The first digit $2$ means the inclusion of a degree $5$ vertex. The second digit refers to the types $1112,122,113,23,14,5$ of degree $5$ vertices. The last alphabet refers to the possible variations. Taking into account of various arrangements, we get total of 172 cases.

The final remaining exceptional case is that $\alpha\beta\gamma$ and $\delta\epsilon^2$ are the only degree $3$ vertices, and there are no vertices of degree $4$ or $5$. In Section \ref{casespecial}, we will show that the only possibility is that $\delta^6$ is a vertex, and $\delta,\epsilon$ are not adjacent in the pentagon.

\section{Calculation of Equilateral Pentagon}
\label{calculation}

Consider the spherical equilateral pentagon in Figure \ref{pent}, with edge length $a$ and five angles $\alpha,\beta,\gamma,\delta,\epsilon$ arranged as $A1$. We have $a<\pi$ because otherwise any two adjacent edges would intersect at two points, contradicting the simple pentagon requirement. By \cite[Lemma 3]{ay1}, we may calculate the great arc $x$ connecting $\beta$ and $\epsilon$ vertices in two ways. From the isosceles triangle above $x$, we get
\[
\cos x
=\cos^2a+\sin\alpha\sin^2a.
\]
From the quadrilateral below $x$, we get
\begin{align*}
\cos x
&=(1-\cos\gamma)(1-\cos\delta)\cos^3a
-\sin\gamma\sin\delta\cos^2a \nonumber \\
&\quad +(\cos\gamma+\cos\delta-\cos\gamma\cos\delta)\cos a
+\sin\gamma\sin\delta.
\end{align*}
Identifying the right side of the two equations, we get a cubic equation for $\cos a$. Dividing $1-\cos a$, we get a quadratic equation for $\cos a$ 
\begin{equation}\label{lmneq}
L\cos^2 a+M\cos a+N=0,
\end{equation}
where the coefficients depend only on $\alpha,\gamma,\delta$
\begin{align*}
L &=(1-\cos\gamma)(1-\cos\delta), \\
M &=\cos\alpha+\cos(\gamma+\delta)-\cos\gamma-\cos\delta, \\
N &=\cos\alpha-\sin\gamma\sin\delta.
\end{align*}
We call \eqref{lmneq} the {\em $LMN$-equation} for $BE$ (connecting $\beta$-vertex and $\epsilon$-vertex). The validity of the $LMN$-equation means exactly the triangle and the quadrilateral can be glued together to form a pentagon. However, we do not yet know the pentagon is simple. 

\begin{figure}[htp]
\centering
\begin{tikzpicture}


\draw[dashed]
	(18:1) -- node[fill=white,inner sep=0.5] {\small $x$} (162:1);
	
\node at (90:0.7) {\small $\alpha$};
\node at (162:0.7) {\small $\beta$};
\node at (234:0.7) {\small $\gamma$};
\node at (-54:0.7) {\small $\delta$};
\node at (18:0.7) {\small $\epsilon$};


\foreach \a in {0,...,4}
{
\draw[rotate=72*\a]
	(18:1) -- (90:1);
\node at (54+72*\a:0.95) {\small $a$};
}

\end{tikzpicture}
\caption{Spherical equilateral pentagon.}
\label{pent}
\end{figure}

The $LMN$-equation can be used for five pairs $(\beta,\epsilon)$, $(\alpha,\gamma)$, $(\beta,\delta)$, $(\gamma,\epsilon)$, $(\alpha,\delta)$, and we get five quadratic equations 
\[
L_i\cos^2 a+M_i\cos a+N_i=0,\quad
i=1,2,3,4,5.
\]
Let $t=\cos a$, $x_1=\cos\alpha$, $y_1=\sin\alpha$, and similarly introduce $x_i,y_i$, $i=2,3,4,5$ for $\beta,\gamma,\delta,\epsilon$. Then the five quadratic equations are polynomial equations of the $11$ variables $t,x_i,y_i$, $i=1,2,3,4,5$. Together with $x_i^2+y_i^2-1=0$, we have total of 10 polynomial equations. The ideal generated by the 10 polynomials has Hilbert dimension $3$. 

Conceptually this means that three more independent polynomial equations can determine isolated solutions. These are provided by the three angle sums in our cases. For example, the angle sum of $\alpha\beta\gamma$ is the same as 
\[
\cos(\alpha+\beta+\gamma)=1,\quad
\sin(\alpha+\beta+\gamma)=0,
\]
which are degree $3$ polynomial equations of $x_1,x_2,x_3,y_1,y_2,y_3$. Of course we may also use simpler degree $2$ polynomial equations such as $\cos(\alpha+\beta)=\cos\gamma$ and $\sin(\alpha+\beta)=-\sin\gamma$. 

We can implement strict inequalities such as $x_i\ne 1$ (i.e., angles are nonzero) by saturating the polynomial ideal with $x_i-1$. We also introduce saturations with $x_i-x_j$ to make sure angles are distinct, for the case of five distinct angles.

For each case, we then get finitely many solutions from a suitable Gr\"obner basis of the saturated ideal. All these can be calculated symbolically. For example, for the first arrangement $A1=[\alpha,\beta,\gamma,\delta,\epsilon]$ of Case 1.1, we get the following Gr\"obner basis
\begin{align*}
&
16x_5^4+8x_5^3-8x_5^2+1, \;
x_5^2+y_5^2-1, \;
8x_5^3+8x_5^2+2x_2-1, \\
&
8x_5^3y_5-4x_5y_5+y_2, \;
16x_5^3+8x_5^2+2x_1-8x_5-1,  \;
8x_5^3y_5+4x_5^2y_5+y_1, \\
&
2x_5y_5+y_4,  \;
-2x_5^2+x_4+1,  \;
-x_5+x_3, \;
y_5+y_3, \\
&
-64x_5^3-96x_5^2+17t+4x_5+21.
\end{align*}
Starting with $4$ roots of $x_5$ from $16x_5^4+8x_5^3-8x_5^2+1$, we may successively find the roots of all $11$ variables. The total number of solutions is $8$.

So far, the calculation has been symbolic. Therefore we may get numerical values to arbitrary accuracy. Then we use the numerical values to exclude non-real solutions. This can be done by calculating the numerical values with sufficiently many digits, and require that the first enough number of digits of the imaginary part to be all zero. For Case 1.1($A1$), 4 of 8 solutions are not real. 
 
For the real solutions, by $0<a<\pi$, the value $t=\cos a$ uniquely determines $a$. Moreover, since we have included $x_i^2+y_i^2-1$ in the polynomial ideal, the real value pairs $x_i,y_i$ uniquely determine five angles within the interval $(0,2\pi)$. Then we verify that the angle sums of all vertices are $2\pi$. For example, we have the following angles for one real solution of Case 1.1($A1$)
\[
\alpha=1.045\pi, \;
\beta=1.564\pi, \;
\gamma=1.391\pi, \;
\delta=0.782\pi, \;
\epsilon=0.609\pi.
\]
We calculate the angle sums for Case 1.1
\[
\alpha+\beta+\gamma=4\pi,\;
\delta+2\epsilon=2\pi,\;
\alpha+\beta+\delta+\epsilon=4\pi,
\]
and get $4\pi$ instead of $2\pi$. Altogether, $3$ real solutions fail this angle sum test. The only real solution passing the test is 
\[
\alpha=0.508\pi, \;
\beta=0.394\pi, \;
\gamma=1.098\pi, \;
\delta=0.197\pi, \;
\epsilon=0.902\pi.
\]
Then we use the angle sum for pentagon (Lemma \ref{anglesum}) to calculate the number $f$ of tiles. The number must be an even integer $f\ge 16$. For the solution above, we get $f=40.644$, which fails this tiling number test. We conclude that Case 1.1($A1$) does not yield any tiling of the sphere by congruent equilateral pentagons.

We remark that the discussion so far can be applied to Cases 3.*. In fact, the calculation is much simpler because we only need to introduce $x_1,x_2,x_3,y_1,y_2,y_3$ for the three distinct angles. We find no real solutions passing the angle sum test and the tiling number test. Therefore Cases 3.* give no suitable solutions.

For five distinct angles, we need to further verify that Lemma \ref{geometry1} is satisfied in five ways for those real solutions that pass both the angle sum test and the tiling number test. It is not hard to reformulate the criterion as the following: The biggest and smallest angles are adjacent. Moreover, if we switch the order of these two angles, then all five angles are ordered monotonically.

For example, Case 4.2b($A1$) has 18 complex solutions, among which 10 are not real, and 3 real solutions fail the angle sum test. All the remaining 5 real solutions pass the tiling number test with $f=36$, and the angles are the following
\begin{align*}
& 
\alpha =0.480\pi, \;
\beta =0.714\pi, \;
\gamma =0.806\pi, \;
\delta =0.667\pi, \;
\epsilon =0.444\pi; \\
& 
\alpha =0.295\pi, \;
\beta =1.625\pi, \;
\gamma =0.080\pi, \;
\delta =0.667\pi, \;
\epsilon =0.444\pi; \\
&
\alpha =0.157\pi, \;
\beta =1.011\pi, \;
\gamma =0.832\pi, \;
\delta =0.667\pi, \;
\epsilon =0.444\pi; \\
& 
\alpha =0.876\pi, \;
\beta =0.430\pi, \;
\gamma =0.694\pi, \;
\delta =0.667\pi, \;
\epsilon =0.444\pi; \\
& 
\alpha =1.094\pi, \;
\beta =0.619\pi, \;
\gamma =0.286\pi, \;
\delta =0.667\pi, \;
\epsilon =0.444\pi.
\end{align*}
We find that the first and fifth solutions violate the criterion because the maximum angle and the minimum angle are not adjacent. Moreover, the angles in the third solution are already in monotonic order $\alpha,\epsilon,\delta,\gamma,\beta$ (this is $A1$). Therefore switching the maximum and minimum angles violates the criterion. The remaining second and fourth solutions satisfy the criterion, and are candidates for tiling the sphere.

After all the calculations and applying tests and criteria to filter out pentagons not suitable for tiling, we get the following complete list of candidate pentagons for tiling the sphere. We provide the approximate value of $\cos a$ only for cases that will yield tilings in Section \ref{tiling}. 

\begin{align*}
\intertext{
Case 4.2b:
$\{\alpha\beta\gamma,\delta^3,\delta\epsilon^3\}$. }  
A1 &\colon 
f=36, \;
\alpha=0.29539\pi, \;
\beta=1.62453\pi, \;
\gamma=0.08007\pi, \;
\delta=\tfrac{2}{3}\pi, \;
\epsilon=\tfrac{4}{9}\pi; \\
A1 &\colon 
f=36, \;
\alpha=0.87573\pi, \;
\beta=0.42997\pi, \;
\gamma=0.69428\pi, \;
\delta=\tfrac{2}{3}\pi, \;
\epsilon=\tfrac{4}{9}\pi; \\
A3 &\colon 
f=36, \;
\alpha=0.85571\pi, \;
\beta=0.45589\pi, \;
\gamma=0.68839\pi, \;
\delta=\tfrac{2}{3}\pi, \;
\epsilon=\tfrac{4}{9}\pi. \\
\intertext{
Case 4.2c:
$\{\alpha\beta\gamma,\delta^3,\epsilon^4\}$. } 
A1 &\colon 
f=24, \;
\alpha=0.27849\pi, \;
\beta=1.59984\pi, \;
\gamma=0.12166\pi, \;
\delta=\tfrac{2}{3}\pi, \;
\epsilon=\tfrac{1}{2}\pi; \\
A1 &\colon 
f=24, \;
\alpha=0.82020\pi, \;
\beta=0.48453\pi, \;
\gamma=0.69526\pi, \;
\delta=\tfrac{2}{3}\pi, \;
\epsilon=\tfrac{1}{2}\pi; \\
A3 &\colon 
f=24, \;
\alpha=0.80106\pi, \;
\beta=0.51139\pi, \;
\gamma=0.68753\pi, \;
\delta=\tfrac{2}{3}\pi, \;
\epsilon=\tfrac{1}{2}\pi, \\
&\cos a=0.85341. \\
\intertext{ 
Case 4.2d: 
$\{\alpha\beta\gamma,\delta^3,\epsilon^5\}$. } 
A1 &\colon 
f=60, \;
\alpha=0.31031\pi, \;
\beta=1.64260\pi, \;
\gamma=0.04708\pi, \;
\delta=\tfrac{2}{3}\pi, \;
\epsilon=\tfrac{2}{5}\pi; \\
A1 &\colon 
f=60, \;
\alpha=0.92294\pi, \;
\beta=0.38907\pi, \;
\gamma=0.68797\pi, \;
\delta=\tfrac{2}{3}\pi, \;
\epsilon=\tfrac{2}{5}\pi; \\
A3 &\colon 
f=60, \;
\alpha=0.90594\pi, \;
\beta=0.40930\pi, \;
\gamma=0.68475\pi, \;
\delta=\tfrac{2}{3}\pi, \;
\epsilon=\tfrac{2}{5}\pi, \\
&\cos a=0.93133. \\
\intertext{
Case 5.5: 
$\{\alpha\beta\gamma,\alpha\delta^2,\delta\epsilon^2\}$. }  
A5 &\colon 
f=24, \;
\alpha=\tfrac{4}{3}\pi, \;
\beta=0.14400\pi, \;
\gamma=0.52265\pi, \;
\delta=\tfrac{1}{3}\pi, \;
\epsilon=\tfrac{5}{6}\pi, \\
&\cos a=0.70688. \\
\intertext{ 
Cases 1.2e, 1.5a, 2.4b: 
$\{\alpha\beta\gamma,\delta\epsilon^2,
\alpha\delta^2\epsilon\text{ or }
\alpha^4\text{ or }
\alpha^2\delta^3\}$. } 
A11 &\colon 
f=24, \;
\alpha=\tfrac{1}{2}\pi, \;
\beta=1.38071\pi, \;
\gamma=0.11928\pi, \;
\delta=\tfrac{1}{3}\pi, \;
\epsilon=\tfrac{5}{6}\pi, \\
&\cos a=0.68125. \\
\intertext{ 
Cases 1.4e, 2.6b: 
$\{\alpha\beta\gamma,\delta\epsilon^2,
\delta^3\epsilon\text{ or }\delta^5\}$. } 
A1 &\colon 
f=20, \;
\alpha=0.60551\pi, \;
\beta=0.50248\pi, \;
\gamma=0.89199\pi, \;
\delta=\tfrac{2}{5}\pi, \;
\epsilon=\tfrac{4}{5}\pi; \\
A3 &\colon 
f=20, \;
\alpha=0.30959\pi, \;
\beta=1.06152\pi, \;
\gamma=0.62888\pi, \;
\delta=\tfrac{2}{5}\pi, \;
\epsilon=\tfrac{4}{5}\pi, \\
&\cos a=0.77680. \\
\intertext{ 
Case 1.5b:
$\{\alpha\beta\gamma,\delta\epsilon^2,
\delta^4\}$. } 
A1 &\colon 
f=16, \;
\alpha=0.10133\pi, \;
\beta=1.56723\pi, \;
\gamma=0.33142\pi, \;
\delta=\tfrac{1}{2}\pi, \;
\epsilon=\tfrac{3}{4}\pi; \\
A1 &\colon 
f=16, \;
\alpha=0.63380\pi, \;
\beta=0.56425\pi, \;
\gamma=0.80193\pi, \;
\delta=\tfrac{1}{2}\pi, \;
\epsilon=\tfrac{3}{4}\pi; \\
A3 &\colon 
f=16, \;
\alpha=0.45368\pi, \;
\beta=0.88238\pi, \;
\gamma=0.66392\pi, \;
\delta=\tfrac{1}{2}\pi, \;
\epsilon=\tfrac{3}{4}\pi; \\
&\cos a=0.77943.
\intertext{
Case 2.5e: 
$\{\alpha\beta\gamma,\delta\epsilon^2,
\delta^4\epsilon\}$. } 
A1 &\colon 
f=28, \;
\alpha=0.55888\pi, \;
\beta=0.43715\pi, \;
\gamma=1.00396\pi, \;
\delta=\tfrac{2}{7}\pi, \;
\epsilon=\tfrac{6}{7}\pi. \\
\intertext{Exceptional Case $\text{AVC}_3=\{\alpha\beta\gamma,\delta\epsilon^2\}$, $v_4=v_5=0$. We calculate for $\{\alpha\beta\gamma,\delta\epsilon^2,\delta^6\}$ and get} 
A1 &\colon 
f=24, \;
\alpha=0.58056\pi, \;
\beta=0.46336\pi, \;
\gamma=0.95606\pi, \;
\delta=\tfrac{1}{3}\pi, \;
\epsilon=\tfrac{5}{6}\pi; \\
A3 &\colon 
f=24, \;
\alpha=0.14400\pi, \;
\beta=\tfrac{4}{3}\pi, \;
\gamma=0.52265\pi, \;
\delta=\tfrac{1}{3}\pi, \;
\epsilon=\tfrac{5}{6}\pi, \\
&\cos a=0.70688; \\
A3 &\colon 
f=24, \;
\alpha=0.11928\pi, \;
\beta=1.38071\pi, \;
\gamma=\tfrac{1}{2}\pi, \;
\delta=\tfrac{1}{3}\pi, \;
\epsilon=\tfrac{5}{6}\pi, \\
&\cos a=0.68125. 
\end{align*}

We remark that we know the exact values of the length $a$ and all five angles, in the sense of the polynomials satisfied by their cosine (and which real root of the polynomial). In fact, it turns out all the polynomials can be solved by radicals, and we have the explicit formula for the cosines in terms of radicals.

In the list above, we choose the accuracy up to five digits because this is sufficient for the subsequent argument in Section \ref{tiling}, that constructs tilings from the approximate values. 

We group 1.2e, 1.5a, 2.4b together because the numerical values are the same up to five digits. Therefore the argument for these cases are the same in Section \ref{tiling}. The same remark applies to 1.4e and 2.6b.

In 4.2b, 4.2c, 4.2d, 1.4e, 2.6b, 1.5b, 2.5e, the precise values of $\delta,\epsilon$ can be calculated from the angle sums of two existing vertices involving only $\delta,\epsilon$. Therefore the precise values can be used in the argument in Section \ref{tiling}. 

In 5.5, 1.2e, 1.5a, 2.4b, the precise values of $\alpha,\delta,\epsilon$ are inferred from numerical calculations, and should be treated here as accurate up to the fifth digit. Although we can rigorously justify these precise values by more symbolic calculation, we will only use five digit approximations of all five angles in constructing tilings in Section \ref{tiling}. After the construction, we further rigorously justify these precise values.

In the exceptional case, we will argue that $\delta^6$ must be a vertex. The precise values of $\delta,\epsilon$ can be calculated from the angle sums of $\delta\epsilon^2,\delta^6$, and the precise value of $\beta$ and $\gamma$ are inferred from numerical calculations.

\section{Tiling by Congruent Equilateral Pentagon}
\label{tiling}

All pentagons that are potentially suitable for tilings are listed near the end of Section \ref{calculation}. Based on the data in the list, we try to construct the tiling. To compare tilings from different cases, we change all the arrangements to $A1=[\alpha,\beta,\gamma,\delta,\epsilon]$ or $A3=[\alpha, \beta, \delta, \gamma, \epsilon]$. This means exchanging $\alpha,\beta$ ($A5$ changed to $A3$) in Case 5.5, and exchanging $\alpha,\gamma$ ($A11$ changed to $A3$) in Cases 1.2e, 1.5a, 2.4b.
\begin{align*}
\intertext{
Case 5.5: 
$\{\alpha\beta\gamma,\beta\delta^2,\delta\epsilon^2\}$. }  
A3 &\colon 
f=24, \;
\alpha=0.14400\pi, \;
\beta=\tfrac{4}{3}\pi, \;
\gamma=0.52265\pi, \;
\delta=\tfrac{1}{3}\pi, \;
\epsilon=\tfrac{5}{6}\pi, \\
&\cos a=0.70688. \\
\intertext{ 
Cases 1.2e, 1.5a, 2.4b: 
$\{\alpha\beta\gamma,\delta\epsilon^2,
\gamma\delta^2\epsilon\text{ or }
\gamma^4\text{ or }
\gamma^2\delta^3\}$. } 
A3 &\colon 
f=24, \;
\alpha=0.11928\pi, \;
\beta=1.38071\pi, \;
\gamma=\tfrac{1}{2}\pi, \;
\delta=\tfrac{1}{3}\pi, \;
\epsilon=\tfrac{5}{6}\pi, \\
&\cos a=0.68125. 
\end{align*}

\subsection{Case 4.2b}
\label{case42b}

To construct tilings for the solution
\[
A1\colon 
f=36, \;
\alpha=0.29539\pi, \;
\beta=1.62453\pi, \;
\gamma=0.08007\pi, \;
\delta=\tfrac{2}{3}\pi, \;
\epsilon=\tfrac{4}{9}\pi,
\]
we first find the AVC. Due to non-precise values for $\alpha,\beta,\gamma$, we cannot find all the possible angle combinations $\alpha^a\beta^b\gamma^c\delta^d\epsilon^e$ at vertices by precisely solving the angle sum equation
\[
\alpha a+\beta b+\gamma c+\delta d+\epsilon e
=2\pi.
\]
Instead, since all angles and multiplicities are non-negative, the approximate values imply 
\[
a\le 6,\;
b\le 1,\;
c\le 24,\;
d\le 3,\;
e\le 4.
\]
By the precise values of $\delta,\epsilon$, this implies that any solution of the exact angle sum equation must satisfy the estimation
\[
|0.29539a+1.62453b+0.08007c+\tfrac{2}{3}d+\tfrac{4}{9}e-2| 
\le 0.00001(a+b+c).
\]
We substitute all combinations of non-negative integers $a,b,c,d,e$ within the bounds to the inequality above and find that only three combinations $\alpha\beta\gamma$, $\delta^3$, $\delta\epsilon^3$ satisfy the inequality. Therefore we conclude
\[
\text{AVC}
=\{\alpha\beta\gamma,\delta^3,\delta\epsilon^3\}.
\]

Applying the same argument to the other two solutions, we get the same AVC. In fact, the bounds for the multiplicities are much smaller for these two solutions, and four or three digits are sufficient for deriving the AVC. 

Next we derive tilings from the AVC. We used the notation adopted in \cite{wy1,wy2}. We denote by $T_i$ the  $i$-th tile, and indicate the tile as circled $i$. We also denote by $\theta_i$ the angle $\theta$ in $T_i$. The notations are unambiguous because the five angles have distinct values.

We will also use the AAD (adjacent angle deduction) introduced in \cite[Section 2.5]{wy1}. We note that, due to the change of adjacency between angles, for the same AVC, the AAD argument for $A1$ and $A3$ are different.

The vertex $\delta\epsilon^3$ has consecutive $\epsilon\epsilon\epsilon$. For the $A1$ arrangement, this implies the AAD $\thin\epsilon\thin^{\alpha}\epsilon^{\delta}\thin\epsilon\thin$ at the vertex. By no $\alpha\delta\cdots$ (due to the AVC), the AAD is $\thin\epsilon\thin^{\alpha}\epsilon^{\delta}\thin^{\delta}\epsilon^{\alpha}\thin$. This implies a vertex $\thin^{\gamma}\delta^{\epsilon}\thin^{\epsilon}\delta^{\gamma}\thin\cdots=\delta^3=\thin^{\epsilon}\delta^{\gamma}\thin\delta\thin^{\gamma}\delta^{\epsilon}\thin=\thin^{\epsilon}\delta^{\gamma}\thin^{\gamma}\delta^{\epsilon}\thin^{\gamma}\delta^{\epsilon}\thin$. Then $\gamma\epsilon\cdots$ is a vertex, contradicting to the AVC. 

For the $A3$ arrangement, by the AVC (we will henceforth omit mentioning AVC), we know $\gamma^2\cdots$ is not a vertex. This implies the unique AAD $\thin^{\beta}\delta^{\gamma}\thin^{\alpha}\epsilon^{\gamma}\thin^{\alpha}\epsilon^{\gamma}\thin^{\alpha}\epsilon^{\gamma}\thin$ of $\delta\epsilon^3$. Therefore we determine four tiles $T_1,T_2,T_3,T_4$ around a vertex $\delta\epsilon^3$ in the first of Figure \ref{case42fig}. By $\alpha_1\gamma_2\cdots=\alpha_1\beta_5\gamma_2$ and no $\alpha_5\delta_2\cdots$, we determine $T_5$. By $\delta_2\delta_5\cdots=\delta_2\delta_5\delta_6$ and no $\gamma_5\gamma_6\cdots$, we determine $T_6$. By $\beta_6\gamma_5\cdots=\alpha_7\beta_6\gamma_5$ and no $\beta_7\epsilon_5\cdots$, we determine $T_7$. By starting with $\alpha_2\gamma_3\cdots$ (in place of $\alpha_1\gamma_2\cdots$) and applying the same reason, we determine $T_8,T_9,T_{10}$. By the unique AAD of $\epsilon_6\epsilon_8\epsilon_{10}\cdots=\delta_{11}\epsilon_6\epsilon_8\epsilon_{10}$, we determine $T_{11}$. By the unique AAD of $\delta_7\epsilon_{11}\cdots=\delta_7\epsilon_{11}\epsilon\epsilon_{12}$, we determine $T_{12}$. By $\alpha_{12}\gamma_7\cdots=\alpha_{12}\beta_{13}\gamma_7$ and no $\alpha_{13}\epsilon_7\cdots$, we determine $T_{13}$. By the unique AAD of $\delta_{13}\epsilon_5\epsilon_7\cdots=\delta_{13}\epsilon_5\epsilon_7\epsilon_{14}$, we determine $T_{14}$. Then $\beta_4\gamma_1\cdots=\alpha\beta\gamma$ and $\delta_1\delta_{14}\cdots=\delta^3$ imply $\alpha,\delta$ are adjacent, a contradiction. 

We conclude that Case 4.2b has no tiling.

\begin{figure}[htp]
\centering
\begin{tikzpicture}[>=latex,scale=1]


\foreach \x in {0,1,2,3}
\draw[rotate=90*\x]
	(0,0) -- (0.8,0) -- (1.2,0.6) -- (0.6,1.2) -- (0,0.8);

\foreach \x in {1,2}
\draw[rotate=-90*\x]
	(-0.6,1.2) -- (-0.6,1.8) -- (1.8,1.8) -- (1.8,0.6) -- (1.2,0.6) 
	(0.6,1.2) -- (0.6,1.8);
	
\draw
	(-1.8,-1.8) -- (-1.8,-2.5) -- (0.6,-2.5)
	(0.6,-1.8) -- (0.6,-2.5) -- (2.5,-2.5) -- (2.5,0.6)
	(0.6,1.2) -- (0.6,1.8) -- (2.5,1.8) -- (3,1.2) -- (2.5,0.6)
	(1.8,1.8) -- (1.8,0.6) 
	(1.8,0.6) -- (2.5,0.6) -- (2.5,-2.5)
	(1.8,-1.8) -- (3,-1.8) -- (3.2,-0.3) -- (3,1.2);
	
\node at (0.2,0.2) {\small $\epsilon$};
\node at (0.7,0.2) {\small $\alpha$};
\node at (0.95,0.55) {\small $\beta$};
\node at (0.55,0.95) {\small $\delta$};
\node at (0.2,0.7) {\small $\gamma$};
\node at (0.2,-0.2) {\small $\epsilon$};
\node at (0.7,-0.2) {\small $\gamma$};
\node at (0.95,-0.6) {\small $\delta$};
\node at (0.6,-0.95) {\small $\beta$};
\node at (0.2,-0.7) {\small $\alpha$};
\node at (-0.2,0.2) {\small $\delta$};
\node at (-0.75,0.2) {\small $\gamma$};
\node at (-0.95,0.6) {\small $\epsilon$};
\node at (-0.6,0.95) {\small $\alpha$};
\node at (-0.2,0.65) {\small $\beta$};
\node at (-0.2,-0.2) {\small $\epsilon$};
\node at (-0.7,-0.2) {\small $\alpha$};
\node at (-0.95,-0.55) {\small $\beta$};
\node at (-0.55,-0.95) {\small $\delta$};
\node at (-0.2,-0.7) {\small $\gamma$};
\node at (1.05,0) {\small $\beta$};
\node at (1.3,-0.4) {\small $\delta$};
\node at (1.3,0.4) {\small $\alpha$};
\node at (1.6,-0.4) {\small $\gamma$};
\node at (1.6,0.4) {\small $\epsilon$};
\node at (1.2,-0.8) {\small $\delta$};
\node at (1.6,-0.85) {\small $\beta$};
\node at (1.6,-1.6) {\small $\alpha$};
\node at (0.8,-1.6) {\small $\epsilon$};
\node at (0.8,-1.25) {\small $\gamma$};
\node at (0,-1.05) {\small $\beta$};
\node at (0.4,-1.25) {\small $\alpha$};
\node at (-0.4,-1.25) {\small $\delta$};
\node at (0.4,-1.6) {\small $\epsilon$};
\node at (-0.4,-1.6) {\small $\gamma$};
\node at (-0.8,-1.6) {\small $\beta$};
\node at (-0.8,-1.25) {\small $\delta$};
\node at (-1.3,-0.8) {\small $\gamma$};
\node at (-1.6,-0.8) {\small $\epsilon$};
\node at (-1.6,-1.6) {\small $\alpha$};
\node at (-0.6,-2) {\small $\alpha$};
\node at (0.4,-2) {\small $\epsilon$};
\node at (-1.6,-2) {\small $\beta$};
\node at (-1.6,-2.3) {\small $\delta$};
\node at (0.4,-2.3) {\small $\gamma$};
\node at (0.8,-2) {\small $\delta$};
\node at (0.8,-2.3) {\small $\beta$};
\node at (1.8,-2) {\small $\gamma$};
\node at (2.3,-2.3) {\small $\alpha$};
\node at (2.3,-2) {\small $\epsilon$};
\node at (2,0.4) {\small $\epsilon$};
\node at (2,-1.6) {\small $\beta$};
\node at (2.3,-1.55) {\small $\delta$};
\node at (2,-0.6) {\small $\alpha$};
\node at (2.3,0.4) {\small $\gamma$};
\node at (2.7,-1.6) {\small $\epsilon$};
\node at (2.7,0.5) {\small $\alpha$};
\node at (2.8,1.2) {\small $\alpha$};
\node at (2,0.8) {\small $\delta$};
\node at (2.45,0.8) {\small $\beta$};
\node at (2,1.6) {\small $\gamma$};
\node at (2.4,1.65) {\small $\epsilon$};
\node at (1.25,0.8) {\small $\gamma$};
\node at (1.6,0.8) {\small $\epsilon$};
\node at (1.6,1.6) {\small $\alpha$};
\node at (0.8,1.6) {\small $\beta$};
\node at (0.8,1.3) {\small $\delta$};
\node at (0.45,1.3) {\small $\delta$};
\node at (0,1.05) {\small $\alpha$};
\node at (2.7,-2) {\small $\epsilon$};

\node[inner sep=0.5,draw,shape=circle] at (0.5,0.5) {\small $1$};
\node[inner sep=0.5,draw,shape=circle] at (0.5,-0.5) {\small $2$};
\node[inner sep=0.5,draw,shape=circle] at (-0.5,-0.5) {\small $3$};
\node[inner sep=0.5,draw,shape=circle] at (-0.5,0.5) {\small $4$};
\node[inner sep=0.5,draw,shape=circle] at (1.45,0) {\small $5$};
\node[inner sep=0.5,draw,shape=circle] at (1.3,-1.3) {\small $6$};
\node[inner sep=0.5,draw,shape=circle] at (2.15,0) {\small $7$};
\node[inner sep=0.5,draw,shape=circle] at (0,-1.45) {\small $8$};
\node[inner sep=0.5,draw,shape=circle] at (-1.3,-1.3) {\small $9$};
\node[inner sep=0,draw,shape=circle] at (0,-2.15) {\small $10$};
\node[inner sep=0,draw,shape=circle] at (1.3,-2.15) {\small $11$};
\node[inner sep=0,draw,shape=circle] at (2.85,-0.3) {\small $12$};
\node[inner sep=0,draw,shape=circle] at (2.3,1.25) {\small $13$};
\node[inner sep=0,draw,shape=circle] at (1.25,1.25) {\small $14$};


\begin{scope}[xshift=5.5cm]

\foreach \x in {0,1,2,3}
{
\begin{scope}[rotate=90*\x]

\draw
	(0,0) -- (0.8,0) -- (1.2,0.6) -- (0.6,1.2) -- (0,0.8);

\node at (0.2,0.2) {\small $\epsilon$};
\node at (0.7,0.2) {\small $\alpha$};
\node at (0.95,0.57) {\small $\beta$};
\node at (0.6,0.95) {\small $\delta$};
\node at (0.2,0.7) {\small $\gamma$};

\end{scope}
}

\draw
	(0.6,1.2) -- (0.6,1.8) -- (1.8,1.8) -- (1.8,-0.6) -- (1.2,-0.6)
	(1.2,0.6) -- (2.3,0.6);

\node at (1.05,0.05) {\small $\beta$};
\node at (1.3,-0.4) {\small $\delta$};
\node at (1.3,0.4) {\small $\alpha$};
\node at (1.65,-0.4) {\small $\gamma$};
\node at (1.65,0.4) {\small $\epsilon$};
\node at (1.25,0.8) {\small $\gamma$};
\node at (1.65,0.8) {\small $\epsilon$};
\node at (1.65,1.6) {\small $\alpha$};
\node at (0.8,1.55) {\small $\beta$};
\node at (0.8,1.25) {\small $\delta$};
\node at (1.95,0.4) {\small $\epsilon$};
\node at (1.95,0.8) {\small $\epsilon$};

\node at (1.25,-0.8) {\small $\delta$};
\node at (0.4,1.25) {\small $\delta$};

\node[inner sep=0.5,draw,shape=circle] at (0.5,0.5) {\small $1$};
\node[inner sep=0.5,draw,shape=circle] at (0.5,-0.5) {\small $2$};
\node[inner sep=0.5,draw,shape=circle] at (-0.5,-0.5) {\small $3$};
\node[inner sep=0.5,draw,shape=circle] at (-0.5,0.5) {\small $4$};
\node[inner sep=0.5,draw,shape=circle] at (1.45,0) {\small $5$};
\node[inner sep=0.5,draw,shape=circle] at (1.3,1.3) {\small $6$};

\end{scope}

\end{tikzpicture}
\caption{Tilings for $\{\alpha\beta\gamma,\delta^3,\delta\epsilon^3\}$ and $\{\alpha\beta\gamma,\delta^3,\epsilon^4\}$.}
\label{case42fig}
\end{figure}

\subsection{Cases 4.2c and 4.2d}
\label{case42c}

By the method in Section \ref{case42b}, we find $\text{AVC}=\{\alpha\beta\gamma,\delta^3,\epsilon^4\}$ for the three pentagons in Case 4.2c, and $\text{AVC}=\{\alpha\beta\gamma,\delta^3,\epsilon^5\}$ for the three pentagons in Case 4.2d.  Again we need to consider two arrangements $A1$ and $A3$. 

Again we have consecutive $\epsilon\epsilon\epsilon$ at $\epsilon^4$ or $\epsilon^5$. For the $A1$ arrangement, we may use the same argument for Case 2.4b to show there is no tiling.

For the $A3$ arrangement and $\text{AVC}=\{\alpha\beta\gamma,\delta^3,\epsilon^4\}$, the AVC implies the unique AAD $\thin^{\alpha}\epsilon^{\gamma}\thin^{\alpha}\epsilon^{\gamma}\thin^{\alpha}\epsilon^{\gamma}\thin^{\alpha}\epsilon^{\gamma}\thin$ of $\epsilon^4$. This determines four tiles $T_1,T_2,T_3,T_4$ around a vertex $\epsilon^4$ in the second of Figure \ref{case42fig}. By $\alpha_1\gamma_2\cdots=\alpha_1\beta_5\gamma_2$ and no $\beta_1\delta_5\cdots$, we determine $T_5$. By $\alpha_5\beta_1\cdots=\alpha_5\beta_1\gamma_6$ and $\epsilon_5\cdots=\epsilon^4$, we determine $T_6$. Similarly, we may start with $\alpha_2\gamma_3\cdots,\alpha_3\gamma_4\cdots,\alpha_4\gamma_1\cdots$ and get more tiles similar to $T_5,T_6$, and get more $\epsilon^4$. We may also apply the argument starting from the original $\epsilon^4$ to these new $\epsilon^4$. More repetitions give the pentagonal subdivision of the octahedron, in the second of Figure \ref{subdivision_tiling}. 

The argument for Case 4.2d is completely similar, starting from the unique AAD $\thin^{\alpha}\epsilon^{\gamma}\thin^{\alpha}\epsilon^{\gamma}\thin^{\alpha}\epsilon^{\gamma}\thin^{\alpha}\epsilon^{\gamma}\thin^{\alpha}\epsilon^{\gamma}\thin$ of $\epsilon^5$. We get the pentagonal subdivision of the icosahedron, in the third of Figure \ref{subdivision_tiling}.

Although we get the tiling, we still need to verify that the solutions of Cases 4.2c($A3$) and 4.2d($A3$) can be realised by simple pentagons. 

The first of Figure \ref{calculation_division} is the pentagonal subdivision of a regular triangular face of the regular tetrahedron ($n=3$), octahedron ($n=4$, Cases 4.2c), icosahedron ($n=5$, Cases 4.2d). The angles of the pentagon are arranged as $A3$, and $\delta=\frac{2}{3}\pi,\epsilon=\frac{2}{n}\pi$ match the angles for the two cases.

\begin{figure}[htp]
\centering
\begin{tikzpicture}[>=latex,scale=1]


\begin{scope}[shift={(-3cm,0.8cm)}]

\foreach \a in {0,1,2}
{
\begin{scope}[rotate=120*\a]

\draw[dotted]
	(-30:2) -- (90:2);

\draw
	(90:2) -- ++(-140:1.2)
	(90:2) -- ++(-80:1.2)
	(0,0) -- (76:0.86)
	(6:1.53) -- (76:0.86);

\fill 
	(90:2) circle (0.1);
	
\end{scope}
}

\draw[gray!50]
	(0,0) -- (0,-1) -- (-30:2) -- cycle;

\draw[dotted]
	(0,2) -- (0,0) -- (30:2) -- (-30:2);

\filldraw[fill=white] 
	(0,0) circle (0.1);

\filldraw[fill=gray] 
	(0,-1) circle (0.1);
			
\node at (1.35,0.05) {\small $\alpha$};	
\node at (0.3,0.5) {\small $\beta$};
\node at (0.7,-0.45) {\small $\gamma$};
\node at (10:0.25) {\small $\delta$};
\node at (-26:1.7) {\small $\epsilon$};

\node at (-0.6,-1.2) {\small $\alpha$};	
\node at (0.3,-0.6) {\small $\beta$};
\node at (-0.7,-0.4) {\small $\gamma$};
\node at (-110:0.3) {\small $\delta$};
\node at (-1.4,-0.95) {\small $\epsilon$};

\node at (0.65,-0.8) {\small $\alpha$};

\node at (1,-0.9) {\scriptsize $\triangle$};
\node at (0.8,0) {\scriptsize $\triangle'$};
\node at (0.45,1) {\scriptsize $\triangle_1$};
\node at (1.5,0.6) {\scriptsize $\triangle_2$};

\end{scope}


\draw[gray!50]
	(0,0) -- (0,2.6) -- (4,0) -- (0,0);

\draw[dashed]
	(2,1.3) -- ++(-1.3,-2);

\draw
	(0,2.6) -- (1.6,0.7) -- (4,0)
	(1.6,0.7) -- (0,0);
	
\fill
	(4,0) circle (0.1);

\filldraw[fill=white]
	(0,2.6) circle (0.1);

\filldraw[fill=gray] 
	(0,0) circle (0.1);
	
\node[fill=white,inner sep=1] at (1.7,1.45) {\small $x$};
\node[fill=white,inner sep=1] at (0,1.3) {\small $y$};
\node[fill=white,inner sep=1] at (2,0) {\small $z$};

\node[gray] at (0.35,0.25) {\small $\frac{1}{2}\pi$};
\node[gray] at (0.35,2.05) {\small $\frac{1}{3}\pi$};
\node[gray] at (3.1,0.25) {\small $\frac{1}{n}\pi$};
	
\node at (0.4,0.5) {\small $\theta$};

\node[fill=white,inner sep=0.5] at (0.8,0.35) {\small $\frac{1}{2}a$};
\node[fill=white, inner sep=0.5] at (2.2,0.5) {\small $a$};
\node[fill=white, inner sep=0.5] at (0.9,1.5) {\small $a$};

\node at (2.2,1.5) {\small $P$};
\node at (1.3,-0.25) {\small $Q$};
\node at (3,-0.5) {triangle $\triangle$};
\node at (1.6,0.5) {\small $\alpha$};
\node at (1.25,0.75) {\small $\beta$};
\node at (1.7,0.85) {\small $\gamma$};

\end{tikzpicture}
\caption{Equilateral pentagon for pentagonal subdivision.}
\label{calculation_division}
\end{figure}

The right triangle $\triangle$ in the second of Figure \ref{calculation_division} is one sixth of the regular triangular face, with respective angles $\frac{1}{n}\pi,\frac{1}{3}\pi,\frac{1}{2}\pi$ at $\bullet,\circ,{\color{gray} \bullet}$. We also recall the key values from the Platonic solids. 
\renewcommand{\arraystretch}{1.3}
\begin{center}
 \begin{tabular}{|c| c|| c c c c c c|} 
 \hline
 $f$ & $n$ & $\cos x$ & $\cos y$ & $\cos z$ & $x$ & $y$ & $z$ \\   
 \hline\hline
 24 & 4 
 & $\frac{1}{\sqrt{3}}$
 & $\frac{\sqrt{2}}{\sqrt{3}}$ 
 & $\frac{1}{\sqrt{2}}$ 
 & $0.3040\pi$
 & $0.1959\pi$
 & $0.2500\pi$ \\
 \hline
 60 & 5  
 & $\frac{\sqrt{5}+1}{\sqrt{6(5-\sqrt{5})}}$
 & $\frac{\sqrt{5}+1}{2\sqrt{3}}$
 & $\frac{\sqrt{2}}{\sqrt{5-\sqrt{5}}}$ 
 & $0.2076\pi$
 & $0.1161\pi$
 & $0.1762\pi$ \\ 
 \hline
\end{tabular}
\end{center}

The distances from the vertex $\alpha\beta\gamma$ to $\bullet,\circ,{\color{gray} \bullet}$ are respectively $a,a,\frac{1}{2}a$.  We draw the line starting from  the middle point $P$ of $x$ and orthogonal to $x$. By $y<z$, the line intersects $z$ at $Q$. For any point $X$ on $PQ$, we have $X\circ=X\bullet$. Then we calculate to find 
\begin{align*}
f=24 &\colon
P{\color{gray} \bullet}=0.1718\pi>\tfrac{1}{4}x,\;
Q\bullet=0.2011\pi>\tfrac{2}{3}z; \\
f=60 &\colon
P{\color{gray} \bullet}=0.1088\pi>\tfrac{1}{4}x,\;
Q\bullet=0.1260\pi>\tfrac{2}{3}z.
\end{align*}
By $P\bullet=\tfrac{1}{2}x$, this implies $X{\color{gray} \bullet}>\frac{1}{2}X\bullet$ at $X=P$, and $X{\color{gray} \bullet}<\frac{1}{2}X\bullet$ at $X=Q$. By the intermediate value theorem, this implies we have $X\circ=X\bullet=2X{\color{gray} \bullet}$ for some $X$ on $PQ$. This $X$ is the vertex $\alpha\beta\gamma$. 

We have the other right triangles $\triangle',\triangle_1,\triangle_2$ that are one sixth of the regular triangular face. We have divided $\triangle$ into three smaller triangles, with respective top angles $\alpha,\beta,\gamma$ and respective base edges $z,y,x$. Then the equilateral pentagon is obtained by glueing the smaller triangles to $\triangle'$ along the base edges. Since $\triangle',\triangle_1,\triangle_2$ do not overlap, the pentagon is simple. The construction of the pentagon shows that it tiles the pentagonal subdivision. 

Next, we calculate the pentagon. Let $\theta$ be the angle between $y$ and $\frac{1}{2}a$. Then we have
\[
\cos a=\cos\tfrac{1}{2}a\cos y+\sin\tfrac{1}{2}a\sin y\cos\theta.
\]
The angle between $z$ and $\frac{1}{2}a$ is $\frac{1}{2}\pi-\theta$. Therefore we have
\[
\cos a=\cos\tfrac{1}{2}a\cos z+\sin\tfrac{1}{2}a\sin z\sin\theta.
\]
We eliminate $\theta$ from the two equations and get an equation for $\cos\frac{1}{2}a$, $\cos a$, $\sin a$. This can be further converted to a precise polynomial equation for $\cos a$ 
\begin{align}
f=24 &\colon
25\cos^4a
+4(1-2\sqrt{3})\cos^3a \nonumber \\
&\quad -2(1+4\sqrt{3})\cos^2a
+4\cos a 
+1=0, \label{edge24} \\
f=60 &\colon
(63-11\sqrt{5})\cos^4a \nonumber \\
&\quad +(-4{\textstyle \sqrt{6}\sqrt{5+\sqrt{5}}}+16-8\sqrt{5})\cos^3a \nonumber \\
&\quad +(-4{\textstyle \sqrt{6}\sqrt{5+\sqrt{5}}}-14+6\sqrt{5})\cos^2a \nonumber \\
&\quad +(4+4\sqrt{5})\cos a 
+3+\sqrt{5}=0. \label{edge60}
\end{align}

The polynomial \eqref{edge24} has four real roots, and the one compatible with the approximate value $\cos a=0.85341$ is the biggest one. By $0<a<\pi$, this uniquely and precisely determines $a$. The approximate value is $a=0.17452\pi$. The polynomial \eqref{edge60} can be factorized as the product of two quadratic polynomials, and the solution compatible with the approximate value $\cos a=0.93133$ is
\begin{align*}
\cos a
&=\frac{1}{6728}( 
(193+7\sqrt{5})\sqrt{6}{\textstyle \sqrt{5+\sqrt{5}}} 
-372+736\sqrt{5}) \\
&+\frac{1+\sqrt{5}}{3364}\sqrt{
(57122+25659\sqrt{5})\sqrt{6}{\textstyle \sqrt{5+\sqrt{5}}} 
-135949\sqrt{5}
+247403
}.
\end{align*} 

We further determine the precise values of $\alpha,\beta,\gamma$ by 
\begin{align*}
\cos x
&=\cos^2 a+\sin^2 a\cos\gamma, \\
\cos y
&=\cos\tfrac{1}{2}a\cos a+\sin\tfrac{1}{2}a\sin a\cos\beta, \\
\cos z
&=\cos\tfrac{1}{2}a\cos a+\sin\tfrac{1}{2}a\sin a\cos\alpha. 
\end{align*}
The approximate values are 
\begin{align*}
f=24 &\colon
\alpha=0.80106\pi,\;
\beta=0.51139\pi,\;
\gamma=0.68753\pi; \\
f=60 &\colon
\alpha=0.90594\pi,\;
\beta=0.40930\pi,\;
\gamma=0.68475\pi.
\end{align*}

\subsection{Case 5.5}
\label{case55}

We start with the version of the solution of Case 5.5 given at the beginning of Section \ref{tiling}. This means we apply the exchange $\alpha\leftrightarrow \beta$ to the original solution for $A5$, and get the new version for $A3$. By the method in Section \ref{case42b}, we get
\[
\text{AVC}= 
\{\alpha\beta\gamma,\beta\delta^2,\delta\epsilon^2,\alpha^3\gamma^3,\alpha^2\gamma^2\delta^2,\alpha\gamma\delta^4,\delta^6\}.
\]
We have $v_i=0$ for $i\ne 3,6$. Then by $f=24$, and the vertex counting equation \cite{wy1,yan}
\[
\tfrac{f}{2}-6
={\textstyle \sum}_{k\ge 4}(k-3)v_k=v_4+2v_5+3v_6+\cdots,
\]
we get $v_6=2$. By \cite[Theorem 6]{yan}, the tiling is the earth map tiling with exactly two degree $6$ vertices as the poles. In particular, there is a $3^5$-tile.

Figure \ref{case55A} is the neighbourhood of a $3^5$-tile, in which all tiles have the $A3$ arrangement. We assume the angles of the center tile $T_1$ are given as indicated. By the AVC, the five degree $3$ vertices of $T_1$ belong to $\text{AVC}_3=\{\alpha\beta\gamma,\beta\delta^2,\delta\epsilon^2\}$. We also know $\alpha\epsilon\cdots$ is not a vertex.

We have $\gamma_1\cdots=\alpha\beta\gamma=\alpha_2\beta_6\gamma_1$ or $\alpha_6\beta_2\gamma_1$, given by the two pictures. In the first picture, by $\beta_6$ and no $\alpha_6\epsilon_1\cdots$, we determine $T_6$. Then by $\alpha_2$ and no $\alpha_6\epsilon_2\cdots$, we determine $T_2$. Then $\delta_1\epsilon_2\cdots=\delta_1\epsilon_2\epsilon_3$, and $\delta_6\epsilon_1\cdots=\delta_6\epsilon_1\epsilon_5$. By $\epsilon_5$, we get $\alpha_1\cdots=\alpha_1\beta_4\gamma_5$. Then $\gamma_5,\epsilon_5$ determine $T_5$. By $\beta_4,\epsilon_3$, we get $\beta_1\cdots=\alpha_4\beta_1\gamma_3$. Then $\gamma_3,\epsilon_3$ determine $T_3$, and $\alpha_4,\beta_4$ determine $T_4$.

\begin{figure}[htp]
\centering
\begin{tikzpicture}[>=latex,scale=1]

\foreach \a in {0,1}
{
\begin{scope}[xshift=3.5*\a cm]

\foreach \x in {1,...,5}
{
\draw[rotate=72*\x]
	(90:0.7) -- (18:0.7) -- (18:1.3) -- (54:1.7) -- (90:1.3);
	
\coordinate (X\x A\a) at (18+72*\x:0.45);
\coordinate (Y\x A\a) at (6+72*\x:0.8);
\coordinate (Z\x A\a) at (30+72*\x:0.8);
\coordinate (U\x A\a) at (10+72*\x:1.14);
\coordinate (V\x A\a) at (26+72*\x:1.14);
\coordinate (W\x A\a) at (-18+72*\x:1.45);
}

\node[draw,shape=circle, inner sep=0.5] at (0,0) {\small 1};
\foreach \x in {2,...,6}
\node[draw,shape=circle, inner sep=0.5] at (198-72*\x:1) {\small \x};

\end{scope}
}


\node at (X1A0) {\small $\gamma$};
\node at (X2A0) {\small $\epsilon$};
\node at (X3A0) {\small $\alpha$};
\node at (X4A0) {\small $\beta$};
\node at (X5A0) {\small $\delta$};

\node at (Y1A0) {\small $\alpha$};
\node at (Y2A0) {\small $\delta$};
\node at (Y3A0) {\small $\gamma$};
\node at (Y4A0) {\small $\alpha$};
\node at (Y5A0) {\small $\epsilon$};

\node at (Z1A0) {\small $\beta$};
\node at (Z2A0) {\small $\epsilon$};
\node at (Z3A0) {\small $\beta$};
\node at (Z4A0) {\small $\gamma$};
\node at (Z5A0) {\small $\epsilon$};

\node at (U1A0) {\small $\beta$};
\node at (U2A0) {\small $\gamma$};
\node at (U3A0) {\small $\delta$};
\node at (U4A0) {\small $\epsilon$};
\node at (U5A0) {\small $\alpha$};

\node at (V1A0) {\small $\alpha$};
\node at (V2A0) {\small $\alpha$};
\node at (V3A0) {\small $\delta$};
\node at (V4A0) {\small $\delta$};
\node at (V5A0) {\small $\gamma$};

\node at (W1A0) {\small $\delta$};
\node at (W2A0) {\small $\epsilon$};
\node at (W3A0) {\small $\beta$};
\node at (W4A0) {\small $\gamma$};
\node at (W5A0) {\small $\beta$};


\node at (X1A1) {\small $\gamma$};
\node at (X2A1) {\small $\epsilon$};
\node at (X3A1) {\small $\alpha$};
\node at (X4A1) {\small $\beta$};
\node at (X5A1) {\small $\delta$};

\node at (Y1A1) {\small $\beta$};
\node at (Y2A1) {\small $\epsilon$};
\node at (Y4A1) {\small $\delta$};
\node at (Y5A1) {\small $\beta$};

\node at (Z1A1) {\small $\alpha$};
\node at (Z2A1) {\small $\delta$};
\node at (Z4A1) {\small $\delta$};
\node at (Z5A1) {\small $\delta$};

\node at (U1A1) {\small $\alpha$};
\node at (U2A1) {\small $\gamma$};
\node at (U5A1) {\small $\alpha$};

\node at (V1A1) {\small $\beta$};
\node at (V4A1) {\small $\gamma$};
\node at (V5A1) {\small $\gamma$};

\node at (W1A1) {\small $\epsilon$};
\node at (W2A1) {\small $\delta$};
\node at (W5A1) {\small $\epsilon$};


\end{tikzpicture}
\caption{Case 5.5: Neighbourhood tiling of a $3^5$-tile.}
\label{case55A}
\end{figure}

In the second picture, by $\beta_2$ and $\alpha_6$, we get $\delta_1\cdots=\delta_1\delta_2\beta_3$ and $\epsilon_1\cdots=\delta_5\epsilon_1\epsilon_6$. Then $\beta_2,\delta_2$ determine $T_2$, and $\alpha_6,\epsilon_6$ determine $T_6$. By $\beta_3$, we have $\beta_1\cdots=\alpha_3\beta_1\gamma_4$ or $\beta_1\delta_3\delta_4$. We also know $\alpha_1\cdots=\alpha\beta\gamma$. We find $\beta_1\cdots=\alpha_3\beta_1\gamma_4$ and $\alpha_1\cdots=\alpha\beta\gamma$ imply either two $\gamma$ in $T_4$, or $\beta,\gamma$ adjacent in $T_4$, a contradiction. Therefore we have $\beta_1\cdots=\beta_1\delta_3\delta_4$, as in the second picture. Then $\alpha_1\cdots=\alpha\beta\gamma=\thin\alpha\thin^{\delta}\beta^{\alpha}\thin^{\epsilon}\gamma^{\delta}\thin$ implies $\alpha\epsilon\cdots$ is a vertex, a contradiction.

We conclude that the first of Figure \ref{case55A} is the only $3^5$-neighbourhood tiling fitting the AVC. 

By \cite[Theorem 6]{yan}, there are five families of earth map tilings, corresponding to distances $5,4,3,2,1$ between the two poles. They are obtained by glueing copies of the ``timezones'' in Figure \ref{timezone} (three timezones are shown for distance $5$) along the ``meridians''. The vertical edges at the top meet at the north pole, and the vertical edges at the bottom meet at the south pole. For $f=24$, the tiling consists of two timezones for distances $4,3,2,1$ and six timezones for distance $5$.

\begin{figure}[htp]
\centering
\begin{tikzpicture}[>=latex, scale=0.7]


\foreach \x in {0,...,3}
\draw
	(-0.3+1.2*\x,1.2) -- 
	(-0.3+1.2*\x,0.7) -- 
	(-0.6+1.2*\x,0.3) -- 
	(-0.6+1.2*\x,-0.3) -- 
	(-0.9+1.2*\x,-0.7) -- 
	(-0.9+1.2*\x,-1.2)
	;
\foreach \x in {0,...,2}
\draw
	(-0.3+1.2*\x,0.7) -- 
	(1.2*\x,0.3) -- 
	(1.2*\x,-0.3) -- 
	(-0.6+1.2*\x,-0.3)
	(1.2*\x,0.3) -- ++(0.6,0)
	(1.2*\x,-0.3) -- ++(0.3,-0.4)
	;
	
\node at (0.9,0.1) {\small $*$};
\node at (1,-1.5) {\small distance 5};


\begin{scope}[xshift=5.7cm]

\foreach \x in {-1,1}
\foreach \y in {-1,1}
\draw[xscale=\x,yscale=\y]
	(0,1.2) -- (0,0.8) -- (0.4,0.5) -- (0.3,0)
	(0.4,0.5) -- (0.9,0.4) -- (1.2,0.7)
	(1.2,1.2) -- (1.2,0.7) -- (1.4,0)
	(0.9,0.4) -- (0.9,0);
\draw
	(-1.9,1.2) -- (-1.9,0.7) -- (-1.7,0)
	(-1.9,-1.2) -- (-1.9,-0.7) -- (-1.7,0)
	(-0.3,0) -- (0.3,0)
	(-1.4,0) -- (-1.7,0);

\node at (0,0.4) {\small $*$};
\node at (-0.3,-1.5) {\small distance 4};

\end{scope}


\begin{scope}[xshift=9.3cm]

\foreach \x in {-1,1}
\draw[scale=\x]
	(0,0) -- (0.16,0.16) -- (0.48,-0.16) -- (0.8,0.16) -- (1.12,-0.16) -- (1.44,0.16) -- (1.76,-0.16)
	(0.64,1.2) -- (0.64,0.7) -- (0.16,0.6) -- (-0.48,0.6) -- (-0.96,0.7) -- (-0.96,1.2)
	(0.16,0.6) -- (0.16,0.16)
	(0.48,-0.6) -- (0.48,-0.16)
	(0.64,0.7) -- (0.8,0.16)
	(0.96,-0.7) -- (1.12,-0.16)
	(1.76,-1.2) -- (1.76,-0.16) -- (1.44,0.16) -- (1.44,1.2)
	;

\node at (-0.15,0.3) {\small $*$};
\node at (0,-1.5) {\small distance 3};
	
\end{scope}


\begin{scope}[xshift=12.5cm]

\foreach \y in {-1,1}
\draw[yscale=\y]
	(0,0.2) -- (-0.4,0.4) -- (-0.1,0.8) -- (0.4,0.7) -- (0.4,0.4) -- cycle
	(-0.1,0.8) -- (-0.1,1.2)
	(0.4,0.4) -- (0.6,0) -- (0.9,0) -- (1,0.5)
	(0.4,0.7) -- (1,0.5) -- (1.3,0.6) -- (1.4,1.2)
	(-0.4,0.4) -- (-0.6,0) -- (-0.9,0)	
	(-1,1.2) -- (-0.9,0) -- (-1,-1.2)
	(1.7,1.2) -- (1.8,0) -- (1.7,-1.2);
\draw
	(0,0.2) -- (0,-0.2)
	(1.3,0.6) -- (1.3,-0.6);
	
\node at (0.3,0) {\small $*$};
\node at (0.3,-1.5) {\small distance 2};
	
\end{scope}


\begin{scope}[xshift=16.2cm]

\foreach \x in {-1,1}
\foreach \y in {-1,1}
\draw[xscale=\x,yscale=\y]
	(0,0) -- (0,0.2) -- (0.4,0.4) -- (0.6,0) -- (0.9,0) -- (1,0.8) -- (1.1,1.2)
	(0.4,0.4) -- (0.4,0.7) 
	(0,0.7) -- (0.4,0.7) -- (1,0.8)
	(-1.4,-1.2) -- (-1.4,1.2)
	(1.4,-1.2) -- (1.4,1.2);

\node at (0.3,0) {\small $*$};
\node at (0,-1.5) {\small distance 1};

\end{scope}

\end{tikzpicture}
\caption{Timezones for earth map tilings.}
\label{timezone}
\end{figure}

Next we carry out a ``propagation'' argument. The idea is to ask which of the five tiles around $T_1$ in the first of Figure \ref{case55A} can be $3^5$-tiles. If one such ``neighbouring'' tile is still a $3^5$-tile, then its neighbourhood is again given by the first of Figure \ref{case55A}. To see whether this is possible, we simplify the presentation of the neighbourhood tiling in the first of Figure \ref{case55A} by keeping only $\gamma$ and indicating the orientations of tiles according to Figure \ref{pentagon}. This is the first of Figure \ref{casespecialC}, in which the other angles can be recovered by $\gamma$ and the orientation. The second of Figure \ref{casespecialC} is the horizontal flip of the first picture.

\begin{figure}[htp]
\centering
\begin{tikzpicture}[>=latex,scale=1]



\foreach \a in {0,2}
\fill[gray!50,xshift=3.2*\a cm]
	(54:1.5) -- (90:1.15) -- (126:1.5) -- (162:1.15) -- (162:0.6) -- (90:0.6) -- (18:0.6) -- (-54:0.6) -- (-54:1.15) -- (-18:1.5) -- (18:1.15) -- cycle;

\fill[gray!50, xshift=3.2 cm]
	(162:0.6) -- (90:0.6) -- (18:0.6) -- (18:1.15) -- (-18:1.5) -- (-54:1.15) -- (-90:1.5)  -- (234:1.15) -- (234:0.6) -- cycle;
	
\foreach \a in {0,...,2}
{
\begin{scope}[xshift=3.2*\a cm]

\coordinate (P\a) at (0,0);

\foreach \x in {1,...,5}
{
\draw[rotate=72*\x]
	(90:0.6) -- (18:0.6) -- (18:1.15) -- (54:1.5) -- (90:1.15);
	
\coordinate (X\x A\a) at (18+72*\x:0.4);
\coordinate (Y\x A\a) at (6+72*\x:0.7);
\coordinate (Z\x A\a) at (30+72*\x:0.7);
\coordinate (U\x A\a) at (9+72*\x:1);
\coordinate (V\x A\a) at (27+72*\x:1);
\coordinate (W\x A\a) at (-18+72*\x:1.3);

\coordinate (P\x A\a) at (-18+72*\x:0.9);
}

\end{scope}
}


\node at (X1A0) {\small $\gamma$};
\node at (Y3A0) {\small $\gamma$};
\node at (Z4A0) {\small $\gamma$};
\node at (U2A0) {\small $\gamma$};
\node at (V5A0) {\small $\gamma$};
\node at (W4A0) {\small $\gamma$};


\node at (X1A1) {\small $\gamma$};
\node at (Y3A1) {\small $\gamma$};
\node at (Z4A1) {\small $\gamma$};
\node at (U2A1) {\small $\gamma$};
\node at (V5A1) {\small $\gamma$};
\node at (W4A1) {\small $\gamma$};



\node at (P1A2) {\small $>\!3$};
\node at (P2A2) {\small $3^5$};
\node at (P3A2) {\small $>\!3$};
\node at (P4A2) {\small $>\!3$};
\node at (P5A2) {\small $3^5$};

\node at (P2) {\small $3^5$};

\node[draw,shape=circle, inner sep=0.5] at (0,0) {\small 1};

\foreach \x in {2,...,6}
\node[draw,shape=circle, inner sep=0.5] at (198-72*\x:0.9) {\small \x};

\end{tikzpicture}
\caption{Propagation of neighbourhood tiling.}
\label{casespecialC}
\end{figure}

The positively oriented tile $T_4$ in the first of Figure \ref{casespecialC} has neighboring tiles $T_3,T_1,T_5$. We ask whether it is possible to embed $T_4,T_3,T_1,T_5$ into the first (because $T_4$ is positively oriented) of Figure \ref{casespecialC}, such that $T_4$ is the center tile $T_1$. We find the identification of $T_4$ with $T_1$ that keeps $\gamma$ and the orientation exchanges $T_3,T_5$. Since $T_3,T_5$ have different orientations, we get a contradiction. This means $T_4$ in the first picture cannot be a $3^5$-tile. We indicate this by labelling $T_4$ as $>3$ in the third picture. 

The negatively oriented tile $T_3$ in the first of Figure \ref{casespecialC} has neighboring tiles $T_2,T_1,T_4$. We ask whether it is possible to embed $T_3,T_2,T_1,T_4$ into the second (because $T_3$ is negatively oriented) of Figure \ref{casespecialC}, such that $T_3$ is the center tile. Indeed, this can be done by taking $T_3,T_2,T_1,T_4$ in the first picture to $T_1,T_3,T_2,T_6$ in the second picture, with matching $\gamma$ and orientation. This means $T_3$ in the first picture can be a $3^5$-tile. We indicate this by labelling $T_3$ as $3^5$ in the third picture. 

Similar argument determines the nature of all the tiles. We get the corresponding labels in the third of Figure \ref{casespecialC}. 

We apply the propagation to the $*$-labeled $3^5$-tiles in Figure \ref{timezone}. For distances $4,3,2,1$, the $*$-labeled tiles have at least three neighbouring $3^5$-tiles. Since the third of Figure \ref{casespecialC} has at most two neighbouring $3^5$-tiles, it cannot be the neighbourhoods of these $*$-labeled tiles. For distance $5$, we note that only the two tiles on the left and right of the $*$-labeled tile are $3^5$-tiles. These two must be the two neighbouring $3^5$-tiles in the third of Figure \ref{casespecialC}. Guided by this observation, it is easy to derive the unique earth map tiling of distance $5$ in Figure \ref{earth} (one more timezone needed). In particular, we find that the full AVC is actually $\{\alpha\beta\gamma,\delta\epsilon^2,\delta^6\}$. 

Note that $\beta\delta^2$ is not a vertex in the tiling, contradicting the original assumption about Case 5.5. In fact, the tiling belongs to the exceptional case of $\text{AVC}_3=\{\alpha\beta\gamma,\delta\epsilon^2\}$ and $v_4=v_5=0$. Since we will rediscover the tiling anyway (although starting from different premise), we still justify the existence of the equilateral pentagon here. Moreover, we also need to justify the precise values of $\beta,\delta,\epsilon$.

We know the approximate value $\cos a=0.70688$. In fact, we obtain the approximate value from the precise value
\begin{equation}\label{earth24A}
\cos a=\tfrac{1}{3}(-1+\sqrt{3}+{\textstyle \sqrt{-5+4\sqrt{3}}}).
\end{equation}
We may also recover the precise value from the $LMN$-equation \eqref{lmneq} for $AC$ 
\begin{align*}
0
&=(1-\cos\tfrac{4}{3}\pi)(1-\cos\tfrac{1}{3}\pi)\cos^2a \\
&\quad +(\cos\tfrac{5}{6}\pi+\cos(\tfrac{4}{3}+\tfrac{1}{3})\pi-\cos\tfrac{4}{3}\pi-\cos\tfrac{1}{3}\pi)\cos a \\
&\quad +(\cos\tfrac{5}{6}\pi-\sin\tfrac{4}{3}\pi\sin\tfrac{1}{3}\pi),
\end{align*}
where the precise values of $\beta,\delta,\epsilon$ are used. The root of the quadratic equation compatible with $\cos a=0.70688$ is \eqref{earth24A}.

We construct $\triangle ACE$ using $a$ from \eqref{earth24A} and $\epsilon=\frac{5}{6}\pi$. We also construct $\square ABDC$ using this $a$ and $\beta=\frac{4}{3}\pi,\delta=\frac{1}{3}\pi$. The validity of the $LMN$-equation for $AC$ means that the triangle and the quadrilateral have matching $AC$ edge. Therefore they can be glued together to form a pentagon. 

Substituting the value of $\cos a$ into the $LMN$-equation for $BC$
\begin{align*}
0
&=(1-\cos\tfrac{5}{6}\pi)(1-\cos\alpha)\cos^2a \\
&\quad +(\cos\tfrac{1}{3}\pi+\cos(\tfrac{5}{6}\pi+\alpha)-\cos\tfrac{5}{6}\pi-\cos\alpha)\cos a \\
&\quad +(\cos\tfrac{1}{3}\pi-\sin\tfrac{5}{6}\pi\sin\alpha),
\end{align*}
we get a linear equation relating $\cos\alpha$ and $\sin\alpha$
\begin{align*}
&(7+6\sqrt{3}+(8+5\sqrt{3}){\textstyle \sqrt{-5+4\sqrt{3}}})\cos\alpha +3(2+\sqrt {3}+{\textstyle \sqrt{-5+4\sqrt{3}}})\sin\alpha \\
&= 19+3\sqrt{3}+5(1+\sqrt{3}){\textstyle \sqrt{-5+4\sqrt{3}}}.
\end{align*}
The equation has two solutions, and the one compatible with the approximate value $\alpha=0.14400\pi$ is
\[
\alpha
= \arctan\tfrac{1}{33}(4+3\sqrt{3}+(-2+4\sqrt{3}){\textstyle \sqrt{-5+4\sqrt{3}}}).
\]

Similarly, the $LMN$-equation for $AD$ is a linear equation relating $\cos\gamma$ and $\sin\gamma$
\begin{align*}
& (7+6\sqrt{3}+(8+5\sqrt{3}){\textstyle \sqrt{-5+4\sqrt{3}}} )\cos\gamma 
+3(2+\sqrt {3}+{\textstyle \sqrt{-5+4\sqrt{3}}})\sin\gamma \\
&= 7-3\sqrt{3}+(-1+5\sqrt{3}){\textstyle \sqrt{-5+4\sqrt{3}}}.
\end{align*}
The solution compatible with $\gamma=0.52265\pi$ is
\[
\gamma
= \pi-\arctan\tfrac{1}{3}(12+7\sqrt{3}+(6+4\sqrt{3}){\textstyle \sqrt{-5+4\sqrt{3}}}).
\]

We may symbolically verify
\[
\tan(\alpha+\gamma)
=\frac{\tan\alpha+\tan\gamma}{1+\tan\alpha \tan\gamma}
=-\sqrt{3}.
\]
The only precise value of $\alpha+\gamma$ compatible with $\alpha=0.14400\pi$ and $\gamma=0.52265\pi$ is $\frac{2}{3}\pi$. Therefore $\alpha+\beta+\gamma=2\pi$. 

We justify the simple property of the pentagon in Section \ref{case15b}.

\subsection{Cases 1.2e, 1.5a, 2.4b}
\label{case12e}

We start with the version of the solution of the cases given at the beginning of Section \ref{tiling}. This means we apply the exchange $\alpha\leftrightarrow\gamma$ to the original solution for $A11$, and get the new version for $A3$. By the method in Section \ref{case42b}, we get
\[
\text{AVC}=
\{\alpha\beta\gamma,
\delta\epsilon^2,
\gamma\delta^2\epsilon,
\gamma^4,
\gamma^2\delta^3,
\delta^6\}.
\]
By $\beta\cdots=\alpha\beta\gamma$, and the total number of times $\beta$ appears in the tiling is $f$, we find that $\alpha\beta\gamma$ appears $f$ times. This implies that $\gamma$ already appears $f$ times at $\alpha\beta\gamma$. Since the total number of times $\gamma$ appears in the tiling is also $f$, we conclude $\gamma\cdots=\alpha\beta\gamma$. Therefore $\gamma\delta^2\epsilon,\gamma^4,\gamma^2\delta^3$ are not vertices, and we get the updated AVC
\[
\text{AVC}=
\{\alpha\beta\gamma,
\delta\epsilon^2,
\delta^6\}.
\]

The AVC shows that we are in none of Cases 1.2e, 1.5a, 2.4b. Still, we may ask whether the reduced AVC admits tiling. Since this AVC is contained in the AVC studied in Section \ref{case55}, the tiling is the earth map tiling obtained in Section \ref{case55}, i.e., given by Figure \ref{earth}. The difference is the pentagon used for tiling. 

We need to justify the existence of the pentagon, including the precise values of $\gamma,\delta,\epsilon$. We note that, by $\alpha+\beta+\gamma=2\pi$ and $\gamma=\frac{1}{2}\pi$, we have $\gamma-\alpha=\beta-\pi$. This means the pentagon is obtained by directly glueing $\triangle ACE$ and $\triangle BCD$ together. In other words, the triangle $\triangle ABC$ is reduced to an arc. See Figure \ref{case2a}.

For $0\le a\le \frac{1}{2}\pi$, therefore, we construct isosceles triangles $\triangle ACE,\triangle BCD$ with the same side length $a$ and the given top angles $\delta=\frac{1}{3}\pi,\epsilon=\frac{5}{6}\pi$. Then the two base angles $\alpha$ and $\theta$ are strictly increasing and continuous functions of $a$, with $\alpha+\theta=\frac{1}{12}\pi+\frac{1}{3}\pi=\frac{5}{12}\pi$ at $a=0$, and $\alpha+\theta=\frac{1}{2}\pi+\frac{1}{2}\pi=\pi$ at $a=\frac{1}{2}\pi$. By the intermediate value theorem, we have $\alpha+\theta=\frac{1}{2}\pi$ for a unique $a$ satisfying $0< a< \frac{1}{2}\pi$. The pentagon $\pentagon ABDCE$ obtained by glueing the two isosceles triangles has $AE=EC=CD=DB=a$, and $\gamma=\alpha+\theta=\frac{1}{2}\pi$, $\delta=\frac{1}{3}\pi$, $\epsilon=\frac{5}{6}\pi$. It remains to show that $AB=AC-BC=a$, so that the pentagon is equilateral.

\begin{figure}[htp]
\centering
\begin{tikzpicture}[>=latex,scale=1]

\draw
	(0,0) arc (150:30:4)
	(0,0) arc (-150:-30:4)
	(2,1.72) node[above] {\small $C$} -- 
	(4,1.02) node[right] {\small $D$} -- 
	(2,0.32) node[left] {\small $B$} -- 
	(2,-1.72) node[below] {\small $A$}
	;

\draw[dashed]
	(2,1.73) -- (2,0.32);

\node at (7.2,0) {\small $E'$};
\node at (-0.2,0) {\small $E$};

\node at (0.4,0) {\small $\frac{5}{6}\pi$};
\node at (2.9,1.72) {\small $\frac{1}{2}\pi$};
\node at (3.2,1.02) {\small $\frac{1}{3}\pi$};
\node at (1.8,-1.4) {\small $\alpha$};
\node at (1.8,1.4) {\small $\alpha$};
\node at (2.2,1.4) {\small $\theta$};
\node at (2.2,0.6) {\small $\theta$};

\node at (0,1.7) {\small $\theta+\alpha=\frac{1}{2}\pi$};

\end{tikzpicture}
\caption{Pentagon for 
$\{
\alpha\beta\gamma,\delta\epsilon^2,\gamma\delta^2\epsilon\text{ or }\gamma^4\text{ or }\gamma^2\delta^3\}$, $A3$ arrangement.}
\label{case2a}
\end{figure}

We have
\begin{equation}\label{earth24C}
\tan\theta=\sec a\cot\tfrac{1}{2}\delta,\quad
\tan\alpha=\sec a\cot\tfrac{1}{2}\epsilon.
\end{equation}
By $\alpha+\theta=\gamma=\frac{1}{2}\pi$, this implies 
\[
\sec a\cot\tfrac{1}{2}\delta\cdot \sec a\cot\tfrac{1}{2}\epsilon=1.
\]
Substituting the precise values of $\delta$ and $\epsilon$, we get
\begin{equation}\label{earth24B}
\cos a={\textstyle \sqrt{-3+2\sqrt{3}}},\quad
\sin a={\textstyle \sqrt{1-(-3+2\sqrt{3})}}=-1+\sqrt{3}.
\end{equation}
Then
\begin{align*}
\cos AC
&=\cos^2a+\sin^2a\cos\tfrac{5}{6}\pi \\
&=(-3+2\sqrt{3})-(4-2\sqrt{3})\tfrac{\sqrt{3}}{2}
=0, \\
\cos BC
&=\cos^2a+\sin^2a\cos\tfrac{1}{3}\pi \\
&=(-3+2\sqrt{3})+(4-2\sqrt{3})\tfrac{1}{2}
=-1+\sqrt{3}
=\sin a.
\end{align*}
The first equality implies $AC=\frac{1}{2}\pi$. The second equality implies $\cos BC>0$, so that $0<BC<\frac{1}{2}\pi$. Since we also have $0<a<\frac{1}{2}\pi$, the second equality further implies $BC+a=\frac{1}{2}\pi=AC$. Moreover, we may use \eqref{earth24C} to further derive
\begin{align*}
\alpha
&=\tfrac{1}{2}\pi-\theta=\tfrac{1}{2}\pi-\arctan{\textstyle \sqrt{3+2\sqrt{3}}},\\
\beta
&=\pi+\theta=\pi+\arctan{\textstyle \sqrt{3+2\sqrt{3}}}.
\end{align*}

\subsection{Cases 1.4e, 2.6b}
\label{case14e}

By the method in Section \ref{case42b}, for both solutions, we get
\[
\text{AVC}
=\{\alpha\beta\gamma,\delta\epsilon^2,\delta^3\epsilon,\delta^5\}.
\]
We need to consider both arrangements $A1$ and $A3$.

For the $A1$ arrangement, we have consecutive $\delta\delta\delta$ at $\delta^3\epsilon,\delta^5$. The AAD of $\delta\delta\delta=\thin\delta\thin^{\gamma}\delta^{\epsilon}\thin\delta\thin$ implies a vertex $\gamma^2\cdots$ or $\gamma\epsilon\cdots$, a contradiction. 

For the $A3$ arrangement, we consider tiles $T_1,T_2,T_3$ containing consecutive $\delta\delta\delta$, in Figure \ref{tiling15b}. We may assume $T_2$ is arranged as indicated. By $\delta_1$, we get $\gamma_2\cdots=\alpha_4\beta_1\gamma_2$. Then $\beta_1,\delta_1$ determine $T_1$. By $\alpha_4$, we get $\epsilon_2\cdots=\delta_5\epsilon_2\epsilon_4$. Then $\alpha_4,\epsilon_4$ determine $T_4$. By $\delta_5$, we get $\gamma_4\cdots=\alpha_6\beta_5\gamma_4$. Then $\beta_5,\delta_5$ determine $T_5$. By $\alpha_6$ and no $\alpha_5\epsilon\cdots$, we determine $T_6$.

\begin{figure}[htp]
\centering
\begin{tikzpicture}[>=latex,scale=1]

\foreach \a in {0,1,-1}
\draw[xshift=2*\a cm]
	(0,1.3) -- (0,0.8) -- (0.5,0.4) -- (1.5,0.4) -- (2,0.8) -- (2,1.3);

\draw
	(0.5,0.4) -- (0.5,-0.4) -- (1,-0.8)
	(0.5,-0.4) -- (-0.5,-0.4)
	(-0.5,0.4) -- (-0.5,-0.4) -- (-1,-0.8) -- (-1,-1.3)
	(1.5,0.4) -- (1.5,-0.4) -- (1,-0.8) -- (1,-1.3);
			
\foreach \a in {0,-1}
{
\begin{scope}[xshift=2*\a cm]

\node at (1.5,0.6) {\small $\alpha$};  
\node at (1.85,0.9) {\small $\beta$};
\node at (0.15,0.95) {\small $\gamma$};
\node at (1,1.1) {\small $\delta$};
\node at (0.5,0.6) {\small $\epsilon$};

\end{scope}
}

\node at (3,1.1) {\small $\delta$};

\node at (0,0.55) {\small $\alpha$}; 
\node at (-0.35,0.3) {\small $\beta$};
\node at (0.35,-0.2) {\small $\gamma$};
\node at (-0.35,-0.2) {\small $\delta$};
\node at (0.35,0.3) {\small $\epsilon$};

\node at (1,-0.55) {\small $\alpha$}; 
\node at (0.65,-0.3) {\small $\beta$};
\node at (1.35,0.2) {\small $\gamma$};
\node at (0.65,0.2) {\small $\delta$};
\node at (1.35,-0.3) {\small $\epsilon$};	

\node at (0.5,-0.6) {\small $\alpha$}; 
\node at (0.85,-0.95) {\small $\beta$}; 
\node at (-0.85,-0.9) {\small $\gamma$}; 
\node at (0,-1.1) {\small $\delta$};
\node at (-0.5,-0.6) {\small $\epsilon$};

\node[inner sep=1,draw,shape=circle] at (-1,0.7) {\small 1};
\node[inner sep=1,draw,shape=circle] at (1,0.7) {\small 2};
\node[inner sep=1,draw,shape=circle] at (3,0.7) {\small 3};
\node[inner sep=1,draw,shape=circle] at (0,0.1) {\small 4};
\node[inner sep=1,draw,shape=circle] at (1,-0.1) {\small 5};
\node[inner sep=1,draw,shape=circle] at (0,-0.7) {\small 6};

\end{tikzpicture}
\caption{Tiling for 
$\{\alpha\beta\gamma,\delta\epsilon^2,\delta^3\epsilon,\delta^5\}$.}
\label{tiling15b}
\end{figure}

We find that $T_2,T_4,T_5,T_6$ form a timezone in the earth map tiling in Figure \ref{earth}. We derive the timezone from $T_1,T_2$. Similarly, for $\delta_3$ on the right of $\delta_2$, by $\delta_3$, we get $\beta_2\cdots=\alpha\beta_2\gamma_3$. This determines $T_3$, and $T_2,T_3$ is the same as $T_1,T_2$. Therefore we get another timezone containing $T_3$. If there are more $\delta$ to the right of $\delta_3$, then the argument continues to give more timezones. If $\delta_1\delta_2\cdots=\delta^5$, then the repetition of the argument gives the earth map tiling in Figure \ref{earth}.

We note that the full AVC for the earth map tiling is $\{
\alpha\beta\gamma,\delta\epsilon^2,\delta^5\}$. Therefore the tiling is actually for Case 2.6b, and not for Case 1.4e.

It remains to consider the $A3$ arrangement, and $\delta^5$ is not a vertex. This means  
\[
\text{AVC}=\{\alpha\beta\gamma,\delta\epsilon^2,\delta^3\epsilon\}.
\]
In other words, we study the tiling for Case 1.4e and not for Case 2.6b. We will show that the tiling is given by Figure \ref{tiling14f}. The tiling has two tiles with two degree $4$ vertices and three degree $3$ vertices, which we call {\em $3^34^2$-tiles}. We draw the two tiles as the north and south ``regions'' (as opposed to the two poles in the earth map tiling) $T_1,T_{10}$. The tiling is obtained by glueing the left and right together.

\begin{figure}[htp]
\centering
\begin{tikzpicture}[>=latex,scale=1.2]

\foreach \a in {1,-1}
\draw[gray!30,line width=3,scale=\a]
	(-0.8,1.5) -- (-1.2,0.7) -- (-1.6,-0.3) -- (-2.4,-0.9) -- (-2.4,-1.5) -- (-4,-1.5)
	(-2.4,1.5) -- (-2.4,0.9) -- (-3.1,0.3) -- (-3.5,-0.7) -- (-4,-1.5);

\foreach \a in {1,-1}
\draw[gray!70,line width=5,scale=\a]
 	(-0.3,0.4) -- (0.3,-0.4) -- (0.8,-1.5) -- (4,-1.5) -- (3.7,-0.6) -- (4.3,0.6) -- (4,1.5)
	;
	
\foreach \a in {1,-1}
{
\begin{scope}[scale=\a]
	
\draw
	(-4.8,1.5) -- (4.8,1.5)
	(-0.8,1.5) -- (-0.3,0.4) 
	(0.8,1.5) -- (0.8,0.6) -- (0.3,-0.4) 
	(2.4,1.5) -- (2.4,0.9) -- (2.5,0.3) 
	(0.8,0.6) -- (1.6,0.3) -- (2.4,0.9)
	(1.6,0.3) -- (1.2,-0.7) -- (2.2,-0.3)
	(4,1.5) -- (3.5,0.7) -- (2.5,0.3) 
	(4,1.5) -- (4.3,0.6) -- (3.7,-0.6)
	(0.3,-0.4) -- (-0.3,0.4) 
	(2.5,0.3) -- (2.2,-0.3) 
	(1.2,-0.7) -- (0.8,-1.5) 
	(2.4,-1.5) -- (2.4,-0.9) -- (2.2,-0.3) 
	(4,-1.5) -- (3.7,-0.6) -- (3.1,-0.3) 
	(2.4,-0.9) -- (3.1,-0.3) -- (3.5,0.7)
	;

\draw[dashed]
	(4.3,0.6) -- (4.9,0.3) -- (4.5,-0.7) -- (4,-1.5)
	(4.5,-0.7) -- (5.5,-0.3);	


\node at (-4,1.7) {\small $\delta$};
\node[gray] at (-4.3,1.3) {\small $\delta$};
\node[gray] at (-4.05,1.15) {\small $\delta$};
\node at (-3.8,1.35) {\small $\epsilon$};

\node at (-2.4,1.7) {\small $\gamma$};
\node at (-2.55,1.3) {\small $\alpha$};
\node at (-2.25,1.3) {\small $\beta$};

\node at (-0.8,1.7) {\small $\epsilon$};
\node at (-0.55,1.3) {\small $\delta$};
\node at (-1.05,1.3) {\small $\delta$};
\node at (-0.83,1.15) {\small $\delta$};

\node at (0.8,1.7) {\small $\alpha$};
\node at (0.65,1.3) {\small $\gamma$};
\node at (1,1.3) {\small $\beta$};

\node at (2.4,1.7) {\small $\beta$};
\node at (2.25,1.3) {\small $\alpha$};
\node at (2.55,1.3) {\small $\gamma$};

\node at (4,1.7) {\small $\delta$};
\node at (3.7,1.3) {\small $\delta$};
\node at (3.95,1.15) {\small $\delta$};
\node[gray] at (4.2,1.35) {\small $\epsilon$};


\node at (-3.55,0.75) {\small $\gamma$};
\node at (-3.6,0.35) {\small $\alpha$};

\node at (-3.1,0.5) {\small $\delta$};
\node at (-3,0.2) {\small $\epsilon$};
\node at (-3.3,0.2) {\small $\epsilon$};

\node at (-2.2,0.9) {\small $\alpha$};
\node at (-2.6,0.95) {\small $\beta$};
\node at (-2.45,0.65) {\small $\gamma$};

\node at (-2.1,0.5) {\small $\epsilon$};
\node at (-2.1,0.15) {\small $\epsilon$};
\node at (-2.4,0.3) {\small $\delta$};

\node at (-1.1,0.6) {\small $\beta$};
\node at (-1.3,0.85) {\small $\gamma$};
\node at (-1.5,0.4) {\small $\alpha$};

\node at (-0.1,0.5) {\small $\beta$};
\node at (-0.5,0.4) {\small $\gamma$};
\node at (-0.25,0.1) {\small $\alpha$};

\node at (0.95,0.75) {\small $\delta$};
\node at (0.65,0.6) {\small $\epsilon$};
\node at (0.85,0.4) {\small $\epsilon$};

\node at (1.6,0.5) {\small $\gamma$};
\node at (1.7,0.2) {\small $\beta$};
\node at (1.4,0.2) {\small $\alpha$};

\node at (2.55,0.9) {\small $\epsilon$};
\node at (2.25,1) {\small $\epsilon$};
\node at (2.3,0.65) {\small $\delta$};

\node at (2.6,0.5) {\small $\alpha$};
\node at (2.6,0.15) {\small $\beta$};
\node at (2.3,0.3) {\small $\gamma$};

\node at (3.65,0.6) {\small $\gamma$};
\node at (3.4,0.85) {\small $\beta$};
\node at (3.2,0.4) {\small $\alpha$};

\node at (4.1,0.6) {\small $\beta$};

\end{scope}
}

\node at (0,2.2) {\small north tile (region)};
\node at (0,-2.2) {\small south tile (region)};

\node[inner sep=0.5,draw,shape=circle] at (0,1.8) {\small $1$};
\node[inner sep=0.5,draw,shape=circle] at (-1.7,1) {\small $2$};
\node[inner sep=0.5,draw,shape=circle] at (-0.9,0.2) {\small $3$};
\node[inner sep=0.5,draw,shape=circle] at (0.2,1) {\small $4$};
\node[inner sep=0,draw,shape=circle] at (-3.1,1) {\small $17$};
\node[inner sep=0,draw,shape=circle,gray] at (4.7,1) {\small $17$};
\node[inner sep=0.5,draw,shape=circle] at (-2,-0.15) {\small $6$};
\node[inner sep=0.5,draw,shape=circle] at (-0.2,-1) {\small $9$};
\node[inner sep=0.5,draw,shape=circle] at (-2.7,0.2) {\small $5$};
\node[inner sep=0.5,draw,shape=circle] at (-1.6,-1) {\small $8$};
\node[inner sep=0,draw,shape=circle] at (0.9,-0.2) {\small $12$};
\node[inner sep=0,draw,shape=circle] at (0,-1.8) {\small $10$};
\node[inner sep=0,draw,shape=circle] at (1.6,1) {\small $11$};
\node[inner sep=0,draw,shape=circle] at (1.7,-1) {\small $13$};
\node[inner sep=0,draw,shape=circle] at (2,0.15) {\small $14$};
\node[inner sep=0,draw,shape=circle] at (2.7,-0.2) {\small $15$};
\node[inner sep=0,draw,shape=circle] at (3,1) {\small $16$};
\node[inner sep=0,draw,shape=circle] at (3.7,0.1) {\small $18$};
\node[inner sep=0,draw,shape=circle,gray] at (-4.3,0.1) {\small $18$};
\node[inner sep=0,draw,shape=circle] at (3.1,-1) {\small $19$};
\node[inner sep=0,draw,shape=circle] at (-3.7,-0.1) {\small $20$};
\node[inner sep=0,draw,shape=circle,gray] at (4.3,-0.1) {\small $20$};
\node[inner sep=0.5,draw,shape=circle] at (-3,-1) {\small $7$};

\end{tikzpicture}
\caption{Tiling for 
$\{ \alpha\beta\gamma,\delta\epsilon^2,\delta^3\epsilon\}$.}
\label{tiling14f}
\end{figure}

We start with four tiles $T_1,T_2,T_3,T_4$ around a vertex $\delta^3\epsilon=\delta_2\delta_3\delta_4\epsilon_1$. By no $\gamma^2\cdots$, we get the unique AAD $\thin^{\beta}\delta^{\gamma}\thin^{\beta}\delta^{\gamma}\thin^{\beta}\delta^{\gamma}\thin^{\alpha}\epsilon^{\gamma}\thin$ of the vertex. This determines $T_1,T_2,T_3,T_4$. We note that the current AVC $\{\alpha\beta\gamma,\delta\epsilon^2,\delta^3\epsilon\}$ is a subset of the AVC $\{\alpha\beta\gamma,\delta\epsilon^2,\delta^3\epsilon,\delta^5\}$ used for deriving Figure \ref{tiling15b}. Therefore we take the pairs $T_3,T_2$ and $T_4,T_3$ as the pair $T_1,T_2$ in Figure \ref{tiling15b}, and get two timezones originating from $T_2$ and $T_3$. The timezone originating from $T_2$ consists of $T_2,T_5,T_6,T_7$. The timezone originating from $T_3$ consists of $T_3,T_8,T_9,T_{10}$.

We have $\epsilon_4\cdots=\delta\epsilon^2$ or $\delta^3\epsilon$. We continue the construction of the tiling by assuming $\epsilon_4\cdots=\delta\epsilon^2$.  By $\alpha_1\gamma_4\cdots=\alpha_1\beta_{11}\gamma_4$ and $\alpha_4\beta_9\cdots=\alpha_4\beta_9\gamma_{12}$, we get $\epsilon_4\cdots=\delta_{11}\epsilon_4\epsilon_{12}$. Then $\beta_{11},\delta_{11},\gamma_{12},\epsilon_{12}$ determine $T_{11},T_{12}$. Then $\delta_9\delta_{12}\epsilon_{10}
\cdots=\delta^3\epsilon$. We may apply the argument starting with the original $\delta^3\epsilon$ to this new $\delta^3\epsilon$, and get $T_{13},T_{14},T_{15},T_{16}$. 

By $\beta_2\gamma_1\cdots=\alpha_{17}\beta_2\gamma_1$ and no $\beta_{17}\delta_1\cdots$, we determine $T_{17}$. Then the unique AAD of $\delta_1\delta_{16}\epsilon_{17}\cdots=\delta_1\delta_{16}\delta_{18}\epsilon_{17}$ determines $T_{18}$. By the same argument, we also determine $T_{19},T_{20}$.

We obtain the tiling by assuming $\epsilon_4\cdots=\delta\epsilon^2$. The alternative assumption is that $\epsilon_4\cdots=\delta^3\epsilon$. This means that $T_4$ is a $3^34^2$-tile, with $\delta_4\cdots=\epsilon_4\cdots=\delta^3\epsilon$. Therefore by taking this $T_4$ as $T_1$, we may assume that the $T_1$ we started with is already a $3^34^2$-tile, and $\delta_1\cdots=\epsilon_1\cdots=\delta^3\epsilon$. Now we repeat the argument with this additional assumption. We still get $T_1,\dots,T_{10}$. We may also repeat the same argument by starting at $\delta_1\cdots=\delta^3\epsilon$. This means $\epsilon_{17},\delta_{18},\delta_{16},\delta_1$ correspond to $\epsilon_1,\delta_2,\delta_3,\delta_4$. The unique AAD of $\delta^3\epsilon$ determines $T_{17}$. Moreover, we get two timezones originating from $\delta_{16},\delta_{18}$. The timezone originating from $\delta_{16}$ consists of $T_{11},T_{12},T_{14},T_{16}$, The timezone originating from $\delta_{18}$ consists of $T_{13},T_{15},T_{18},T_{19}$. Then the unique AAD of $\delta_7\delta_{10}\epsilon_{19}\cdots=\delta_7\delta_{10}\delta_{20}\epsilon_{19}$ determines $T_{20}$.

The dark and thickly shaded lines divide the tiling in Figure \ref{tiling14f} into two equal halves. Each half has two timezones and two additional tiles, and can be identified with the halves in Figure \ref{earth} between the thick lines. By looking at the boundary vertices along the thickly shaded lines in Figure \ref{tiling14f}, we find the tiling is the flip modification described by Figure \ref{flipmod}.

We still need to justify the existence of the pentagon. We leave the discussion to the next section.

\subsection{Cases 1.5b, 2.5e}
\label{case15b}

By the method in Section \ref{case42b}, we get 
\[
\text{AVC}
=\{\alpha\beta\gamma,\delta\epsilon^2,\delta^4\}
\]
for (all three solutions of) Case 1.5b, and get
\[
\text{AVC}
=\{\alpha\beta\gamma,\delta\epsilon^2,\delta^4\epsilon,\delta^7\}
\]
for Case 2.5e. We need to consider $A1$ and $A3$ for Case 1.5b, and consider $A1$ for Case 2.5e.

For the $A1$ arrangement, we have consecutive $\delta\delta\delta$ at $\delta^4,\delta^4\epsilon,\delta^7$. Similar to the earlier Cases 1.4e and 2.6b, the AAD of $\delta\delta\delta=\thin\delta\thin^{\gamma}\delta^{\epsilon}\thin\delta\thin$ implies a vertex $\gamma^2\cdots$ or $\gamma\epsilon\cdots$, a contradiction. 

For the $A3$ arrangement, we may also use Figure \ref{tiling15b} to prove that the $A3$ arrangement of Case 1.5b gives the earth map tiling with four timezones.

We obtained four earth map tilings, all with the $A3$ arrangement:
\begin{enumerate}
\item Case 1.5b: $f=16$, $\text{AVC}
=\{\alpha\beta\gamma,\delta\epsilon^2,\delta^4\}$.
\item Case 2.6b (and 1.4e): $f=20$, $\text{AVC}
=\{\alpha\beta\gamma,\delta\epsilon^2,\delta^5\}$.
\item Exceptional 1 (and 5.5): $f=24$, $\text{AVC}
=\{\alpha\beta\gamma,\delta\epsilon^2,\delta^6\}$.
\item Exceptional 2 (and 1.2e, 1.5a, 2.4b): $f=24$, $\text{AVC}
=\{\alpha\beta\gamma,\delta\epsilon^2,\delta^6\}$.
\end{enumerate}
Note that Case 1.4e has the same pentagon as 2.6b, but with the flip modification as the tiling. Moreover, Case 5.5 has the same pentagon as the first exceptional case, but with a vertex $\beta\delta^2$ that does not show up in the tiling. The same happens to Cases 1.2e, 1.5a, 2.4b, compared with the second exceptional case. Therefore, although we calculated the pentagons for the two exceptional cases in Sections \ref{case55} and \ref{case12e}, the calculations were based on the premise not immediately implied by the earth map tiling.   

The rest of the section is devoted to the calculation of the pentagons for the four cases, starting from the premise of earth map tiling. We also justify the existence of the pentagons.

The first of Figure \ref{emt_calculate} describes part of the earth map tiling. We know $\alpha\beta\gamma,\delta\epsilon^2,\delta^n$ are vertices, and the poles $N$ and $S$ are $\delta^n$. We know the precise values $\delta=\frac{2}{n}\pi$, $\epsilon=\left(1-\frac{1}{n}\right)\pi$, and all the normal lines have the same length $a$. Then both $AB$ and $AC$ can be reached by combining three segments of length $a$ at alternating angles $\beta$ and $\epsilon$. Therefore $AB$ and $AC$ have the same length. This implies $NAS$ is a great arc connecting the two poles. Then we get $NA=\pi-a$, and the second of Figure \ref{emt_calculate}. We may calculate $AB$ by three $a$ at alternating angles $\beta,\epsilon$, or by the triangle $\triangle ABN$
\begin{align*}
\cos AB
&=(1-\cos\beta)(1-\cos\epsilon)\cos^3a
+\sin\beta\sin\epsilon\cos^2a \\
&\quad +(\cos\beta+\cos\epsilon-\cos\beta\cos\epsilon)\cos a
-\sin\beta\sin\epsilon; \\
\cos AB
&=\cos a\cos(\pi-a)+\sin a\sin(\pi-a)\cos\tfrac{1}{2}\delta \\
&=-\cos^2 a+\sin^2a\cos\tfrac{1}{2}\delta
=-(1+\cos\tfrac{1}{2}\delta)\cos^2 a+\cos\tfrac{1}{2}\delta.
\end{align*}
Identifying the two gives one equation relating $\beta,\delta,\epsilon,a$. We may further divide a factor $1+\cos a$ to get a quadratic equation of $\cos a$. On the other hand, we have the $LMN$-equation \eqref{lmneq} for $AC$
\begin{align*}
&(1-\cos\beta)(1-\cos\delta)\cos^2a \\
&+(\cos\epsilon+\cos(\beta+\delta)-\cos\beta-\cos\delta)\cos a \\
&+(\cos\epsilon-\sin\beta\sin\delta)=0.
\end{align*}
This gives another equation relating $\beta,\delta,\epsilon,a$. 

Both equations are of the form $\lambda_i\cos\beta+\mu_i\sin\beta=\nu_i$, where $\lambda_i,\mu_i,\nu_i$ are quadratic polynomials of $\cos a$ with coefficients involving only $\delta,\epsilon$. Then we may solve $\cos\beta$ and $\sin\beta$ from the system, and get the equality $\cos^2\beta+\sin^2\beta=1$, which is
\[
(\lambda_1\nu_2-\lambda_2\nu_1)^2
+(\mu_1\nu_2-\mu_2\nu_1)^2
=(\lambda_1\mu_2-\lambda_2\mu_1)^2.
\]
This is actually a polynomial of $\cos a$ of degree $6$, with coefficients determined by the precise values of $\delta,\epsilon$.

\begin{figure}[htp]
\centering
\begin{tikzpicture}[>=latex,scale=1]


\foreach \a in {-1,1}
\draw[xscale=\a]
	(-0.8,0.5) -- (0.8,0.5)
	(1.6,1.4) -- (1.6,0.8) -- (0.8,0.5) -- (0.8,-0.5) -- (0,-0.8) -- (0,-1.4);

\draw[dashed]
	(1.6,0.8) -- (0,-0.8) -- (-1.6,0.8)
	(0,1.4) -- (0,-0.8);
	
\node at (-0.75,0.65) {\small $\epsilon$};	
\node at (0.75,0.65) {\small $\alpha$};
\node at (-1.45,0.9) {\small $\gamma$};
\node at (1.45,0.9) {\small $\beta$};

\node at (0,-0.6) {\small $\alpha$};
\node at (-0.65,-0.35) {\small $\beta$};
\node at (0.65,0.3) {\small $\gamma$};
\node at (-0.65,0.3) {\small $\delta$};
\node at (0.65,-0.4) {\small $\epsilon$};

\node at (0.95,0.35) {\small $\beta$};
\node at (-0.95,0.4) {\small $\epsilon$};

\node at (0.2,-0.9) {\small $A$};
\node at (1.8,0.8) {\small $B$};
\node at (-1.8,0.8) {\small $C$};
\node at (0,1.6) {\small $N$};
\node at (0,-1.6) {\small $S$};


\begin{scope}[xshift=3cm]

\draw
	(0,-1.6) -- ++(20:1.2) -- ++(100:1.2) -- ++(30:1.2) -- ++(130:1.2);

\draw[dashed]
	(1.96,0.59) -- (0,-1.6) -- (1.19,1.51);

\node at (-0.2,-1.6) {\small $A$};
\node at (2.1,0.6) {\small $B$};
\node at (1,1.6) {\small $N$};

\node[fill=white,inner sep=0.5] at (1.6,1) {\small $a$};
\node[fill=white,inner sep=0.5] at (1.45,0.3) {\small $a$};
\node[fill=white,inner sep=0.5] at (1.05,-0.7) {\small $a$};
\node[fill=white,inner sep=0.5] at (0.6,-1.35) {\small $a$};
\node[fill=white,inner sep=0.5, rotate=70] at (0.5,0.2) {\small $\pi-a$};

\node at (1.1,-0.1) {\small $\beta$};
\node at (1,-1.1) {\small $\epsilon$};

\node at (1.2,1.15) {\small $\frac{1}{2}\delta$};

\end{scope}
	
\end{tikzpicture}
\caption{Calculate the pentagon for earth map tiling.}
\label{emt_calculate}
\end{figure}

For $f=16$ (Case 1.5b), the degree $6$ polynomial of $\cos a$ is
\[
(49t^4+(16-18\sqrt{2})t^3+(48-54\sqrt{2})t^2 
+(52-34\sqrt{2})t+43-30\sqrt{2})(t^2+3-2\sqrt{2}).
\]
The root compatible with $\cos a=0.77943$ is the biggest root of the quartic factor
\begin{align}
49\cos^4a &+(16-18\sqrt{2})\cos^3a+(48-54\sqrt{2})\cos^2a \nonumber \\
& +(52-34\sqrt{2})\cos a+43-30\sqrt{2}=0. \label{length16}
\end{align}
For $f=20$ (Case 2.6b), the degree $6$ polynomial for $\cos a$ is 
\[
(5\sqrt{5}t^3+t^2-(4+\sqrt{5})t-1)
(\sqrt{5}t^2+2+\sqrt{5})
(\sqrt{5}t-1).
\]
The root compatible with $\cos a=0.77680$ is the biggest root of the cubic factor
\begin{equation}\label{length20}
5\sqrt{5}\cos^3a+\cos^2a-(4+\sqrt{5})\cos a-1=0.
\end{equation}
After getting the precise value of $a$ from \eqref{length16} and \eqref{length20}, we may use suitable $LMN$-equations to derive the precise values of $\alpha,\beta,\gamma$ similar to Section \ref{case55}. Their approximate values are the same as the ones given near the end of Section \ref{calculation}.

For $f=24$ (two exceptional cases), the degree $6$ polynomial for $\cos a$ is 
\[
(3t^2+(2-2\sqrt{3})t+3-2\sqrt{3})
(t^2+3-2\sqrt{3})
(t^2+7-4\sqrt{3}).
\]
The roots compatible with $\cos a=0.70688$ and $0.68125$ are
\[
\cos a
=\tfrac{1}{3}(-1+\sqrt{3}+{\textstyle \sqrt{-5+4\sqrt{3}}}),\quad
\cos a={\textstyle \sqrt{-3+2\sqrt{3}}}.
\]
Then we recover the two pentagons in Sections \ref{case55} and \ref{case12e}.

Next, we justify the existence of the pentagon for $f=16$,  $f=20$, and the $\cos a=0.70688$ case of $f=24$. The main concern is that the pentagon should be simple. The simple property is clear for the $\cos a=0.68125$ case of $f=24$, by the construction of the pentagon in Section \ref{case12e}. 

For $a$ in certain interval $[a_0,a_1]\subset (0,\frac{1}{2}\pi)$, we construct isosceles triangles $\triangle ACE,\triangle BCD$ with edges $a$ and top angles $\epsilon=(1-\frac{1}{n})\pi,\delta=\frac{2}{n}\pi$. By $\delta<\epsilon$, we have $u=AC>v=BC$. If we have $u-v<a$, then we may enlarge $\angle DCE$ until $AB=a$. Then we get an equilateral pentagon (in the $A3$ arrangement) with given $\delta,\epsilon$. See Figure \ref{emt_construct}. We need to argue that, for suitable $a$, we get $\alpha+\beta+\gamma=2\pi$. This gives an equilateral pentagon, such that $\alpha\beta\gamma,\delta^2\epsilon,\delta^n$ can be vertices (or strictly speaking, the angle sums are $2\pi$).  

\begin{figure}[htp]
\centering
\begin{tikzpicture}[>=latex,scale=1]

\draw
	(0,1.2) -- (-0.9,0) -- (0,-1.2)
	(0,1.2) -- (1.2,0.3) -- (0,-0.6) -- (0,-1.2);

\draw[dashed]
	(0,1.2) -- (0,-0.6);
	
\node at (-0.7,0) {\small $\epsilon$};
\node at (0.9,0.3) {\small $\delta$};

\node at (0,-1.4) {\small $A$};
\node at (0.2,-0.7) {\small $B$};
\node at (0,1.4) {\small $C$};
\node at (1.4,0.3) {\small $D$};
\node at (-1.1,0) {\small $E$};

\node[fill=white, inner sep=1.5] at (-0.45,0.6) {\small $a$};
\node[fill=white, inner sep=1.5] at (-0.45,-0.6) {\small $a$};
\node[fill=white, inner sep=1.5] at (0.6,0.75) {\small $a$};
\node[fill=white, inner sep=1.5] at (0.6,-0.15) {\small $a$};

\draw[->,very thick]
	(1.8,0) -- ++(1,0);

\begin{scope}[xshift=4.5cm]

\draw
	(0,1.2) -- (-0.9,0) -- (0,-1.2) -- (1.12,-0.2);
	
\draw[yshift=1.2cm, rotate=38.5]
	(0,0) -- (1.2,-0.9) -- (0,-1.8);

\draw[dashed]
	(1.12,-0.2) -- (0,1.2) -- (0,-1.2);
	
\node at (-0.7,0) {\small $\epsilon$};
\node at (1.3,1.05) {\small $\delta$};

\node at (0,-1.4) {\small $A$};
\node at (1.25,-0.3) {\small $B$};
\node at (0,1.4) {\small $C$};
\node at (1.7,1.3) {\small $D$};
\node at (-1.1,0) {\small $E$};

\node[fill=white, inner sep=1.5] at (-0.45,0.6) {\small $a$};
\node[fill=white, inner sep=1.5] at (-0.45,-0.6) {\small $a$};
\node[fill=white, inner sep=1.5] at (0.75,1.2) {\small $a$};
\node[fill=white, inner sep=1.5] at (1.3,0.55) {\small $a$};
\node[fill=white, inner sep=1.5] at (0.55,-0.7) {\small $a$};

\node[fill=white, inner sep=1.5] at (0,0) {\small $u$};
\node[fill=white, inner sep=1.5] at (0.62,0.4) {\small $v$};

\node at (0.05,-0.9) {\small $\alpha$};
\node at (0.95,-0.1) {\small $\beta$};
\node at (0.1,0.95) {\small $\gamma$};

\end{scope}

\end{tikzpicture}
\caption{Construct the pentagon for earth map tiling.}
\label{emt_construct}
\end{figure}

We have 
\begin{align*}
\cos u 
&=\cos^2 a(1-\cos\epsilon)+\cos\epsilon, \\
\cos v 
&=\cos^2 a(1-\cos\delta)+\cos\delta.
\end{align*}
The equalities show that, for $a\in (0,\frac{1}{2}\pi)$, $u$ and $v$ are strictly increasing functions of $a$. Recall the following approximate values
\begin{align*}
f=16 &\colon 
\cos a=0.77943,\;
a=0.21550\pi; \\
f=20 &\colon 
\cos a=0.77680,\;
a=0.21683\pi;\\
f=24 &\colon 
\cos a=0.70688,\;
a=0.25010\pi;\\
f=24 &\colon 
\cos a=0.68125,\;
a=0.26143\pi.
\end{align*}
Therefore we take $[a_0,a_1]=[0.21\pi,0.22\pi]$ for $f=16,20$, and take $[a_0,a_1]=[0.2501\pi,0.2502\pi]$ for the $\cos a=0.70688$ case of $f=24$. Then we calculate the approximate values of $u(a_0)$ and $v(a_1)$
\begin{align*}
f=16 &\colon 
u(a_1)=0.4008\pi,\;
v(a_0)=0.2853\pi; \\
f=20 &\colon 
u(a_1)=0.4146\pi,\;
v(a_0)=0.2346\pi;\\
f=24 &\colon 
u(a_1)=0.4790\pi,\;
v(a_0)=0.2301\pi.\; (\cos a=0.70688)
\end{align*}
We verify that $u(a_1)\le a_0+v(a_0)$ holds in all three cases. Then we have
\[
u(a)\le u(a_1)\le a_0+v(a_0)\le a+v(a)\text{ for }a\in [a_0,a_1].
\]
Therefore we have $u-v\le a$ for $a\in [a_0,a_1]$, and we can construct the corresponding pentagons. 

Next, we calculate the pentagons for $a=a_0$ and $a=a_1$. This can be done by first using $\triangle ACE,\triangle BCD$ to calculate $u,v$, and then using $a,u,v$ to calculate $\triangle ABC$. Then we calculate the approximate values of the corresponding $\sigma=\alpha+\beta+\gamma$.
\begin{align*}
f=16 &\colon 
\sigma(a_0)=1.9861\pi,\;
\sigma(a_1)=2.0116\pi; \\
f=20 &\colon 
\sigma(a_0)=1.9868\pi,\;
\sigma(a_1)=2.0062\pi;\\
f=24 &\colon 
\sigma(a_0)=1.9999993\pi,\;
\sigma(a_1)=2.0000805\pi.\; (\cos a=0.70688)
\end{align*}
We find $\sigma(a_0)<2\pi<\sigma_1(a_1)$. By the intermediate value theorem, we conclude that $\sigma(a)=2\pi$ for some $a\in [a_0,a_1]$.

Since the equilateral pentagons constructed above can have the vertices $\alpha\beta\gamma,\delta^2\epsilon,\delta^n$, we may use the pentagons to construct earth map tilings. Then the earlier calculation for precise values is still valid, with ``compatible with $\cos a=0.77943$'' replaced by ``compatible with $a\in [a_0,a_1]$''. In particular, we get the precise values of the three triangles constituting the pentagon. Their approximate values are given in Figure \ref{emt_simple}. We conclude that the four angles of $\square ABCE$ are $<\pi$. Therefore $\square ABCE$ is convex, and $\square ABCE$, $\triangle BCD$ are in the separate hemispheres divided by the great circle $\bigcirc BC$. This implies the pentagon is simple.

\begin{figure}[htp]
\centering
\begin{tikzpicture}[>=latex,scale=1]

\foreach \a in {0,1,2}
{
\begin{scope}[xshift=4.5*\a cm]

\foreach \x in {1,...,5}
\draw[rotate=72*\x]
	(90:1.8) -- (18:1.8);
	
\draw[dashed]
	(90:1.8) -- (234:1.8)
	(90:1.8) -- (-54:1.8);

\node at (90:2) {\small $C$};
\node at (18:2) {\small $D$};
\node at (162:2) {\small $E$};
\node at (234:2) {\small $A$};
\node at (-54:2) {\small $B$};
		
\end{scope}
}



\node[rotate=90] at (0,1.1) {\footnotesize $0.21\pi$};
\node at (158:1.25) {\footnotesize $0.75\pi$};
\node at (22:1.25) {\footnotesize $0.50\pi$};
\node at (248:1.3) {\footnotesize $0.29\pi$};
\node at (-68:1.3) {\footnotesize $0.59\pi$};
\node[rotate=-70] at (1,-0.5) {\footnotesize $0.28\pi$};
\node[rotate=70] at (-1,-0.5) {\footnotesize $0.15\pi$};

\node at (0,-2.2) {$f=16$};


\begin{scope}[xshift=4.5cm]


\node[rotate=90] at (0,1.1) {\footnotesize $0.16\pi$};
\node at (158:1.25) {\footnotesize $0.80\pi$};
\node at (22:1.25) {\footnotesize $0.40\pi$};
\node at (248:1.3) {\footnotesize $0.18\pi$};
\node at (-68:1.3) {\footnotesize $0.72\pi$};
\node[rotate=-70] at (1,-0.5) {\footnotesize $0.33\pi$};
\node[rotate=70] at (-1,-0.5) {\footnotesize $0.12\pi$};

\node at (0,-2.2) {$f=20$};

\end{scope}


\begin{scope}[xshift=9cm]


\node[rotate=90] at (0,1.1) {\footnotesize $0.03\pi$};
\node at (158:1.25) {\footnotesize $0.83\pi$};
\node at (22:1.25) {\footnotesize $0.33\pi$};
\node at (248:1.3) {\footnotesize $0.02\pi$};
\node at (-68:1.3) {\footnotesize $0.95\pi$};
\node[rotate=-70] at (1,-0.5) {\footnotesize $0.37\pi$};
\node[rotate=70] at (-1,-0.5) {\footnotesize $0.11\pi$};

\node at (0,-2.2) {$f=24$ ($\cos a=0.70688$)};

\end{scope}

\end{tikzpicture}
\caption{Pentagon for the earth map tiling.}
\label{emt_simple}
\end{figure}

\begin{figure}[htp]
\centering
\begin{tikzpicture}[>=latex]

\pgftext{
	\includegraphics[scale=0.25]{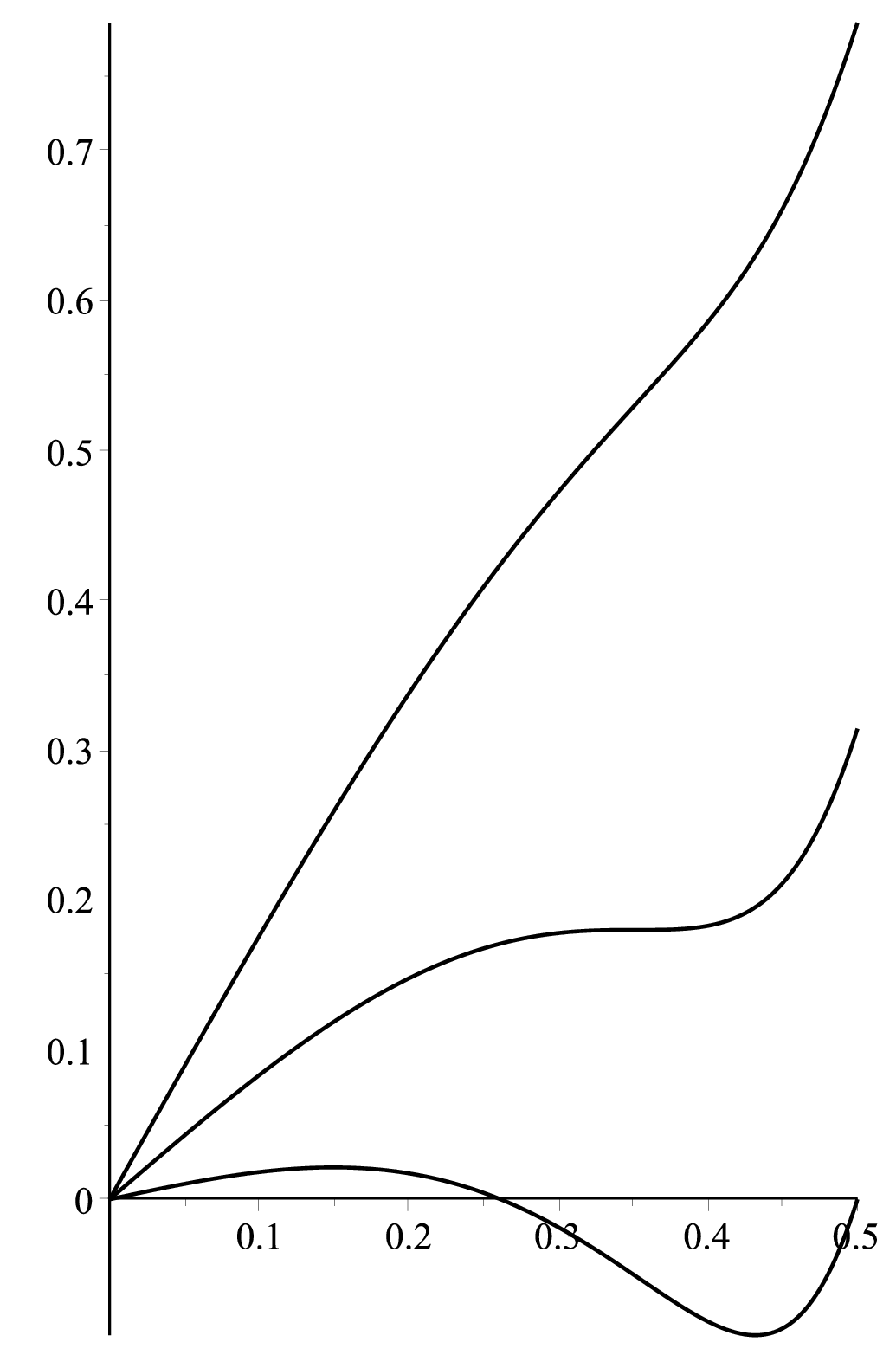}};

\node at (2.3,-2.4) {\small $\frac{a}{\pi}$};
\node at (-0.7,3) {\small $a+v-u$};
\node at (0.1,1.1) {\small $f=16$};
\node at (0.1,-1.1) {\small $f=20$};
\node at (0.1,-2.2) {\small $f=24$};

\begin{scope}[shift={(5cm, 0cm)}]

\pgftext{
	\includegraphics[scale=0.25]{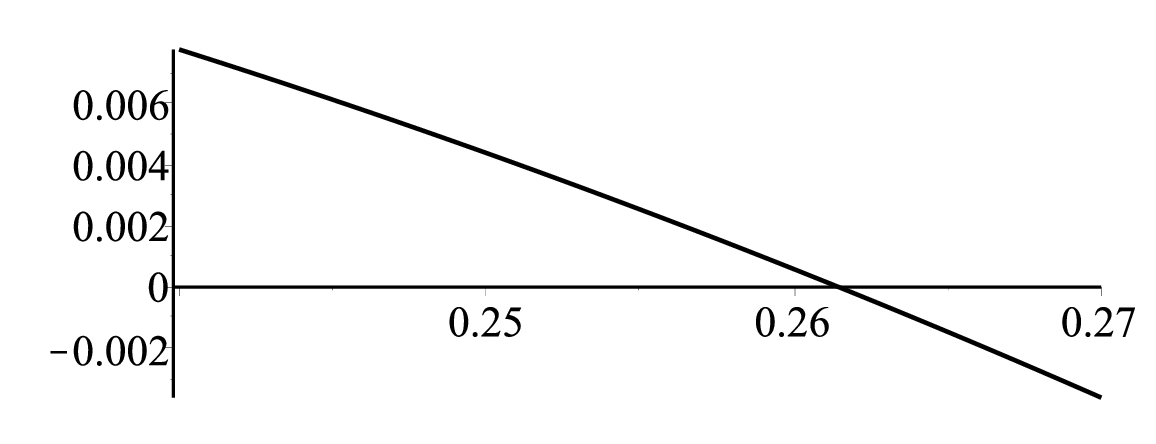}};

\end{scope}

\draw[gray]
	(0.05,-2.65) rectangle (0.5,-2.4);

\draw[->]
	(0.5,-2.5) -- (4,-0.7);

\draw[<-]
	(6.17,-0.28) -- ++(0,1);

\node at (6.6,0.85) {\tiny $a_{\max}=0.26143\pi$};

\end{tikzpicture}
\caption{When can we construct the equilateral pentagon?}
\label{emt_construct2}
\end{figure}

Finally, we examine the construction of the equilateral pentagon in more detail. Figure \ref{emt_construct2} is the graph of $a+v-u$ as a function of $a$. The graph shows that, for $f=16,20$, we have $u-v<a$ for all $a\in (0,\frac{1}{2}\pi)$. Therefore we may always construct the equilateral pentagon. For $f=24$, we may construct the equilateral pentagon only for $a$ in certain range $(0,a_{\max})$. Here we have $u-v=a$ at $a_{\max}=0.26143\pi$, i.e., the middle triangle $\triangle ABC$ is reduced to a line. This is the case discussed in Section \ref{case12e}, and the precise value of $a$ is given by \eqref{earth24B}.

If we draw the similar graph for $f=28,32,\dots$, then we find the curve to be below the $a$-axis. This means we cannot construct the equilateral pentagon suitable for tiling.

\subsection{Exceptional Case $\text{AVC}_3=\{\alpha\beta\gamma,\delta\epsilon^2\}$, $v_4=v_5=0$}
\label{casespecial}

By \cite[Lemma 1]{wy1} and $v_4=v_5=0$, there is a $3^5$-tile. The neighbourhood of a $3^5$-tile is given in Figure \ref{casespecialA}. Up to the symmetry of $\text{AVC}_3=\{\alpha\beta\gamma,\delta\epsilon^2\}$, we only need to consider $A1$ and $A3$ arrangements, as in $T_1$ of the two pictures.  

The first of Figure \ref{casespecialA} has the $A1$ arrangement. We have $\delta_1\cdots=\delta_1\epsilon_3\epsilon_4$. By $\epsilon_3$, we get $\epsilon_1\cdots=\delta_3\epsilon_1\epsilon_2$. Then $\epsilon_2$ contradicts $\alpha_1\cdots=\alpha\beta\gamma$.

\begin{figure}[htp]
\centering
\begin{tikzpicture}[>=latex,scale=1]

\foreach \a in {0,1}
{
\begin{scope}[xshift=3.5*\a cm]

\foreach \x in {1,...,5}
{
\draw[rotate=72*\x]
	(90:0.7) -- (18:0.7) -- (18:1.3) -- (54:1.7) -- (90:1.3);
	
\coordinate (X\x A\a) at (18+72*\x:0.45);
\coordinate (Y\x A\a) at (6+72*\x:0.8);
\coordinate (Z\x A\a) at (30+72*\x:0.8);
\coordinate (U\x A\a) at (10+72*\x:1.14);
\coordinate (V\x A\a) at (26+72*\x:1.14);
\coordinate (W\x A\a) at (-18+72*\x:1.45);
}

\node[draw,shape=circle, inner sep=0.5] at (0,0) {\small 1};
\foreach \x in {2,...,6}
\node[draw,shape=circle, inner sep=0.5] at (198-72*\x:1) {\small \x};

\end{scope}
}


\node at (X1A0) {\small $\alpha$};
\node at (X2A0) {\small $\beta$};
\node at (X3A0) {\small $\gamma$};
\node at (X4A0) {\small $\delta$};
\node at (X5A0) {\small $\epsilon$};

\node at (Y4A0) {\small $\epsilon$};
\node at (Y5A0) {\small $\delta$};

\node at (Z4A0) {\small $\epsilon$};
\node at (Z5A0) {\small $\epsilon$};





\node at (X1A1) {\small $\gamma$};
\node at (X2A1) {\small $\epsilon$};
\node at (X3A1) {\small $\alpha$};
\node at (X4A1) {\small $\beta$};
\node at (X5A1) {\small $\delta$};

\node at (Y1A1) {\small $\alpha$};
\node at (Y2A1) {\small $\delta$};
\node at (Y3A1) {\small $\gamma$};
\node at (Y4A1) {\small $\alpha$};
\node at (Y5A1) {\small $\epsilon$};

\node at (Z1A1) {\small $\beta$};
\node at (Z2A1) {\small $\epsilon$};
\node at (Z3A1) {\small $\beta$};
\node at (Z4A1) {\small $\gamma$};
\node at (Z5A1) {\small $\epsilon$};

\node at (U1A1) {\small $\beta$};
\node at (U2A1) {\small $\gamma$};
\node at (U3A1) {\small $\delta$};
\node at (U4A1) {\small $\epsilon$};
\node at (U5A1) {\small $\alpha$};

\node at (V1A1) {\small $\alpha$};
\node at (V2A1) {\small $\alpha$};
\node at (V3A1) {\small $\delta$};
\node at (V4A1) {\small $\delta$};
\node at (V5A1) {\small $\gamma$};

\node at (W1A1) {\small $\delta$};
\node at (W2A1) {\small $\epsilon$};
\node at (W3A1) {\small $\beta$};
\node at (W4A1) {\small $\gamma$};
\node at (W5A1) {\small $\beta$};

\end{tikzpicture}
\caption{$\text{AVC}_3=\{\alpha\beta\gamma,\delta\epsilon^2\}$: Neighbourhood tiling of a $3^5$-tile.}
\label{casespecialA}
\end{figure}

The second of Figure \ref{casespecialA} has the $A3$ arrangement. We have $\delta_1\cdots=\delta_1\epsilon_2\epsilon_3$. By $\epsilon_2$, we get $\gamma_1\cdots=\alpha_2\beta_6\gamma_1$. Then $\alpha_2,\epsilon_2$ determine $T_2$. By $\beta_6$, we get $\epsilon_1\cdots=\delta_6\epsilon_1\epsilon_5$. Then $\beta_6,\delta_6$ determine $T_6$. By $\epsilon_5$, we get $\alpha_1\cdots=\alpha_1\beta_4\gamma_5$. Then $\gamma_5,\epsilon_5$ determine $T_5$. By $\epsilon_3,\beta_4$, we get $\beta_1\cdots=\alpha_4\beta_1\gamma_3$. Then $\gamma_3,\epsilon_3,\alpha_4,\beta_4$ determine $T_3,T_4$. The neighbourhood tiling is the same as the first of Figure \ref{case55A}.

Next we argue $f\le 24$. Since $f$ is even, it is sufficient to show that $f<26$. The angle sums of $\alpha\beta\gamma,\delta\epsilon^2$ imply
\[
\tfrac{1}{2}\delta
=\alpha+\beta+\gamma+\delta+\epsilon-3\pi
=\tfrac{4}{f}\pi,\quad
\delta=\tfrac{8}{f}\pi,\;
\epsilon=\pi-\tfrac{1}{2}\delta.
\]
By $f\ge 16$, we have $\delta\le\frac{1}{2}\pi$. We will have two inequality restrictions on $f$. 

\begin{figure}[htp]
\centering
\begin{tikzpicture}[>=latex]

\pgftext{
	\includegraphics[scale=0.4]{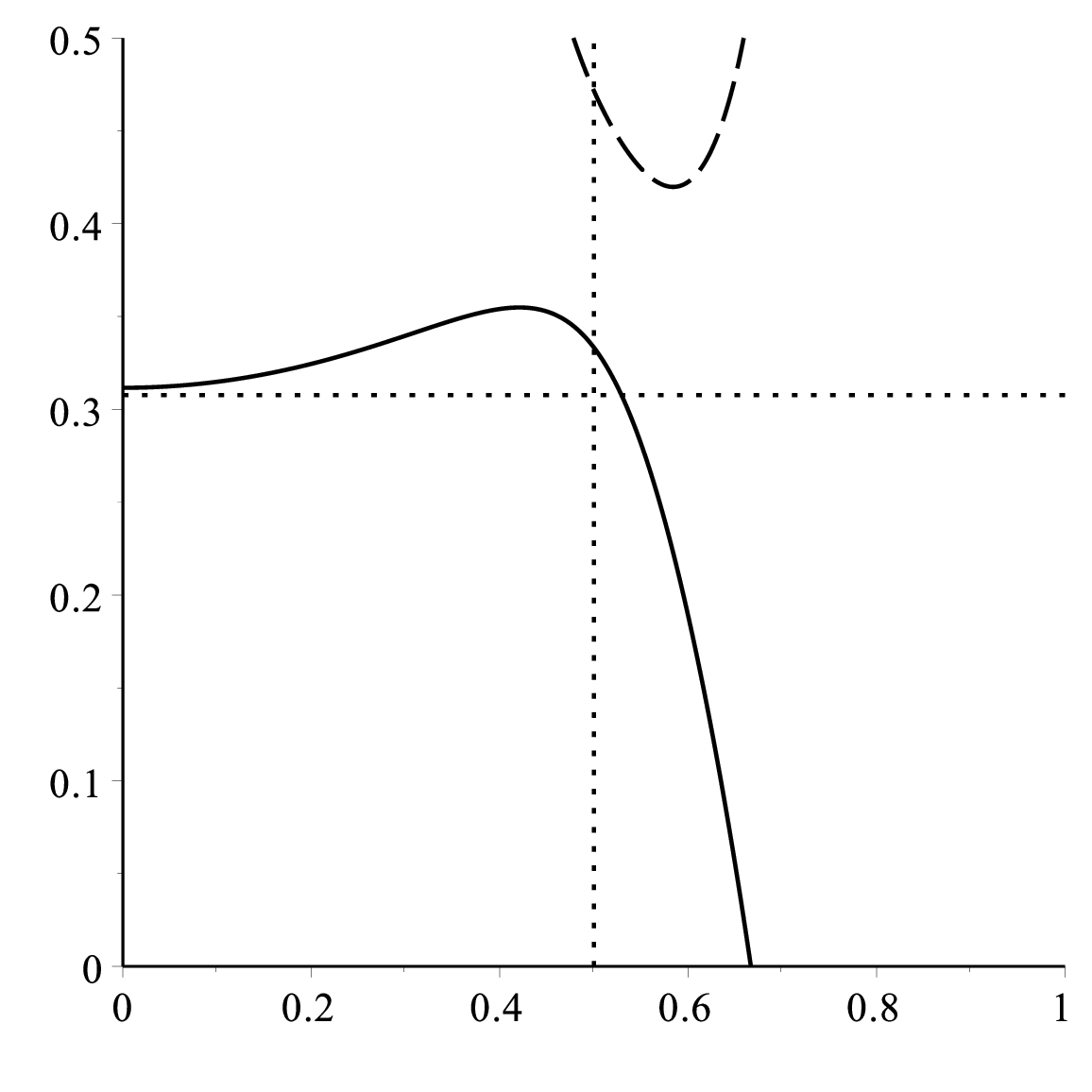}};

\node at (0.9,3.8) {\small $\Sigma\ge \frac{22}{13}\pi$};
\node at (1.2,2) {\small $\Sigma<\frac{22}{13}\pi$};
\node at (-1,2.4) {\small $u-v\le a$};
\node at (-1,0) {\small $u-v>a$};
\node at (-2.9,3.8) {\small $\delta$};
\node at (4,-3) {\small $\frac{a}{\pi}$};

\draw[<->]
	(3,-3) -- node[fill=white] {\small $f\ge 26$} (3,1.1);

\draw[<-]
	(3,1.3) -- node[fill=white] {\small $f\le 24$} (3,4);
	
\end{tikzpicture}
\caption{$f\le 24$ for the exceptional case.}
\label{casespecialB}
\end{figure}

The $A3$ arrangement means the pentagon is $\pentagon ABDEC$, given by the first of Figure \ref{emt_construct}. We recall $a<\pi$ by the simple pentagon requirement. Let $u$ and $v$ be the arcs of length $<\pi$ connecting $A,C$ and connecting  $B,C$. We may determine the arcs $u$ and $v$ by the cosine laws
\begin{align*}
\cos u &=\cos^2a+\sin^2a\cos\epsilon
=\cos^2a-\sin^2a\cos\tfrac{1}{2}\delta, \\
\cos v &=\cos^2a+\sin^2a\cos\delta.
\end{align*}
The inequality $u-v\le a$ then defines a region on the rectangle $(a,\delta)\in (0,\pi)\times (0,\frac{1}{2}\pi]$. As shown by Figure \ref{casespecialB}, for $0<a\le\frac{1}{2}\pi$, this already implies $f<26$.

\begin{figure}[htp]
\centering
\begin{tikzpicture}[>=latex,scale=1]


\begin{scope}[xshift=-3cm]

\foreach \a in {0,...,4}
\draw[rotate=72*\a]
	(18:1.2) -- (-54:1.2);

\draw[dashed]
	(234:1.2) -- (90:1.2) -- (-54:1.2);

\node at (90:0.9) {$\gamma$};
\node at (162:0.9) {$\epsilon$};
\node at (15:0.9) {$\delta$};
\node at (234:0.9) {$\alpha$};
\node at (-54:0.9) {$\beta$};

\node at (90:1.4) {$C$};
\node at (162:1.4) {$E$};
\node at (15:1.4) {$D$};
\node at (234:1.4) {$A$};
\node at (-54:1.4) {$B$};

\node[fill=white, inner sep=2] at (18:0.4) {\small $v$};
\node[fill=white, inner sep=2] at (162:0.4) {\small $u$};

\end{scope}


\draw[dashed]
	(0,0) arc (135:45:4) 
	(0,0) arc (-135:-45:4)
	(0,0) arc (-135:-155:4)
	(0,0) arc (135:75:4) arc (10:-10:6);
	
\draw
	(0,0) arc (135:75:4) arc (-30:-58:6)
	(0,0) arc (-135:-75:4)
	;

\node at (-0.2,-0.1) {$D$};
\node at (0.3,0) {$\delta$};

\node at (3.9,1.3) {$B$};
\node at (3.9,-1.3) {$C$};
\node[fill=white, inner sep=2] at (3.95,0) {$v$};

\node at (1.8,1.2) {$a$};
\node at (1.75,-1.3) {$X$};
\node at (1.8,-0.8) {$\xi$};

\node at (-0.8,1.3) {$P$};

\end{tikzpicture}
\caption{Pentagon in $A3$ arrangement, and $\triangle BCD$ lies outside $\square ABCE$.}
\label{casespecialE}
\end{figure}

For $\frac{1}{2}\pi<a<\pi$, another inequality may be obtained by estimating the area of the pentagon. First we argue that $\delta\le\frac{1}{2}\pi$ implies that $\triangle BCD$ lies outside $\square ABCE$. In the second of Figure \ref{casespecialE}, we draw the isosceles triangle $\triangle BCD$, where $BC$ has length $<\pi$. For any $X$ on $CD$, we connect an arc $BX$ inside the isosceles triangle, and consider the change of angle $\xi=\angle BXP$. By $\frac{1}{2}\pi<a<\pi$, we have $\angle BCP=\angle CBD>\frac{1}{2}\pi$. We also know the outside angle $\angle BDP=\pi-\delta\ge \frac{1}{2}\pi$. This means that $\xi>\frac{1}{2}\pi$ when $X=C$ and $\xi\ge \frac{1}{2}\pi$ when $X=D$. It is then a fact of the spherical geometry that $\xi>\frac{1}{2}\pi$ for any $X$ on the interior of $CD$. (We can see this, for example, by using stereographic projection.) Then by $\xi>\frac{1}{2}\pi\ge \delta$, we conclude $a>BX$. Since the pentagon is simple, this implies that $BX$ cannot be part of the edge $BA$. In other words, the edge $BA$ must point outside the triangle $\triangle BCD$. By the same reason, the edge $CE$ must point outside the triangle $\triangle BCD$. This proves that $\triangle BCD$ lies outside $\square ABCE$.

Since $\triangle BCD$ lies outside $\square ABCE$, we have
\[
\tfrac{4}{f}\pi
=\text{Area}(\pentagon ABDCE)
\ge \text{Area}(\square ABCE).
\]
The area of the quadrilateral can be further estimated
\[
\text{Area}(\square ABCE)
\ge \text{Area}(\triangle ACE) - \text{Area}(\triangle ABC).
\]
By the assumption $\frac{1}{2}\pi<a<\pi$, we have
\[
\text{Area}(\triangle ACE)
\ge \epsilon=\pi-\tfrac{1}{2}\delta=(1-\tfrac{4}{f})\pi.
\]
Moreover, $\text{Area}(\triangle ABC)+\pi$ is the sum $\sum$ of the three angles of $\triangle ABC$. Combining all the inequalities together, we get
\[
{\textstyle \sum}\ge 2(1-\tfrac{4}{f})\pi.
\]
The sides of $\triangle ABC$ are $u,v,a$, and its three angles can be calculated by the cosine law. Then $\sum$ may be explicitly expressed as a function of $a$ and $\delta$.

If $f\ge 26$, then $\sum\ge \frac{22}{13}\pi$. The region $\sum\ge \frac{22}{13}\pi$ is indicated by dashed boundary in Figure \ref{casespecialB}, and the picture shows that $\sum\ge \frac{22}{13}\pi$ implies $f\le 24$. Therefore we conclude $f\le 24$ for the exceptional case.

Recall that we obtained the neighbourhood tiling in the second of Figure \ref{casespecialA} for $\text{AVC}_3=\{\alpha\beta\gamma,\delta\epsilon^2\}$. The neighbourhood tiling is the same as the first of Figure \ref{case55A} for Case 5.5. Similar to Section \ref{case55}, we may use the vertex counting equation and \cite[Theorem 6]{yan} to show that $v_4=v_5=0$ and $16\le f\le 24$ imply that the tiling is the earth map tiling with exactly two degree $6$ vertices. Then by the propagation argument in Section \ref{case55}, we get the earth map tiling in Figure \ref{earth}. The calculation in Section \ref{case15b} shows that there are two possible pentagons suitable for the earth map tiling.


\begin{thebibliography}{1}

\bibitem{ay1}
Y.~Akama, M.~Yan.
\newblock On deformed dodecahedron tiling.
\newblock {\em preprint}, arXiv:1403.6907, 2014.

\bibitem{gsy}
H.~H.~Gao, N.~Shi, M.~Yan.
\newblock Spherical tiling by $12$ congruent pentagons.
\newblock {\em J. Combinatorial Theory Ser. A}, 120(4):744--776, 2013.

\bibitem{luk}
H.~P.~Luk.
\newblock Angles in Spherical Pentagon Tilings.
\newblock MPhil Thesis, Hong Kong University of Science and Technology, 2012.

\bibitem{ly1}
H.~P.~Luk, M.~Yan.
\newblock Tilings of the sphere by congruent almost equilateral pentagons I: five distinct angles.
\newblock {\em preprint}, 2021.

\bibitem{wy1}
E.~X.~Wang, M.~Yan.
\newblock Tilings of the sphere by congruent pentagons I: edge combinations $a^2b^2c$ and $a^3bc$.
\newblock {\em preprint}, arXiv:1804.03770, 2021.

\bibitem{wy2}
E.~X.~Wang, M.~Yan.
\newblock Tilings of the sphere by congruent pentagons II: edge combination $a^3b^2$.
\newblock {\em preprint}, arXiv:1903.02712, 2021.

\bibitem{wy3}
E.~X.~Wang, M.~Yan.
\newblock Moduli of pentagonal subdivision tiling.
\newblock preprint, arXiv: 1907.08776, 2019.

\bibitem{yan}
M.~Yan.
\newblock Combinatorial tilings of the sphere by pentagons.
\newblock {\em Elec. J. of Combi.}, 20(1):\#P54, 2013.

\end{thebibliography}
\end{document}